\documentclass[11pt]{amsart}

\usepackage{amsfonts,amssymb,amsmath,mathrsfs,graphicx}
\usepackage{float}
\usepackage{epsfig}

\usepackage{CJK,color}

\definecolor{red}{rgb}{1.00,0.00,0.00}
\definecolor{blue}{rgb}{0.00,0.00,0.63}
\definecolor{black}{rgb}{0.00,0.00,0.00}

\newtheorem{theorem}{Theorem}[section]
\newtheorem{lemma}{Lemma}[section]

\newtheorem{proposition}{Proposition}[section]
\newtheorem{remark}{Remark}[section]

\newcommand{\bbr}{\mathbb R}

\newcommand{\bbt}{\mathbb T}

\newcommand{\ba}{\begin{aligned}}
\newcommand{\ea}{\end{aligned}}
\newcommand{\be}{\begin{equation}}
\newcommand{\ee}{\end{equation}}
\renewcommand{\div}{ {\rm div }  }
\def\na{\nabla}

\def\charf {\mbox{{\text 1}\kern-.30em {\text l}}}

\def\mb{\mathbf}

\def\f{\frac}
\def\t{\theta}
\def\d{\delta}
\def\a{\alpha}

\def\z{\zeta}

\def\di{\displaystyle}
\def\r{\rho}
\def\wt{\widetilde}
\def\k{\kappa}

\def\v{\varepsilon}

\def\charf {\mbox{{\text 1}\kern-.30em {\text l}}}

\def\mb{\mathbf}
\def\f{\frac}
\def\t{\theta}

\def\d{\delta}
\def\s{\sigma}
\def\a{\alpha}

\def\di{\displaystyle}
\def\r{\rho}
\def\wt{\widetilde}
\def\k{\kappa}
\def\pa{\partial}

\begin{document}

\title[Stability of planar rarefaction wave to 3D Boltzmann equation]{Nonlinear stability of planar rarefaction wave to the three-dimensional Boltzmann equation}

\author[Wang]{Teng Wang}
\address[Teng Wang]{\newline  College of Applied Sciences, Beijing University of Technology, Beijing 100124, P. R. China}
\email{tengwang@amss.ac.cn}

\author[Wang]{Yi Wang}
\address[Yi Wang]{\newline  CEMS, HCMS, NCMIS, Academy of Mathematics and Systems Science, Chinese Academy of Sciences, Beijing 100190, China
\newline and School of Mathematical Sciences, University of Chinese Academy of Sciences, Beijing 100049, China}
\email{wangyi@amss.ac.cn}

\date{\today}



\maketitle

\begin{abstract}
We investigate the time-asymptotic stability of planar rarefaction wave for the three-dimensional Boltzmann equation, based on the micro-macro decomposition introduced in \cite{LY, LYY} and our new observations on the underlying wave structures of the equation to overcome the difficulties due to the wave propagation along the transverse directions and its interactions with the planar rarefaction wave.  Note that this is the first stability result of planar rarefaction wave for 3D Boltzmann equation, while the corresponding results for the shock and contact discontinuities are still completely open.\end{abstract}

\section{Introduction and main result}
\setcounter{equation}{0}
We  investigate the time-asymptotic stability of planar rarefaction wave for the three-dimensional Boltzmann equation which  takes the form
\begin{equation}\label{B}
f_{t} +v\cdot\nabla_x f=\frac1\kappa Q(f, f),
\end{equation}
with the initial data
\begin{equation}\label{initial}
f(0,x,v)=f_0(x,v),
\end{equation}
where the time variable $t\in \mathbb {R}^+$, spatial variables $x=(x_1,x_2,x_3)\in {\mathbb R}\times {\mathbb T}^2:=\mathbb{D}$ with ${\mathbb  T}^2:=(\mathbb {R}/\mathbb {Z})^2$ being the two-dimensional unit flat torus, the particle velocity $v=(v_1,v_2,v_3)\in {\mathbb R}^3$  and $f=f(t,x,v)$ represents the distributional
density of particles at time-space $(t,x)$ with velocity $v$.  In the Boltzmann equation \eqref{B}, the physical parameter $\kappa$ called the Knudsen
number is proportional to the mean free path of the interacting particles. Since we are concerned with the time-asymptotic behavior of the solution to the Boltzmann equation \eqref{B}, the Knudsen number $\kappa$ will be fixed to be 1 in the following.  For the hard sphere model, the collision operator $Q(f,f)$ takes the following bilinear and symmetric form
$$
\begin{array}{l}
 \di Q(f,g)(v) =\f{1}{2}
\int_{{\mathbb R}^3}\!\!\int_{{\mb S}_+^2} \big[f(v_*')g(v')+f(v')g(v_*')\\
\di ~\qquad\qquad\qquad\qquad\qquad -f(v_*) g(v)-f(v) g(v_*) \big]
|(v-v_*)\cdot\Omega|
 \; d v_* d\Omega,
 \end{array}
$$
where the unit vector $\Omega \in {\mb S}^2_+=\{\Omega\in {\mb S}^2:\ (v-v_*)\cdot \Omega\geq 0\}$, and $(v,v_*)$ and $(v',v_*')$ are the pair velocities of the two particles before and after a binary elastic collision respectively, which together with the conservation laws of momentum and energy during the collisions, imply the following relations
$$
 v'= v -[(v-v_*)\cdot \Omega] \; \Omega, \qquad
 v_*'= v_* + [(v-v_*)\cdot \Omega] \; \Omega, \quad \Omega \in {\mb S}^2.
$$

It is well-known that Boltzmann equation is closely related to the system of classical fluid mechanics, in particular, the systems of compressible Euler and Navier-Stokes equations, as described by the famous Hilbert expansion and Chapman-Enskog expansion, respectively. Either the Hilbert expansion or the Chapman-Enskog expansion yields the compressible Euler equations in the leading order with respect to the Knudsen number $\kappa$. The system of compressible Euler equations is a typical example of hyperbolic conservation laws system, which has the distinguished feature that the solution may blow up in finite-time, i. e., the formation of the shock wave, no matter how smooth or small the initial values are. In fact, there are three basic wave patterns to the system of hyperbolic conservation laws, that is, two nonlinear waves, shock and rarefaction waves, in the genuinely nonlinear characteristic fields and a linear wave, contact discontinuity, in the linearly degenerate field. These dilation invariant solutions, and their linear superposition in the increasing order of characteristic speed, called Riemann solutions, govern both the local and large time asymptotic behavior of general solutions to the inviscid Euler system. Hence, one can expect the wave phenomena to Boltzmann equation as for the macroscopic fluid dynamics and it is interesting and important both mathematically and physically to prove the dynamic stability of these basic wave patterns for the Boltzmann equation.

For the Boltzmann equation with slab symmetry, i. e., the spatial one-dimensional case, the stability of the these three basic wave patterns are well-understood until now. For instance, the pioneering study on the stability of viscous shock wave was first proved by Liu and Yu \cite{LY} with the zero total macroscopic mass condition by the energy method based on the micro-macro decomposition while the existence of viscous shock profile to the Boltzmann equations is given by Caflish and Nicolaenko \cite{Caflish-Nicolaenko} and Liu and Yu \cite{Liu-Yu-1}. Then stability of rarefaction wave fan is proved by  Liu, Yang, Yu, and Zhao \cite{Liu-Yang-Yu-Zhao} and the stability of viscous contact wave, which is the viscous version of contact discontinuity, by Huang and Yang \cite{Huang-Yang} with the zero mass condition and Huang, Xin, and Yang \cite{Huang-Xin-Yang} without the zero mass condition.  Furthermore, Yu \cite{Yu} proved the stability of single viscous shock profile without zero mass condition by the elegant point-wise method based on the Green function around the shock profile. Recently, Wang and Wang \cite{Wang-Wang} proved the stability of the superposition of two viscous shock profiles to the Boltzmann equation without the zero mass condition by the weighted characteristic energy method.
Nevertheless, the hydrodynamic limit of Boltzmann equation to Euler system in the setting of  1D Riemann solution is proved in \cite{Yu1, Huang-Wang-Yang, xin-zeng, Huang-Wang-Yang-2, Huang-Wang-Yang-3} and finally in \cite{hwwy} for 1D generic Riemann solution.

However, for the three-dimensional Boltzmann equation \eqref{B}, there is no any result on the stability of  these three basic wave patterns, as far as we know. The main purpose of the present paper is to establish the first result on the time-asymptotic stability of planar rarefaction wave to 3D Boltzmann equation in the infinite long flat nozzle domain $\mathbb{D}:={\mathbb R}\times {\mathbb T}^2$. In our stability proof,  we first use the macro-micro decomposition invented by Liu and Yu \cite{LY} and  Liu, Yang, and Yu \cite{LYY} to rewrite Boltzmann equation \eqref{B} by a fluid-type system, which is a compressible Navier-Stokes equations with temperature-dependent viscosities and heat-conductivities coupled with the microscopic terms, and a non-fluid part equation where the linearized collision operator around the local Maxwellian are strongly dissipative.  Then we construct a smooth profile to the 1D self-similar rarefaction wave fan and we look for the solution to the 3D Boltzmann equation around this profile and prove the stability of this planar rarefaction wave.

Compared with the stability problem in one-dimensional case, the main difficulties in the proof of the stability of planar rarefaction wave lie in the wave propagations along the transverse directions $(x_2, x_3)\in \mathbb{T}^2$ and their interactions with the planar rarefaction wave in the direction $x_1\in \mathbb{R}$. Motivated by the previous work of the second author and his collaborator in \cite{LW} for the stability of planar rarefaction wave of the two-dimensional compressible Navier-Stokes equations,  we use some underlying wave structures to overcome the difficulties mentioned above. However, different from \cite{LW}, here we need to consider the case that both viscosities and heat conductivity depend on the temperature and moreover, we need to cope with the microscopic terms and their interactions with the fluid part. Furthermore, we need to make full use of the advantage of the different forms of the fluid system in the micro-macro decomposition for Boltzmann equation when we estimate the fluid parts. For the stability of planar rarefaction wave to the multi-dimensional viscous conservation laws, one can also refer \cite{X, LWW} and the references therein.

Now we formulate our problem. Since we are concerned with the time-asymptotic stability of planar rarefaction wave to 3D Boltzmann equation \eqref{B}, it is assumed that the far fields conditions of initial data on the $x_1-$direction
\begin{equation}\label{far-field}
 f(0,x,v)=f_0(x,v)\rightarrow \mathbf{M}_{[\rho_\pm,u_\pm,\t_\pm]}(v)
 :=\frac{\rho_\pm}{\sqrt{(2\pi R\t_\pm)^3}}\exp\left(-\frac{|v-u_\pm|^2}{2R\t_\pm}\right),
\\
\end{equation}
as $x_1\rightarrow\pm\infty,$ with $\rho_\pm>0$, $u_\pm=(u_{1\pm},0,0)^t$, $\t_\pm>0$ being
prescribed constant states.  Moreover, the periodic boundary conditions are imposed on $(x_2,x_3)\in\mathbb{T}^2$ for the solution $f(t,x,v)$. Here the two  end states $(\rho_\pm,u_\pm,\t_\pm)$ are connected by the rarefaction wave solution to the Riemann problem of the
corresponding 1D compressible Euler system
\begin{equation}\label{euler}
\left\{
\begin{array}{l}
\di\rho_t+(\rho u_1)_{x_1}=0,   \quad\qquad \qquad x_1 \in \mathbb{R}, ~t > 0, \\
\di(\rho u_1)_t+(\rho u_1^2+p)_{x_1}=0,\\
\di(\rho u_i)_t+(\rho u_1u_i)_{x_1}=0,~i=2,3,\\
\di[\rho(e+\f{|u|^2}{2})]_t+[\rho u_1(e+\f{|u|^2}{2})+pu_1]_{x_1}=0,
\end{array}
\right.
\end{equation}
with the Riemann initial data
\begin{equation}\label{R-in}
(\r,u,\t)(0,x_1)=(\r_0^r,u_0^r,\t_0^r)(x_1)=\left\{
\begin{array}{ll}
\di (\r_-,u_-,\t_-), &\di x_1<0 \\
\di (\r_+,u_+,\t_+), &\di x_1>0.
\end{array}
 \right.
\end{equation}

It could be expected that the large-time behavior of the solution to the 3D Boltzmann equation \eqref{B}-\eqref{far-field} is closely related to the Riemann problem to the corresponding  3D compressible Euler equations
\begin{equation}  \label{2ES}
\begin{cases}
\displaystyle \rho_t + \div (\rho u)= 0,               \quad\qquad \qquad x \in \bbr\times\bbt^2, ~t > 0, \\
\displaystyle (\rho u)_t + \div(\rho u\otimes u) + \nabla p = 0,\\
\di (\rho (e+\f{|u|^2}{2}))_t + \div (\rho u(e+\f{|u|^2}{2}) + pu)= 0,
\end{cases}
\end{equation}
with the Riemann initial data
\begin{equation} \label{2ini}
(\rho_0, u_0,\theta_0)(x) = \begin{cases}
(\rho_-, u_-,\theta_-),     \quad x_1 < 0,\\
(\rho_+, u_+,\theta_+),     \quad x_1 > 0.
\end{cases}
\end{equation}

There are essential differences between the one-dimensional Riemann problem \eqref{euler}-\eqref{R-in} and the multi-dimensional Riemann problem \eqref{2ES}-\eqref{2ini}. In the 2D isentropic regime, that is, the system \eqref{2ES} with the constant entropy and then the energy equation $\eqref{2ES}_3$ being satisfied trivially,  it is first proved by Chiodaroli, DeLellis, and Kreml \cite{CDK}  and Chiodaroli and Kreml \cite{CK}
that there are infinitely many bounded admissible weak solutions to \eqref{2ES}-\eqref{2ini} satisfying the natural entropy condition for shock Riemann initial data by using the convex integration methods in DeLellis and Sz$\acute{\rm{e}}$kelyhid \cite{DS} while the construction of infinitely many admissible weak solutions in \cite{CDK,CK} seems essential to the multi-dimensional system and could not be applied to one-dimensional problem \eqref{euler}-\eqref{R-in}. Then Klingenberg and Markfelder \cite{KM} and Brezina, Chiodaroli, and Kreml \cite{BCK} extend the results in \cite{CDK,CK} to the case when the corresponding Riemann initial data contain shock or contact discontinuity. On the other hand, similar to the one-dimensional case, for the Riemann solution only containing rarefaction waves to  \eqref{2ES}-\eqref{2ini}, Chen and Chen\cite{CC} and Feireisl and Kreml \cite{FK},  Feireisl, Kreml, and Vasseur \cite{FKV} independently proved the uniqueness of the uniformly bounded admissible weak solution even the rarefaction waves are connected with vacuum states (cf. \cite{CC}).

As mentioned before, we can expect that the large-time behavior of the solution to the Boltzmann equation \eqref{B}-\eqref{far-field} is determined by the Riemann problem to the corresponding inviscid Euler system \eqref{2ES} or \eqref{euler}, which contains planar shock wave and rarefaction wave in the genuinely nonlinear characteristic fields and contact discontinuity in the linearly degenerate field.  Our goal in the paper is to prove the above expectations in mathematics rigor and to investigate the dynamic stability of planar rarefaction wave for the 3D Boltzmann equation \eqref{B} at a first step.

Now we first carry out the micro-macro decomposition around the local\break Maxwellian to the Boltzmann equation \eqref{B} as introduced by Liu and Yu \cite{LY} and Liu, Yang, and Yu \cite{LYY}.
In fact, for any solution $f(t,x,v)$ to equation $\eqref{B}$, there are five macroscopic (fluid) quantities: the mass density $\r(t,x)$, the momentum $m(t,x)=\r u(t,x)$,
and the total energy $E(t,x)=\r\big(e+\frac12|u|^2\big)(t,x)$ defined by
\begin{equation}
\left\{
\begin{array}{l}
\di\rho(t,x)=\int_{\mathbb{R}^3}\xi_0(v)f(t,x,v)dv,\\
\di\rho u_i(t,x)=\int_{\mathbb{R}^3}\xi_i(v)f(t,x,v)dv,~i=1,2,3,\\
\di\rho(e+\f{|u|^2}{2})(t,x)=\int_{\mathbb{R}^3}\xi_4(v)f(t,x,v)dv,
\end{array}
\right. \label{macro}
\end{equation}
where  $\xi_i(v)$ $(i=0,1,2,3,4)$ are the collision
invariants given by
\begin{equation}
 \xi_0(v) = 1,~~~
 \xi_i(v) = v_i~  (i=1,2,3),~~~
 \xi_4(v) = \f{1}{2} |v|^2,
\label{collision-invar}
\end{equation}
and satisfy
$$
\int_{{\mathbb R}^3} \xi_i(v) Q(g_1,g_2) d v =0,\quad {\textrm
{for} } \ \  i=0,1,2,3,4.
$$

For a solution $f(t,x,v)$ to the Boltzmann equation \eqref{B}, we decompose it into the macroscopic (fluid) component, i.e., the local Maxwellian $\mb{M}=\mb{M}_{[\rho, u, \theta]}(t,x,v)$,
and the microscopic (non-fluid) component, i.e., $\mb{G}=\mb{G}(t,x,v)$ as follows (cf. \cite{LYY})
$$
f(t,x,v)=\mb{M}(t,x,v)+\mb{G}(t,x,v).
$$
Here, the local Maxwellian $\mb{M}$ is associated to the solution $f(t,x,v)$ of the equation $\eqref{B}$  in terms of the five fluid quantities $(\rho, u, \theta)(t,x)$ defined by
\begin{equation}
\mb{M}:=\mb{M}_{[\rho, u,\t]} (t,x,v) = \f{\rho(t,x)}{\sqrt{ (2 \pi
R \t(t,x))^3}} e^{-\f{|v-u(t,x)|^2}{2R\t(t,x)}},  \label{maxwellian}
\end{equation}
where $\t(t,x)$ is the temperature which is related to the internal
energy $e(t,x)$ by $e=\frac{3}{2}R\t$ with $R>0$ being the gas constant,  and $u(t,x)=\big(u_1(t,x),u_2(t,x),u_3(t,x)\big)^t$
is the fluid velocity.

For some given global or local Maxwellian $\widetilde{\mb{M}}$, the weighted inner product in $L^2_v$ space with respect to the Maxwellian $\wt{\mb{M}}$ is defined by:
\begin{equation}\label{product}
 \langle g_1,g_2\rangle_{\widetilde{\mb{M}}}\equiv \int_{{\mathbb R}^3}
 \f{1}{\widetilde{\mb{M}}}g_1(v)g_2(v)d v,
\end{equation}
for any functions $g_1(v),~g_2(v)$ such that the above integral is well-defined, and $\|g_1\|_{L^2_{v}(\frac{1}{\sqrt{\wt{\mb{M}}}})}:=\langle g_1,g_1\rangle^{\frac12}$.
If $\widetilde{\mb{M}}$ is the local Maxwellian $\mb{M}$ in \eqref{maxwellian}, we shall use the simplified notation $\langle\cdot,\cdot\rangle$
instead of $\langle\cdot,\cdot\rangle_{\mb{M}}$ if without confusions.  With respect to this inner product $\langle\cdot,\cdot\rangle$,
the following five pairwise orthogonal bases span the macroscopic space  $\mathfrak{N}$
\begin{equation}\label{orthogonal-base}
\left\{
\begin{array}{l}
 \chi_0(v) \equiv {\di\f1{\sqrt{\rho}}\mb{M}}, \quad
 \chi_i(v) \equiv {\di\f{v_i-u_i}{\sqrt{R\t\rho}}\mb{M}} \ \ {\textrm {for} }\ \  i=1,2,3, \\[2mm]
 \chi_4(v) \equiv
 {\di\f{1}{\sqrt{6\rho}}(\f{|v-u|^2}{R\t}-3)\mb{M}},\quad \langle\chi_i,\chi_j\rangle=\delta_{ij}, ~i,j=0,1,2,3,4.
 \end{array}
\right.
\end{equation}
In terms of above orthogonal bases, the macroscopic projection $\mb{P}_0$ from $L^2(\mathbb{R}^3_v)$ to
$\mathfrak{N}$ and the microscopic projection $\mb{P}_1$  from $L^2({\mathbb R}^3_v)$ to $\mathfrak{N}^\bot$ can be defined as
\begin{equation*}
 \mb{P}_0g = {\di\sum_{j=0}^4\langle g,\chi_j\rangle\chi_j},\qquad
 \mb{P}_1g= g-\mb{P}_0g.
\end{equation*}

Based on the above preparation, the solution $f(t,x,v)$  to the  equation $\eqref{B}$ can be decomposed into
the macroscopic (fluid) part, i.e., the local Maxwellian $\mb{M}=\mb{M}(t,x,v)$
defined in \eqref{maxwellian}, and the microscopic (non-fluid) part, i.e.
$\mb{G}=\mb{G}(t,x,v)$:
$$
f(t,x,v)=\mb{M}(t,x,v)+\mb{G}(t,x,v),\quad
\mb{P}_0f=\mb{M},~~~\mb{P}_1f=\mb{G},
$$
and the  equation $\eqref{B}$ becomes
\begin{equation}
(\mb{M}+\mb{G})_t+v\cdot\nabla_x(\mb{M}+\mb{G})=2Q(\mb{M},\mb{G})+Q(\mb{G},\mb{G}). \label{M+G}
\end{equation}
Taking the inner product of the equation \eqref{M+G} and the collision invariants $\xi_i(v)$ $(i=0,1,2,3,4)$  with respect to
$v$ over ${\mathbb R}^3$, one has the following  system for the fluid variables $(\rho, u, \theta)$:
\begin{equation}\label{ns-1}
\left\{
\begin{array}{lll}
\di \rho_{t}+{\rm div}_x (\rho u)=0, \\
\di (\rho u)_t+{\rm div}_x (\rho u\otimes u) +\nabla_x  p =-\int v\otimes v\cdot\nabla_x\mb{G}dv,  \\
\di [\rho(e+\f{|u|^2}{2})]_t+{\rm div}_x  [\rho u(e+\f{|u|^2}{2})+pu]=-\int\f12|v|^2v\cdot\nabla_x\mb{G}dv,
\end{array}
\right.
\end{equation}
where $
p=\frac23\rho e=R\rho\t
$ is the pressure for the mono-atomic gas. However, the above fluid-type system~\eqref{ns-1} is not self-contained and a equation for the microscopic component ${\mb{G}}$
is needed, which can be  derived by applying the projection operator $\mb{P}_1$ into   the equation\eqref{M+G}:
\begin{equation}\label{non-fluid}
\mb{G}_t+\mb{P}_1(v\cdot\nabla_x\mb{M})+\mb{P}_1(v\cdot\nabla_x\mb{G})=\mb{L}_\mb{M}\mb{G}+ Q(\mb{G}, \mb{G}),
\end{equation}
where $\mb{L}_\mb{M}$ is the linearized operator around the local Maxwellian $\mb{M}$ given by
$$
\mb{L}_\mb{M}g=Q(\mb{M}+g,\mb{M}+g)-Q(g,g)=2Q(\mb{M},g).
$$
Recall that the linearized collision operator $\mb{L}_{\mb{M}}$ is symmetric in $L^2_v$ space and the null space $\mathfrak{N}$ of $\mb{L}_\mb{M}$ is exactly spanned by $\xi_i(v)~(i=0,1,2,3,4)$.
For the hard sphere model, $\mb{L}_\mb{M}$ takes form (cf. \cite{Grad, Liu-Yang-Yu-Zhao})
\be\label{lm}
(\mb{L}_\mb{M}h)(v)=-\nu_{\mb{M}}(v)h(v)+\sqrt{\mb{M}(v)}K_{\mb{M}}\left(\left(\frac{h}{\sqrt{\mb{M}}}\right)(v)\right).
\ee
Here $K_{\mb{M}}(\cdot)=K_{2\mb{M}}(\cdot)-K_{1\mb{M}}(\cdot)$ is a symmetric and compact operator in the above weighted $L^2_v$ space. Let $k_{i\mb{M}}(v,v_*)~(i=1,2)$ be the kernel of the operator $K_{i\mb{M}}(\cdot)~(i=1,2)$. Then the collision frequency
$\nu_{\mb{M}}(v)$ and $k_{i\mb{M}}(v,v_*)$ $(i=1,2)$ have the following expressions
\begin{equation}\label{nkk}
\left\{
\begin{array}{ll}
\di \nu_{\mb{M}}(v)=\frac{2\rho}{\sqrt{2\pi R\t}}\Big\{\Big(\frac{R\t}{|v-u|}+|v-u|\Big)\int^{|v-u|}_0\exp\Big(-\frac{y^2}{2R\t}\Big)dy\\[4mm]
\di\qquad\qquad\qquad\qquad\quad +R\t\exp\Big(-\frac{|v-u|^2}{2R\t}\Big)\Big\},\\[4mm]
\di  k_{1\mb{M}}(v,v_*)=\frac{\pi\rho}{\sqrt{(2\pi R\t)^3}}|v-v_*|\exp\Big(-\frac{|v-u|^2}{4R\t}-\frac{|v_*-u|^2}{4R\t}\Big), \\[4mm]
\di  k_{2\mb{M}}(v,v_*)=\frac{4\pi\rho}{\sqrt{(2\pi R\t)^3}}|v-v_*|^{-1}\exp\Big(-\frac{|v-v_*|^2}{8R\t}-\frac{(|v|^2-|v_*|^2)^2}{8R\t|v-v_*|^2}\Big).
\end{array}
 \right.
\end{equation}
By $\eqref{nkk}_1$, we have $\nu_{\mb{M}}(v)\sim(1+|v|)$
as $|v|\rightarrow+\infty$. Furthermore, the celebrated H-theorem implies the strongly dissipative property of the linearized collision operator $\mb{L}_\mb{M}$ on the non-fluid component, i. e., there exists a positive constant $\wt\sigma_1>0$ such that for any function $g(v)\in \mathfrak{N}^\bot$ (cf. \cite{CC,Grad}), it holds that
\begin{equation}\label{H-thm}
\langle g,\mb{L}_\mb{M}g\rangle
\le
 -\wt\sigma_1\langle \nu(v)g,g\rangle,
\end{equation}
where and in the sequel the subscript $_{\mb{M}}$ in the collision frequency $\nu_{\mb{M}}(v)$ will often be omitted as the simplified form $\nu(v)$ if without confusions. By \eqref{H-thm}, the inverse of the operator $\mb{L}_\mb{M}$ exists in $\mathfrak{N}^\bot$, which together with \eqref{non-fluid}  implies that
\begin{equation}\label{Pi}
\ba
\mb{G} &=\mb{L}_\mb{M}^{-1}[\mb{P}_1(v\cdot\nabla_x\mb{M})]+\mb{L}_\mb{M}^{-1}[\mb{G}_t+\mb {P}_1(v\cdot\nabla_x\mb{G})-Q(\mb{G}, \mb{G})]\\
 :&= \mb{L}_\mb{M}^{-1}[\mb{P}_1(v\cdot\nabla_x\mb{M})] +\Pi.
\ea
\end{equation}

Substituting \eqref{Pi} into \eqref{ns-1},
we can obtain the following compressible Navier-Stokes-type equations for  the macroscopic fluid quantities $(\rho,u,\theta)$:
\begin{equation}\label{ns-2}
\left\{
\begin{array}{l}
\di \rho_{t}+{\rm div}_x (\rho u)=0,\\
\di (\rho u)_t+{\rm div}_x  (\rho u\otimes u) +\nabla_x  p\\
\di\qquad =-\int v\otimes v\cdot\nabla_x\big( \mb{L}_\mb{M}^{-1}[\mb{P}_1(v\cdot\nabla_x\mb{M})] \big)dv-\int v\otimes v\cdot\nabla_x\Pi dv,  \\
\di [\rho(\t+\f{|u|^2}{2})]_t+{\rm div}_x[\rho
u(\t+\f{|u|^2}{2})+pu]\\
\di\qquad  =-\int\f12|v|^2v\cdot\nabla_x\big( \mb{L}_\mb{M}^{-1}[\mb{P}_1(v\cdot\nabla_x\mb{M})] \big) dv
-\int\f12|v|^2v\cdot\nabla_x\Pi dv,
\end{array}
\right.
\end{equation}
A direct computation gives rise to the diffusion terms
$$
-\int v_i v_j\cdot\nabla_{x_j}\big( \mb{L}_\mb{M}^{-1}[\mb{P}_1(v\cdot\nabla_x\mb{M})] \big)dv=\sum_{j=1}^3\Big[\mu(\theta)(u_{ix_j}+u_{jx_i}-\frac23\delta_{ij}{\rm div}_x u)\Big]_{x_j},
$$
and
$$
\begin{array}{ll}
\di -\int\f12|v|^2v\cdot\nabla_x\big( \mb{L}_\mb{M}^{-1}[\mb{P}_1(v\cdot\nabla_x\mb{M})] \big) dv\\
\di\qquad =\sum_{j=1}^3(\k(\t)\t_{x_j})_{x_j}+\sum_{i,j=1}^3\Big\{\mu(\theta)u_i(u_{ix_j}+u_{jx_i}-\frac23\delta_{ij}{\rm div}_x u)\Big\}_{x_j},
\end{array}
$$
with the viscosity coefficient $\mu(\t)>0$ and the heat
conductivity coefficient $\k(\t)>0$ being the smooth functions of the
temperature $\t$.

Now we describe the rarefaction wave solution to the Euler system \eqref{euler}-\eqref{R-in} or \eqref{2ES}-\eqref{2ini} with the state equation
$$
p=\frac{2}{3}\rho\theta=k\rho^{\frac53}\exp(S),
$$
where the constant $k=\frac{1}{2\pi e}$ and $S$ is the macroscopic entropy.
It is straight to calculate that the Euler system \eqref{euler} for $(\r,u_1,\t)$ has three distinct eigenvalues
$$
\lambda_i(\r,u_1,S)=u_1+(-1)^{\f{i+1}{2}}\sqrt{p_\r(\r,S)},~ i=1,3,\qquad
\lambda_2(\r,u_1,S)=u_1,~
$$
with corresponding right eigenvectors
$$
r_i(\r,u_1,S)=((-1)^{\f{i+1}{2}}\r,\sqrt{p_\r(\r,S)},0)^t,~i=1,3, \qquad r_2(\r,u_1,S)=(p_{\scriptscriptstyle S},0,-p_\rho)^t,$$
such that
$$
r_i(\r,u_1,S)\cdot \nabla_{(\r,u_1,S)}\lambda_i(\r,u_1,S)\neq 0,~ i=1,3,
$$
and
$$
r_2(\r,u_1,S)\cdot \nabla_{(\r,u_1,S)}\lambda_2(\r,u_1,S)\equiv 0.
$$
Thus the  two $i$-Riemann invariants $\Sigma_i^{(j)}(i=1,3, j=1,2)$ can be defined by (cf. \cite{Smoller})
\begin{equation}\label{RI}
\Sigma_i^{(1)}=u_1+(-1)^{\f{i-1}{2}}\int^{\r}\f{\sqrt{p_z(z,S)}}{z}dz,\qquad
\Sigma_i^{(2)}=S,
\end{equation}
such that
$$
\nabla_{(\r,u_1,S)} \Sigma_i^{(j)}(\r,u_1,S)\cdot r_i(\r,u_1,S)\equiv0,\quad i=1,3,~ j=1,2.
$$
Given the right state $(\rho_+, u_{1+}, \theta_+)$ with $\r_+>0,\t_+>0$,  the $i$-Rarefaction wave curve $(i=1,3)$ in the phase space $(\r,u_1,\t)$ with $\r>0$ and $\t>0$ can be defined by (cf. \cite{L}):
{\small\begin{equation}\label{R3-curve}
 R_i (\rho_+, u_{1+}, \theta_+):=\Big{ \{} (\rho, u_1, \theta)\Big{ |}\lambda_{ix_1}>0, \Sigma_i^{(j)}(\r,u_1,S)=\Sigma_i^{(j)}(\r_+,u_{1+},S_+),~~j=1,2\Big{ \}}.
\end{equation}}Without the loss of generality, we consider the stability of $3-$rarefaction wave to the Euler system \eqref{euler}-\eqref{R-in} in the present paper and the stability of $1-$rarefaction wave can be considered similarly.  The $3-$rarefaction wave to the Euler system \eqref{euler}-\eqref{R-in} can be expressed explicitly through
 the Riemann solution to the inviscid Burgers equation:
\begin{equation}\label{bur}
\left\{\begin{array}{ll}
w_t+ww_{x_1}=0,\\
w(0,x_1)=\left\{\begin{array}{ll}
w_-,&x_1<0,\\
w_+,&x_1>0.
\end{array}
\right.
\end{array}
\right.
\end{equation}
If $w_-<w_+$, then the Riemann problem $(\ref {bur})$ admits a self-similar
rarefaction wave fan solution $w^r(t, x_1) = w^r(\f {x_1}{t})$ given by
\begin{equation}\label{abur}
w^r\Big(\f {x_1}{t}\Big)=\left\{\begin{array}{lr}
w_-,&\f {x_1}{t}\leq w_-,\\
\f {x_1}{t},&w_-\leq \f {x_1}{t}\leq w_+,\\
 w_+,&\f {x_1}{t}\geq w_+.
\end{array}
\right.
\end{equation}
Then the 3-rarefaction wave solution $(\r^r,u^r,\t^r)(\f {x_1}{t})$ to the
compressible Euler equations \eqref{euler}-\eqref{R-in} can be defined explicitly by
\begin{eqnarray}   \label{3-rw}
\left\{
\begin{array}{l}
\di w_\pm=\lambda_3(\r_\pm,u_{1\pm},\t_\pm), \qquad w^r(\f {x_1}{t})= \lambda_3(\r^r,u_1^r,\t^r)(\f {x_1}{t}),\\
\di
\Sigma_3^{(j)}(\r^r,u_1^r,\t^r)(\f{x_1}{t})=\Sigma_3^{(j)}(\rho_\pm,u_{1\pm},\t_\pm),\quad j=1,2,\quad
u_2^r= u_3^r=0,
\end{array} \right.
\end{eqnarray}
where $\Sigma_3^{(j)}~(j=1,2)$ are the 3-Riemann invariants defined in \eqref{RI}.

Next, we construct a smooth profile to the 3-rarefaction wave defined in \eqref{3-rw}.
Motivated by Matsumura-Nishihara \cite{MN2},  the smooth rarefaction wave can be constructed by the Burgers equation
\begin{equation}\label{dbur}
\left\{
\begin{array}{l}
\di \bar w_{t}+\bar w\bar w_{x_1}=0,\\
\di \bar w( 0,x_1
)=\bar w_0(x_1)=\f{w_++w_-}{2}+\f{w_+-w_-}{2}k_q\int_0^{\varepsilon x_1}(1+y^2)^{-q}dy,
\end{array}
\right.
\end{equation}
where $\v>0$ is a small parameter to be determined and $k_q$ is a positive constant such that $k_q\int_0^{\infty}(1+y^2)^{-q}dy=1$ for each $q\geq2$.
Note that the solution $\bar w(t,x_1)$ of the problem \eqref{dbur} can be given by
\begin{equation}\label{b-s}
\bar w(t,x_1)=\bar w_0(x_0(t,x_1)),\qquad x=x_0(t,x_1)+\bar w_0(x_0(t,x_1))t.
\end{equation}

Correspondingly, the smooth 3-rarefaction wave
$(\bar{\r},\bar{u},\bar{\t})(t,x_1)$ to compressible Euler equations
$\eqref {euler}-\eqref{R-in}$ can be defined by
\begin{eqnarray}
\left\{
\begin{array}{l}
\di w_\pm=\lambda_3(\r_\pm,u_{1\pm},\t_\pm), \qquad \bar w(x_1,1+t)= \lambda_3(\bar{\r},\bar{u}_1,\bar{\t})(t,x_1),\\
\di
\Sigma_3^{(j)}(\bar{\r},\bar{u}_1,\bar{\t})(t,x_1)=\Sigma_3^{(j)}(\rho_\pm,u_{1\pm},\t_\pm),\quad j=1,2,\quad
\bar u_2=\bar u_3=0,
\end{array} \right.\label{au}
\end{eqnarray}
where $\bar w(t, x_1)$ is the solution of Burger's equation $(\ref
{dbur})$ defined in \eqref{b-s}.
Then the planar 3-rarefaction wave $(\bar{\r},\bar{u},\bar{\t})(t,x_1)$ satisfies the Euler system
\begin{equation}\label{rare-s}
\left\{
\begin{array}{ll}
\di \bar\r_t+(\bar \r \bar u_1)_{x_1}=0,\\
\di (\bar \r \bar u_1)_t+(\bar \r\bar u_1^2+\bar p)_{x_1}=0,\\
\di  (\bar \r \bar u_i)_t+(\bar \r\bar u_1\bar u_i)_{x_1}=0,\qquad i=2,3,\\
\di (\bar\r\bar\t)_t+(\bar\r \bar u_1\bar\t)_{x_1}+\bar p\bar u_{1x_1}=0.
\end{array}
\right.
\end{equation}

The solution space to the Boltzmann equation \eqref{B}-\eqref{far-field} considered in this paper is
$H^3_x\Big(L^2_v\Big(\frac{1}{\sqrt{\mb{M_*}}}\Big)\Big)\big(\mathbb{D}\times\mathbb{R}^3\big)$
with $H^3_x$ being the standard Sobolev space and  $L^2_v\big(\frac{1}{\sqrt{\mb{M_*}}}\big)$ being weighted $L^2_v$ space defined in \eqref{product} for some global Maxwellian $\mb{M_*}$.  Now we can state our main result as follows.

\begin{theorem}\label{thm}
Suppose  $(\bar\r, \bar u, \bar\t)(t,x_1)$ is the 3-rarefaction wave defined in \eqref{au} satisfying $\di \frac12\sup_{(t,x_1)}\bar\t(t,x_1)<\inf_{(t,x_1)}\bar\t(t,x_1)$,
then there exist positive constants $\d_0, \v_0<1$ and a global Maxwellian $\mb{M_*}=\mb{M}_{[\r_*,u_*,\t_*]}$ with $\r_*>0,~\t_*>0$, such that  if $\varepsilon\leq\varepsilon_0$ and the initial values satisfy
\begin{equation}
\mathcal{E}(0):=\Big\|f_0(x,v)-\mb{M}_{[\bar\r(0,x_1),\bar u(0,x_1),\bar\t(0,x_1)]}\Big\|^2_{H^3_x\big(L^2_{v}\big(\frac{1}{\sqrt{\mb{M_*}}}\big)\big)}\leq \d_0,
\end{equation}
then 3D Boltzmann equation \eqref{B}-\eqref{far-field} admits a unique global solution $f(t,x,v)$ satisfying
\begin{equation}\label{ue}
\Big\|f(t,x,v)-\mb{M}_{[\bar\r(t,x_1),\bar u(t,x_1),\bar\t(t,x_1)]}\Big\|^2_{H^3_x\big(L^2_{v}\big(\frac{1}{\sqrt{\mb{M_*}}}\big)\big)}\leq C(\d_0+\v^{\f18}),
\end{equation}
with some uniform-in-time constant $C$ and the time-asymptotic stability of planar 3-rarefaction wave:
\begin{equation}\label{large-time}
\lim_{t\rightarrow\infty} \Big\|f(t,x,v)-\mb{M}_{[\bar\r(t,x_1),\bar u(t,x_1),\bar\t(t,x_1)]}\Big\|_{L^{\infty}_x\big(L^2_{v}\big(\frac{1}{\sqrt{\mb{M_*}}}\big)\big)}=0.
\end{equation}
Here $g(v)\in L_v^2(\f{1}{\sqrt{\mb{M}_*}})$ means that
$\frac{g(v)}{\sqrt{\mb{M_*}}}\in L_v^2(\mathbb{R}^3)$.

\begin{remark}
Theorem 1.1 is the first result on the time-asymptotic stability of basic wave patterns for the three-dimensional Boltzmann equation, even though the corresponding stability results for shock wave or contact discontinuity are still completely open.
\end{remark}


\end{theorem}

The rest part of the paper is arranged as follows. First, we present the local-in-time existence of the solution to 3D Boltzmann equation \eqref{B}-\eqref{far-field},
and list some properties for the rarefaction wave and Boltzmann equation's microscopic H-theorem in Section \ref{Preliminaries}. Then, we will prove our main result Theorem \ref{thm} based on the a priori energy estimates in Section \ref{Energy Estimates}.
Finally, we give the proof of local-in-time existence of solution to Boltzmann equation \eqref{B}-\eqref{far-field} in Appendix.

\section{Preliminaries}\label{Preliminaries}
\setcounter{equation}{0}

In this section, we first present the local-in-time existence of solution to Boltzmann equation \eqref{B}-\eqref{far-field}, list some properties for the rarefaction wave and linearized collision operator, and then show the celebrated microscopic H-theorem for Boltzmann equation.

We start from the local-in-time existence of solution to Boltzmann equation \eqref{B}-\eqref{far-field}, whose proof will be given in Appendix.
\begin{lemma}\label{local}{\rm{(Local-in-time existence)}}
For any suitable small constant $\Xi>0$, there exists a positive constant $T^*(\Xi)>0$, such that if the initial values $f_0(x,v)\geq0$ and
\begin{equation}
\mathcal{E}(0)=\Big\|f_0(x,v)-\mb{M}_{[\bar\r(0,x_1),\bar u(0,x_1),\bar\t(0,x_1)]}\Big\|_{H^3_x\big(L^2_{v}\big(\frac{1}{\sqrt{\mb{M_*}}}\big)\big)}\leq \f{\Xi}{2\sqrt{C_0}},
\end{equation}
where $C_0:=\f{1}{\min\{1,\bar\sigma\}}\geq1$ and the positive constant $\bar\sigma$ is defined in Lemma \ref{k} in Appendix.
Then 3D Boltzmann equation \eqref{B}-\eqref{far-field} admits a unique global solution
$f(t,x,v)\in H^3_x\big(L^2_{v}\big(\frac{1}{\sqrt{\mb{M_*}}}\big)\big)$ on $[0,T^*(\Xi))\times\mathbb{D}\times\mathbb{R}^3$
satisfying $f(t,x,v)\geq0$ and
\begin{equation}
\sup_{0\leq t\leq T^*(\Xi)}\Big\|f(t,x,v)-\mb{M}_{[\bar\r(t,x_1),\bar u(t,x_1),\bar\t(t,x_1)]}\Big\|_{H^3_x\big(L^2_{v}\big(\frac{1}{\sqrt{\mb{M_*}}}\big)\big)}\leq \Xi.
\end{equation}

\end{lemma}

Then we list some properties of the smooth $3-$rarefction wave constructed in \eqref{au} in the following two lemmas.
\begin{lemma}[\cite{MN, MN2}]\label{appr}
The Burgers equation $(\ref{dbur})$ has a unique smooth global solution $\bar w(t,x_1)$ such that
\begin{itemize}
\item[(1)] $w_-<\bar w(t,x_1)<w_+, \ \partial_{x_1} \bar w(t,x_1)>0,$
 \ for  $x_1\in\mathbb{R}, \ t\geq 0.$
\item[(2)] For any $\ t> 0$ and p $\in[1,\infty]$, there exists a constant $C_{p,q}$ such that
\begin{equation*}
\begin{array}{ll}
\di \|\bar w(t,\cdot)-w^r(\frac\cdot t)\|_{L^p}\leq  C_p\v^{-\frac1p}(w_+-w_-),\\[3mm]
\di
\| \bar w_{x_1}(t,\cdot)\|_{L^p}\leq  C_{p,q}\min\{\v^{1-\frac{1}{p}}(w_+-w_-),
(w_+-w_-)^{\frac1p}t^{-1+\frac1p}\}, \\[3mm]
\di
  \| \bar w_{x_1x_1}(t,\cdot)\|_{L^p}\leq
C_{p,q}\min\{\v^{2-\frac1p}(w_+-w_-),\v^{(1-\frac{1}{2q})(1-\frac1p)}(w_+-w_-)^{-\frac{p-1}{2pq}}t^{-1-\frac{p-1}{2pq}}\},\\[3mm]
\di
|\bar w_{x_1x_1}(t,x_1)|\leq C\v\bar w_{x_1}(t, x_1).
\end{array}
\end{equation*}
\item[(3)] The smooth rarefaction wave $\bar w(t,x_1)$ and the original rarefaction wave $w^r(\frac{x_1}{t})$ are time-asymptotically equivalent, i.e.,
$$
\lim_{t\rightarrow+\infty}\sup_{x_1\in \mathbb{R}} |\bar w(t,x_1)-w^r(\f{x_1}{t})|=0.
$$
\end{itemize}
\end{lemma}

\begin{lemma}[\cite{MN, MN2}]\label{appu}
Let $\d=|(\r_+-\r_-, u_+-u_-, \t_+-\t_-)|$ is the strength of the 3-rarefaction wave
$(\bar\r,\bar u, \bar\t)$ defined in (\ref{au}), then the following properties hold:
\begin{itemize}
\item[(i)] $\bar u_{1x_1}(t,x_1)>0,$
 \ for  $x_1\in\mathbb{R}, \ t\geq 0.$
\item[(ii)] The following estimates hold for all $t> 0$ and
p $\in[1,\infty]$:
\begin{equation*}\begin{array}{l}
\di \|(\bar \r, \bar u,\bar \t)(t,\cdot)-(\r^r,u^r,\t^r)(\frac \cdot t)\|_{L^p}\leq C_p\delta\v^{-\frac1p}, \\
\di \| (\bar\r,\bar u,\bar\t)_{x_1}(t,\cdot)\|_{L^p}\leq
C_{p,q} \min\{\d\v^{1-\frac1p},\d^{\frac1p}(1+t)^{-1+\frac1p}\},\\
  \|(\bar\r,\bar u,\bar\t)_{x_1x_1}(t,\cdot)\|_{L^p}\leq
C_{p,q} \min\{\d\v^{2-\frac1p},\d^{-\frac{p-1}{2pq}}\v^{(1-\frac{1}{2q})(1-\frac1p)}
(1+t)^{-1-\frac{p-1}{2pq}}\\
\di\qquad\qquad\qquad\qquad\qquad\qquad\qquad\qquad\qquad\qquad\qquad\qquad  +\d^{\frac1p}(1+t)^{-2+\frac1p}\},
\end{array}\end{equation*}
\item[(iii)] Time-asymptotically, the smooth 3-rarefaction wave and the inviscid\break 3-rarefaction wave are equivalent, i.e.,
\begin{equation*}
\lim_{t\rightarrow+\infty}\sup_{x_1\in\mathbb{R}}| (\bar\r, \bar
u, \bar\t)(t,x_1)-(\r^r,u^r,\t^r)(\f{x_1}{t})|=0.
\end{equation*}
\end{itemize}
\end{lemma}

Next, we list some lemmas on the estimates and dissipative properties of the linearized collision operator in the weighted $L^2$ space, based on the celebrated H-theorem. The first lemma can be found from \cite{LY}.

\begin{lemma}\label{Lemma 4.1} There exists a positive
constant $C$ such that
$$
\int\f{\nu(v)^{-1}Q(f,g)^2}{\wt{\mb{M}}}dv\leq
C\left\{\int\f{\nu(v)f^2}{\wt{\mb{M}}}dv\cdot\int\f{g^2}{\wt{\mb{M}}}dv+
\int\f{f^2}{\wt{\mb{M}}}dv\cdot\int\f{\nu(v)g^2}{\wt{\mb{M}}}dv\right\},
$$
where $\wt{\mb{M}}$ can be any Maxwellian so that the above
integrals are well-defined.
\end{lemma}

Based on Lemma \ref{Lemma 4.1}, the following three lemmas are taken
from \cite{Liu-Yang-Yu-Zhao}. Their proofs are straightforward by
using Lemma \ref{Lemma 4.1} and Cauchy inequality.

\begin{lemma}\label{Lemma 4.2} If $\t/2<\t_ *<\t$, then there exist two
positive constants $\wt\sigma=\wt\sigma(\r,u,\t;$ $\r_ *,u_ *,\t_ *)$ and
$\eta_0=\eta_0(\r,u,\t;\r_ *,u_ *,\t_ *)$ such that if
$|\r-\r_ *|+|u-u_ *|+|\t-\t_ *|<\eta_0$, then for
$g( v)\in  \mathfrak{N}^\bot$,
$$
-\int\f{g\mb{L}_\mb{M}g}{\mb{M}_ *}dv\geq
\wt\sigma\int\f{\nu(v)g^2}{\mb{M}_ *}dv.
$$
\end{lemma}

 \begin{lemma}\label{Lemma 4.3} Under the assumptions in Lemma \ref{Lemma 4.2}, we
have  for each $g( v)\in  \mathfrak{N}^\bot$,
$$
 \int\f{\nu(v)}{\mb{M}}|\mb{L}_\mb{M}^{-1}g|^2dv
\leq \wt\sigma^{-2}\int\f{\nu(v)^{-1}g^2}{\mb{M}}dv,$$
and
$$
 \int\f{\nu(v)}{\mb{M}_ *}|\mb{L}_\mb{M}^{-1}g|^2dv
\leq \wt\sigma^{-2}\int\f{\nu(v)^{-1}g^2}{\mb{M}_ *}dv.
$$
\end{lemma}

\begin{lemma}\label{Lemma 4.4} Under the assumptions in Lemma \ref{Lemma 4.2},  for any
positive constants $k$ and $\lambda$, it holds that
$$
|\int\f{g_1\mb{P}_1(| v|^kg_2)}{\mb{M}_ *}d v-\int\f{g_1| v|^kg_2}{\mb{M}_ *}d v|\le
C_{k}\int\f{\lambda|g_1|^2+\lambda^{-1}|g_2|^2}{\mb{M}_ *}d v,
$$
where the constant $C_{k}$ only depends on $k$.
\end{lemma}

\begin{remark}
In Lemmas \ref{Lemma 4.2}-\ref{Lemma 4.4}, $\eta_0$ may not be sufficiently small positive constant. However, in the proof of Theorem \ref{thm} in the following sections, the smallness of $\eta_0$ is crucially used to close the a priori assumptions \eqref{assumption}.
\end{remark}

\begin{lemma}[\cite{AF}]\label{sobolev} 
There exists some positive constant $C$ such that for $g\in H^2(\mathbb{D})$ with $\mathbb{D}:=\mathbb{R}\times \mathbb{T}^2$, it holds that
\begin{equation}\label{sobolev-inequality}
\|g\|^2_{L^{\infty}(\mathbb{D})}\leq C\Big[\|g\|_{L^2(\mathbb{D})}\|\nabla g\|_{L^2(\mathbb{D})}
+ \|\na g\|_{L^2(\mathbb{D})}\|\na^2 g\|_{L^2(\mathbb{D})}\Big].
\end{equation}
\end{lemma}

%
%
%
%

\section{The proof of Theorem \ref{thm}}\label{Energy Estimates}
\setcounter{equation}{0}

In this section, we prove our main result Theorem \ref{thm}, based on the local-in-time existence of the solution in Lemma \ref{local} and the a priori estimates carried out in the following. First set the perturbation around the 3-rarefaction wave by
\be\label{perturb}
\begin{cases}
(\phi,\psi,\z)(t,x)=(\r-\bar\r, u-\bar u,\t-\bar\t)(t,x),\\
\widetilde{\mb{G}}(t,x,v)=\mb{G}(t,x,v)-\bar{\mb{G}}(t,x,v),
\end{cases}
\ee
with the correction function $\bar{\mb{G}}$ as
\begin{equation}\label{bG}
\bar{\mb{G}}=\f{3}{2\t}\mb{L}_{\mb{M}}^{-1}\Big\{\mb{P}_1\Big[v_1\Big(v_1\bar
u_{1x_1}+\f{|v-u|^2}{2\t}\bar\t_{x_1}\Big)\mb{M}\Big]\Big\},
\end{equation}
due to the fact that $\|(\bar u_{1x_1}, \bar\theta_{x_1})\|^2\sim (1+t)^{-1}$ is not integrable with respect to the time $t$ but uniformly bounded in $t$. Note that the correction function $\bar{\mb{G}}$ in \eqref{bG} was first introduced in \cite{Liu-Yang-Yu-Zhao} for the stability of the rarefaction wave to the 1D Boltzmann equation, see also \cite{xin-zeng}.

Based on the local-in-time existence of the solution in Lemma \ref{local} and the standard continuum argument,  to prove the global existence
on the time interval $[0, T]$ with $T>0$ being any positive time and the uniform-in-time estimates \eqref{ue} and then Theorem \ref{thm},  it is sufficient to close the following a-priori assumptions:
\begin{equation}\label{assumption}
\begin{array}{ll}
\di \mathcal{N}(T)=\sup_{0\leq t \leq
T}\Big\{
\|(\phi,\psi,\z)(t,\cdot)\|^2_{H^2}\\[4mm]
\di \qquad\qquad +\int\int\frac{1}{\mb{M}_*}\Big(|\wt{\mb{G}}|^2+\sum_{1\leq|\alpha|\leq 2}|\partial^{\alpha}\mb{G}|^2
+\sum_{|\beta|=3}|\partial^\beta f|^2\Big)dv dx\Big\}\leq \chi^2,
\end{array}
\end{equation}
and verify \eqref{large-time}, where and in the sequel $(\partial^\alpha,\partial^\beta)=(\partial_{t,x}^\alpha,\partial_{t,x}^\beta)$
and $\chi$ is a small positive constant depending on the initial data but independent of the time $T$. Note that the global Maxellian $\mb{M}_*$
in \eqref{assumption} is determined in Theorem \ref{thm}. It can be seen easily from \eqref{assumption} and Sobolev's inequality that
$$
\mathcal{N}(0)=O(1)(\mathcal{E}(0)+\varepsilon) \quad{\rm and} \quad \mathcal{E}(0)= O(1)(\mathcal{N}(0)+\varepsilon),
$$
and
\begin{equation}\label{assum-1}
\sup_{(\tau,x)\in[0,t]\times\mathbb{D}}\sum_{0\leq |\alpha'|\leq 1}\Big(|\pa^{\alpha'}(\phi,\psi,\z)|
+\big(\int\frac{|\pa^{\alpha'}\mb{G}|^2}{\mb{M}_*} dv\big)^{\f12}\Big)
\leq C(\chi+\v),
\end{equation}
with some positive constant $C$. Under the a priori assumption \eqref{assumption}, we can prove that

\begin{theorem}\label{priori}
{\rm (A priori estimates)} Under the assumptions in Theorem \ref{thm}, then any solution $f(t,x,v)$ to
Boltzmann equation \eqref{B}-\eqref{far-field} on the time interval $[0,T]$ satisfies the following uniform-in-time a priori estimates:
\begin{equation}\label{full-es}
\begin{array}{ll}
\di \|(\phi,\psi,\z)(\cdot,t)\|_{H^2}^2+\int_0^t\Big[\|\sqrt{\bar u_{1x_1}}(\phi,\psi_1,\zeta)\|^2+\sum_{1\leq|\beta|\leq3}\|\pa^{\beta}(\phi,\psi,\z)\|^2\Big]d\tau
\\
\di\qquad +\int\int\Big(\f{|\widetilde{\mb{G}}|^2}{\mb{M}_*}+\sum_{1\leq|\a|\leq2}\f{|\pa^{\a}\mb{G}|^2}{\mb{M}_*}
+\sum_{|\beta|=3}\f{|\pa^{\beta}f|^2}{\mb{M}_*}\Big)(t,x,v) dxdv
\\
\di\qquad +\int_0^t\int\int\f{\nu(v)}{\mb{M}_*}\Big(|\wt{\mb{G}}|^2+\sum_{1\leq|\beta|\leq3}|\pa^{\beta}\mb{G}|^2\Big)dxdvd\tau
\leq C(\mathcal{N}(0)+\v^{\frac18}).
\end{array}
\end{equation}

\end{theorem}

The proof of \eqref{full-es} in Theorem \ref{priori} will be done through the suitable combination of Proposition \ref{Prop3.1} and Proposition \ref{Prop3.2} below by
multiplying \eqref{higher} with a large constant $C$ and then adding the resulting equality and \eqref{lower} together.

Once we proved Theorem \ref{priori},  we can finish the proof of Theorem \ref{thm}. The global-in-time
existence of solution  follows immediately from Lemma \ref{local} (Local-in-time existence) and Theorem \ref{priori} (A priori estimates).
Then we only need to justify the time-asymptotic stability of planar rarefaction wave as in \eqref{large-time}.
In fact, from \eqref{full-es} it holds that
\be\label{4.1}
\begin{cases}
\di \int_0^{\infty}\int\int \f{|\na_x\mb{G}|^2}{\mb{M}_*}dv dxdt+
\int_0^{\infty}\int\int \f{|\na_x(\mb{M}-\mb{M}_{[\bar\r,\bar u,\bar\t]})|^2}{\mb{M}_*}dv dxdt\leq C,\\[4mm]
\di\int_0^{\infty}\Big|\f{d}{dt}\int\int \f{|\na_x\mb{G}|^2}{\mb{M}_*}dvdx\Big|dt
\leq \int_0^{\infty}\int\int \f{\nu(v)}{\mb{M}_*}(|\na_x\mb{G}|^2+|\na_x\mb{G}_t|^2)dv dx dt\leq C,\\[4mm]
\di\int_0^{\infty}\Big|\f{d}{dt}\int\int \f{|\na_x(\mb{M}-\mb{M}_{[\bar\r,\bar u,\bar\t]})|^2}{\mb{M}_*}dvdx\Big|dt\\[4mm]
\di\qquad\leq \int_0^{\infty}\int\int \f{1}{\mb{M}_*}\big(|\na_x(\mb{M}-\mb{M}_{[\bar\r,\bar u,\bar\t]})|^2
+|\pa_t\na_x(\mb{M}-\mb{M}_{[\bar\r,\bar u,\bar\t]})|^2\big)dv dx dt\\[4mm]
\di\qquad\leq C\sum_{1\leq|\alpha|\leq 2}\int_0^{\infty}\|\pa^{\alpha}(\phi,\psi,\z)\|^2dt\leq C,
\end{cases}
\ee
which implies
\be\label{4.2}
\di \lim_{t\rightarrow\infty}\int_{\mathbb{R}^3}\int_{\mathbb{D}}\f{|\na_x\mb{G}|^2+|\na_x(\mb{M}-\mb{M}_{[\bar\r,\bar u,\bar\t]})|^2}{\mb{M}_*}dx dv=0.
\ee
Then by three-dimensional Sobolev's inequality \eqref{sobolev-inequality}, one has
\begin{equation*}
\begin{array}{l}
\di \Big\|\int \f{|\mb{G}|^2+|\mb{M}-\mb{M}_{[\bar\r,\bar u,\bar\t]}|^2}{\mb{M}_*} dv\Big\|_{L^{\infty}_x}
 \leq C\Big(\int\int \f{|\mb{G}|^2}{\mb{M}_*} dv dx\Big)^{\f12}\Big(\int\int \f{|\na_x\mb{G}|^2}{\mb{M}_*} dv dx\Big)^{\f12}\\[4mm]
\di\quad
 +C\Big(\int\int \f{|\na_x\mb{G}|^2}{\mb{M}_*} dv dx\Big)^{\f12}\Big(\int\int \f{|\na^2_x\mb{G}|^2}{\mb{M}_*} dv dx\Big)^{\f12}\\[4mm]
 \di\quad +C\Big(\int\int \f{|(\mb{M}-\mb{M}_{[\bar\r,\bar u,\bar\t]})|^2}{\mb{M}_*} dv dx\Big)^{\f12}
\Big(\int\int \f{|\na_x(\mb{M}-\mb{M}_{[\bar\r,\bar u,\bar\t]})|^2}{\mb{M}_*} dv dx\Big)^{\f12}\\[4mm]
\di\quad +C\Big(\int\int \f{|\na_x(\mb{M}-\mb{M}_{[\bar\r,\bar u,\bar\t]})|^2}{\mb{M}_*} dv dx\Big)^{\f12}
\Big(\int\int \f{|\na^2_x(\mb{M}-\mb{M}_{[\bar\r,\bar u,\bar\t]})|^2}{\mb{M}_*} dv dx\Big)^{\f12},
\end{array}
\end{equation*}
which together with \eqref{full-es} and \eqref{4.2} yields
$$
\lim_{t\rightarrow\infty}\sup_{x\in\mathbb{D}}\int \f{|\mb{G}|^2+|\mb{M}-\mb{M}_{[\bar\r,\bar u,\bar\t]}|^2}{\mb{M}_*}
(t,x,v)dv=0.
$$
Thus
\begin{equation*}
\begin{array}{ll}
\di \lim_{t\rightarrow\infty}\sup_{x\in\mathbb{D}}\int_{\mathbb{R}^3}\f{|f-\mb{M}_{[\bar\r,\bar u,\bar\t]}|^2}{\mb{M}_*}(t,x,v)dv\\
\di
\leq 2\lim_{t\rightarrow\infty}\sup_{x\in\mathbb{D}}\int_{\mathbb{R}^3} \f{|\mb{G}|^2+|\mb{M}-\mb{M}_{[\bar\r,\bar u,\bar\t]}|^2}{\mb{M}_*}
(t,x,v)dv=0,
\end{array}
\end{equation*}
which verifies \eqref{large-time}, hence the proof of Theorem \ref{thm} is completed. \hfill $\Box$

In the following subsections, we will prove the a priori estimates in Theorem \ref{priori} by the suitable combinations of the lower order estimates in  Proposition \ref{Prop3.1}  and the higher order estimates in  Proposition \ref{Prop3.2}.


\subsection{Lower order estimates}\label{low-es}

We start from the lower order estimates.

\begin{proposition} \label{Prop3.1}
Under the a priori assumption \eqref{assumption}, it holds that
\begin{equation}\label{lower}
\begin{array}{ll}
\di \|(\phi,\psi,\z)(t)\|^2+\int\int\f{|\widetilde{\mb{G}}|^2}{\mb{M}_*}(t,x,v) dxdv+\int_0^t\|\sqrt{\bar u_{1x_1}}(\phi,\psi_1,\z)\|^2d\tau\\
\di \quad+\sum_{|\alpha^\prime|=1}\int_0^t\|\partial^{\alpha^\prime}(\phi,\psi,\z)\|^2d\tau
+\int_0^t\int\int\f{\nu(v)}{\mb{M}_*}|\widetilde{\mb{G}}|^2dvdxd\tau\\
\di \leq C\|(\phi_0,\psi_0,\z_0,\na\phi_0)\|^2+C\v^{\f18}+C\int\int\f{|\widetilde{\mb{G}}|^2}{\mb{M}_*}(0,x,v) dxdv\\
\di \quad+C\|\na\phi(t)\|^2+C\int_0^t\int\int\f{\nu(v)}{\mb{M}_*}|(\mb{G}_t,\na_x\mb{G})|^2dvdxd\tau.
\end{array}
\end{equation}

\end{proposition}

\vspace*{4pt}\noindent\emph{Proof.} First, define the macroscopic entropy by
\begin{equation*}
-\frac{3}{2}\rho S=\int\mathbf{M}\ln \mathbf{M}d v.
\end{equation*}
Multiplying the equation \eqref{B} by $\ln\mathbf{M}$ and integrating
over $v$, it holds that
$$
-\left(\frac{3}{2}\rho S\right)_t-\div_x\left(\frac{3}{2}\rho u S\right)+\na_x\left(\int (v\ln\mathbf{M})\mathbf{G}d v\right)
=\int \frac{\mathbf{G}v\cdot\na_x\mathbf{M}}{\mb{M}}d v.
$$
Direct computations yields
$$
 S=-\frac{2}{3}\ln\rho+\ln (\frac{4\pi}{3}\theta)+1,\\
\qquad  p=\frac{2}{3}\rho\theta=k\rho^{\frac53}\exp(S).
$$
Denote
$$
\begin{array}{l}
\displaystyle \mathbf{X}=\Big(\rho,m_1, m_2, m_3,E\Big)^t,~{\rm with} ~ m_i=\r u_i\ (i=1,2,3), ~E=\r \Big(\t+\f{|u|^2}{2}\Big),\\
\displaystyle \mathbf{Y}=\Big(\r u,\r u u_1+p\mathbb{I}_1,\r uu_2+p\mathbb{I}_2,\r u u_3+p\mathbb{I}_3,\r u\Big(\t+\f{|u|^2}{2}\Big)+pu\Big)^t,
\end{array}
$$
where $\mathbb{I}_1=(1,0,0)^t, \mathbb{I}_2=(0,1,0)^t, \mathbb{I}_3=(0,0,1)^t$. Then the conservation law \eqref{ns-2} can be rewritten as
$$
\begin{array}{ll}
\di \mathbf{X}_t+\div\mathbf{Y}\\
\di\quad  =
\left(
\begin{array}{c}
0\\
\displaystyle \pa_j\Big[\mu(\t)\big(\pa_j u_1+\pa_1 u_j-\frac23\d_{j1}\div u\big)\Big]-\int v_1 v\cdot\na_x\Pi d v\\
\displaystyle \pa_j\Big[\mu(\t)\big(\pa_j u_2+\pa_2 u_j-\frac23\d_{j2}\div u\big)\Big]-\int v_2 v\cdot\na_x\Pi d v\\
\displaystyle \pa_j\Big[\mu(\t)\big(\pa_j u_3+\pa_3 u_j-\frac23\d_{j3}\div u\big)\Big]-\int v_3 v\cdot\na_x\Pi d v\\
\displaystyle \div\big(\kappa(\theta)\na\theta\big)+\div\Big[\mu(\t)\big(\na u+(\na u)^t-\frac23\mathbb{I}\div u\big)\Big]
-\int\frac12 | v|^2v\cdot\na_x\Pi d v
\end{array}
\right),
\end{array}
$$
where $\pa_j=\pa_{x_j}~(j=1,2,3)$ and $\mathbb{I}$ is the $3\times3$ unit matrix. Here and in the sequel $\div$ and $\na$
denote the divergence and gradient operator with respect to the spatial variable $x$ if without confusions.  Define a relative entropy-entropy flux pair $(\eta,q)$ around the local
Maxwellian $\mathbf{M}_{[\bar\rho,\bar u,\bar\theta]}$ as
$$
\left\{
\begin{array}{l}
\displaystyle \eta=\bar\theta\left\{-\frac32\rho
S+\frac32\bar\rho\bar{S}+\frac32\nabla_\mathbf{X}(\rho
S)\Big|_{\mathbf{X}=\bar{\mathbf{X}}}\cdot(\mathbf{X}-\bar{\mathbf{X}})\right\},\\
\displaystyle q_j=\bar\theta\left\{-\frac32\rho u_j
S+\frac32\bar\rho\bar{u}_j\bar{S}+\frac32\nabla_\mathbf{X}(\rho
S)\Big|_{\mathbf{X}=\bar{\mathbf{X}}}\cdot(\mathbf{Y}_j-\bar{\mathbf{Y}}_j)\right\},\quad j=1,2,3.
\end{array}
\right.
$$
Here, we can compute that
$$
\displaystyle (\rho S)_{\rho}=S+\frac{|u|^2}{2\theta}-\frac53,\qquad
\displaystyle (\rho S)_{m_i}=-\frac{u_i}{\theta},~i=1,2,3,\qquad
\displaystyle (\rho S)_{E}=\frac{1}{\theta},
$$
and then
$$
\left\{\begin{aligned}
 \eta&=\frac{3}{2}\left\{\rho\theta-\bar\theta\rho
S+\rho\Big[\Big(\bar{S}-\frac53\Big)\bar\theta+\frac{|u-\bar u|^2}{2}\Big]
+\frac23\bar\rho\bar\theta\right\} \\
&=\frac32\left(\f12\rho|u-\bar u|^2+\frac23\rho\bar\t\Psi\left(\f{\bar\r}{\rho}\right)+\rho\bar\t\Psi\left(\f{\t}{\bar\t}\right)\right),\\
q_j&=u_j\eta +(u_j-\bar u_j)(\rho\theta-\bar\rho\bar\theta),\quad j=1,2,3,
\end{aligned}
\right.
$$
where $\Psi(s)=s-\ln s-1$ is a strictly convex function around $s=1$. Then, for any $\mathbf{X}$ in the closed and bounded region of
$\sum=\{\mathbf{X}:\rho>0,\theta>0\}$, there exists a positive
constant $\wt C$ such that
$$
 {\wt C}^{-1}|(\phi,\psi,\z)|^2\leq\eta\leq
\wt C|(\phi,\psi,\z)|^2.
$$
Direct computations yield that
{\small\begin{equation}
\begin{array}{l}
\displaystyle
\eta_t+\div q+\f{3\k(\t)\bar\t}{2\t^2}|\na\z|^2
+\f{3\mu(\t)\bar\t}{2\t}\Big[\frac{\big(\na\psi+(\na\psi)^t\big)^2}{2}-\f23(\div\psi)^2\Big]
\\[3mm]
\di-\Big[\nabla_{(\bar\rho,\bar{u},\bar{S})}\eta\cdot(\bar\rho,\bar{u},\bar{S})_t
+\nabla_{(\bar\rho,\bar{u},\bar{S})}q\cdot(\bar\rho,\bar{u},\bar{S})_{x_1}\Big]
=\f32\Big\{\div\left(\f{\k(\t)\z\na\z}{\t}\right)
\nonumber\end{array}
\end{equation}\begin{equation}\label{E1}
\begin{array}{l}\di+\div\Big[\psi\mu(\t)\Big(\na u+(\na u)^t-\f23\mathbb{I}\div u\Big)\Big]
+\div\Big(\f23\mu(\t)\psi\bar u_{1x_1}\Big)-\pa_1\big(2\mu(\t)\psi_1\bar u_{1x_1}\big)\Big\}
\\[3mm]
\di +\f32\Big\{\big(\f{\k(\t)}{\t^2}+\f{\k'(\t)}{\t}\big)\z\pa_1\z\bar\t_{x_1}
+\f{\k'(\t)}{\t}\z\bar\t^2_{x_1}+\f{\k(\t)}{\t}\z\bar\t_{x_1x_1}
+\f{4\mu(\t)\z}{\t}\big(\pa_1\psi_1-\f{\div\psi}{3}\big)\bar u_{1x_1}
\\[3mm]
\di +\f{4\mu(\t)}{3\t}\z\bar u_{1x_1}^2+\f43\mu'(\t)\psi_1\bar\t_{x_1}\bar u_{1x_1}+\f43\mu(\t)\psi_1\bar u_{1x_1x_1}
+2\mu'(\t)\big[\psi_1\pa_1\z-\f13\psi\cdot\na\z\big]\bar u_{1x_1}
\Big\}\\[3mm]
\di +\f32\Big\{-\div\Big[\f{\z}{\t}\int\f12|v|^2v\Pi dv\Big]+\Big(\int\f12|v|^2v\Pi dv\Big)\cdot\na\Big(\f{\z}{\t}\Big)
-\pa_j\Big(\f{\bar\t}{\t}\psi_i\int v_iv_j\Pi dv\Big)\\[3mm]
\di +\pa_j\Big(\f{\bar\t}{\t}\psi_i\Big)\int v_iv_j\Pi dv+\div\Big[\f{\z\bar u_1}{\t}\int v_1v\Pi dv\Big]
-\Big(\int v_1 v\Pi dv\Big)\cdot\na \Big(\f{\z\bar u_1}{\t}\Big)
\Big\}.
\end{array}
\end{equation}}There exists a positive constant $C>0$ such that 
$$
\begin{array}{ll}
\di\quad -\left[\nabla_{(\bar\rho,\bar{u},\bar{S})}\eta\cdot(\bar\rho,\bar{u},\bar{S})_t
+\nabla_{(\bar\rho,\bar{u},\bar{S})}q\cdot(\bar\rho,\bar{u},\bar{S})_{x_1}\right]\\[3mm]
\di =\bar u_{1x_1}\Big[\f32\rho  (u_1-\bar u_1)^2+\f23\rho\bar\t \Psi\Big(\f{\bar\rho}\rho\Big)
+\rho\bar\t\Psi\Big(\f{\t}{\bar\t}\Big)\Big]
+\f32\bar\t_{x_1}\rho(u_1-\bar
u_1)(\f23\ln\f{\bar\r}{\r}+\ln\f{\t}{\bar\t})\\[4mm]
\di  \geq C^{-1} \bar u_{1x_1}(\phi^2+\psi_1^2+\z^2).
\end{array}
$$
Integrating \eqref{E1} with respect to $t,x$ over $[0,t]\times\mathbb{D}$ yields that
\begin{equation}\label{E2}
\ba
&\quad \|(\phi,\psi,\z)(t)\|^2+\int_0^t\Big[\|(\na\psi,\na\z)\|^2+
\|\sqrt{\bar u_{1x_1}}(\phi,\psi_1,\z)\|^2\Big]d\tau
\\
& \leq C\|(\phi_0,\psi_0,\z_0)\|^2+C\Big|\int_0^t\int\Big[\Big(\f{\k(\t)}{\t^2}+\f{\k'(\t)}{\t}\Big)\z\pa_1\z\bar\t_{x_1}
\\
&~~ +\f{4\mu(\t)\z}{\t}\Big(\pa_1\psi_1-\f{\div\psi}{3}\Big)\bar u_{1x_1}
+2\mu'(\t)\Big(\psi_1\pa_1\z-\f13\psi\cdot\na\z\Big)\bar u_{1x_1}\Big]
dxd\tau\Big|\\
&~~ +C\Big|\int_0^t\int\Big[\f{\k'(\t)}{\t}\z\bar\t^2_{x_1}+\f{\k(\t)}{\t}\z\bar\t_{x_1x_1}
+\f{4\mu(\t)}{3\t}\z\bar u_{1x_1}^2
+\f43\mu'(\t)\psi_1\bar\t_{x_1}\bar u_{1x_1}\\
&~~ +\f43\mu(\t)\psi_1\bar u_{1x_1x_1}\Big]
dxd\tau\Big|
+C\Big|\int_0^t\int\Big[\Big(\int\f12|v|^2v\Pi dv\Big)\cdot\na\Big(\f{\z}{\t}\Big)
\\
&~~  +\Big(\int v_iv_j\Pi dv\Big)\pa_j\Big(\f{\bar\t}{\t}\psi_i\Big)
-\Big(\int v_1 v\Pi dv\Big)\cdot\na \Big(\f{\z\bar u_1}{\t}\Big)\Big]dxd\tau\Big|
\\
&
:=C\|(\phi_0,\psi_0,\z_0)\|^2+\sum_{i=1}^3I_i.
\ea
\end{equation}
First, by the Cauchy's inequality and \eqref{assumption}, it holds that
\begin{equation}\label{e1}
\begin{array}{ll}
\di\quad\int_0^t\int\bar u^2_{1x_1}|\psi|^2dxd\tau
\leq C\v^{\f18}\int_0^t\int_{\mathbb{T}^2}(1+\tau)^{-\f78}\|\psi\|^2_{L^{\infty}_{x_1}(\mathbb{R})}dx_2dx_3d\tau\\
\di \leq C\v^{\f18}\int_0^t\int_{\mathbb{T}^2}(1+\tau)^{-\f78}\|\psi\|_{L_{x_1}^2(\mathbb{R})}\|\pa_1\psi\|_{L_{x_1}^2(\mathbb{R})}dx_2dx_3d\tau\\
\di \leq C\v^{\f18}\int_0^t(1+\tau)^{-\f78}\|\psi\|\|\pa_1\psi\|d\tau
 \leq C\v^{\f18}\int_0^t\|\pa_1\psi\|^2d\tau+C\v^{\f18},
\end{array}
\end{equation}
which yields
\begin{equation}\label{I1}
\begin{array}{ll}
\di
I_1\leq \f18\int_0^t\|(\na\psi,\na\z)\|^2d\tau+C\int_0^t\int\bar u_{1x_1}^2\big(|\psi|^2+\z^2\big)dxd\tau\\
\di\quad \leq \Big(\f18+C\v^{\f18}\Big)\int_0^t\|(\na\psi,\na\z)\|^2d\tau+C\v^{\f18}.
\end{array}
\end{equation}
Similar to \eqref{e1}, one has
\begin{equation}\label{I2}
\begin{array}{ll}
\di I_2 \leq C\v^{\f18} \int_0^t\int(|\psi_1|+|\z|)(\bar u^2_{1x_1}+\bar u_{1x_1x_1})dxd\tau\\
\di\quad \leq C\v^{\f18}\int_0^t\|\pa_1(\psi_1,\z)\|^2d\tau+C\v^{\f18}.
\end{array}
\end{equation}
It follows from Cauchy's inequality and \eqref{assumption} that
\begin{equation}\label{I3}
\begin{array}{ll}
\di I_3 \leq \f18\int_0^t\|(\na\psi,\na\z)\|^2d\tau+C\int_0^t\|\bar u_{1x_1}(\psi,\z)\|^2d\tau\\
\di\quad +C\ \int_0^t\int\Bigg[\Big|\int\f12|v|^2v\Pi dv\Big|^2+\sum_{i,j=1}^3\Big|\int v_iv_j\Pi dv\Big|^2\Bigg]dxd\tau.
\end{array}
\end{equation}
Note that by \eqref{Pi}, it holds that
\begin{equation}\label{I3+}
\begin{array}{ll}
\di \Big|\int \f12|v|^2v\Pi dv\Big|=\Big|\int \f12|v|^2v\mb{L}_\mb{M}^{-1}[\mb{G}_t+\mb {P}_1(v\cdot\na_x\mb{G})
-2Q(\mb{G}, \mb{G})] dv\Big| \\
\di \qquad\qquad\qquad\quad :=\sum_{i=1}^3I_{3i}.
\end{array}
\end{equation}
Choose the global Maxellian $\mb{M}_*=\mb{M}_{[\rho_*,u_*,\t_*]}$ such that
\begin{equation}\label{ma1}
\rho_*>0,\qquad \f12\t(t,x)<\t_*<\t(t,x),
\end{equation}
and
\begin{equation}\label{ma2}
|\rho(x,t)-\r_*|+|u(x,t)-u_*|+|\t(x,t)-\t_*|<\eta_0
\end{equation}
with $\eta_0$ being the small positive constant in Lemma \ref{Lemma 4.2}.
Then with such chosen $\mb{M}_*$, it holds that
\begin{equation}\label{I31}
\begin{array}{ll}
\di I_{31}=\Big|\int \f12|v|^2v\mb{L}_\mb{M}^{-1}\mb{G}_t dv\Big|\\
\di~\quad  \leq
\Big(\int\f{\nu(v)}{\mb{M}_*}|\mb{L}_\mb{M}^{-1}\mb{G}_t|^2
dv\Big)^{\f12} \Big(\int\mb{M}_*\nu(v)^{-1}\Big(\f12|v|^2v\Big)^2
dv\Big)^{\f12}\\
 \di~\quad \leq C \Big(\int\f{\nu(v)^{-1}}{\mb{M}_*}|\mb{G}_t|^2dv\Big)^{\f12},
\end{array}
\end{equation}
and
\begin{equation}\label{I32}
\begin{array}{ll}
\di I_{32}=\Big|\int \f12|v|^2v \mb{L}_\mb{M}^{-1}[\mb{P}_1(v\cdot\na_x\mb{G})] dv\Big|\\
\di~\quad
\leq C\Big(\int\f{\nu(v)}{\mb{M}_{[\r_*,u_*,2\t_*]}}|\mb{L}_\mb{M}^{-1}[\mb
{P}_1(v\cdot\na_x\mb{G})]|^2 dv\Big)^{\f12} \\[3mm]
 \di~~\quad  \leq C \Big(\int\f{\nu(v)^{-1}}{\mb{M}_{[\r_*,u_*,2\t_*]}}|\mb
{P}_1(v\cdot\na_x\mb{G})|^2 dv\Big)^{\f12}\\
\di~\quad
\leq C \Big(\int\f{\nu(v)^{-1}}{\mb{M}_*}|\na_x\mb{G}|^2 dv\Big)^{\f12}.
\end{array}
\end{equation}
Furthermore, one has
\begin{equation*}
\begin{array}{ll}
\di I_{33}=\Big|\int \f12|v|^2v \mb{L}_\mb{M}^{-1}[Q(\mb{G},\mb{G})] dv\Big|
\leq C \Big(\int\f{\nu(v)}{\mb{M}_*}|\mb{L}_\mb{M}^{-1}[Q(\mb{G}, \mb{G})]|^2 dv\Big)^{\f12} \\[3mm]
\di\quad ~~ \leq C\Big(\int\f{\nu(v)^{-1}}{\mb{M}_*}|Q(\mb{G}, \mb{G})|^2 dv\Big)^{\f12}
\leq C\Big(\int\f{\nu(v)}{\mb{M}_*}|\mb{G}|^2dv\Big)^{\f12}\Big(\int\f{|\mb{G}|^2}{\mb{M}_*}
dv\Big)^{\f12} \\
\di \leq C\Big(\int\f{\nu(v)}{\mb{M}_*}|\widetilde{\mb{G}}|^2
dv\Big)^{\f12}\Big(\int\f{|\widetilde{\mb{G}}|^2}{\mb{M}_*}
dv\Big)^{\f12}
\end{array}
\end{equation*}\begin{equation}\label{I33}
\begin{array}{ll}\di\quad+C|(\bar\t_{x_1},\bar
u_{1x})|\Big(\int\f{\nu(v)}{\mb{M}_*}|\widetilde{\mb{G}}|^2
dv\Big)^{\f12}+C|(\bar\t_{x_1},\bar u_{1x_1})|^2\\
\di \leq C(\chi+\v)\Big(\int\f{\nu(v)}{\mb{M}_*}|\widetilde{\mb{G}}|^2
dv\Big)^{\f12}+C|(\bar\t_{x_1},\bar u_{1x_1})|^2.
\end{array}
\end{equation}
Substituting \eqref{I31}-\eqref{I33} into \eqref{I3+} and then into \eqref{I3}, one can obtain
\begin{equation}\label{I3-e}
\begin{array}{ll}
\di I_3 \leq \Big(\f18+C\v^{\f18}\Big) \int_0^t\|(\na\psi,\na\z)\|^2d\tau
+C\int_0^t\int\int\f{\nu(v)}{\mb{M}_*}|(\mb{G}_t,\na_x\mb{G})|^2
dv dxd\tau\\
\di \quad+C(\chi+\v)\int_0^t\int\int\f{\nu(v)}{\mb{M}_*}|\widetilde{\mb{G}}|^2
dv dxd\tau +C\v^{\f18}.
\end{array}
\end{equation}
Substituting the estimates for $I_i(i=1,2,3)$ in \eqref{I1}, \eqref{I2} and \eqref{I3-e} into \eqref{E2}
gives the first-step lower order estimates
\begin{equation}\label{entropy-es}
\begin{array}{ll}
\di
\|(\phi,\psi,\z)(t)\|^2+\int_0^t\|(\na\psi,\na\z)\|^2d\tau
+\int_0^t\|\sqrt{\bar u_{1x_1}}(\phi,\psi_1,\z)\|^2d\tau\\
\di \leq C\|(\phi_0,\psi_0,\z_0)\|^2+C\v^{\f18}
+C\int_0^t\int\int\f{\nu(v)}{\mb{M}_*}|(\mb{G}_t,\na_x\mb{G})|^2dvdxd\tau\\
\di~\quad +C(\chi+\v)\int_0^t\int\int\f{\nu(v)}{\mb{M}_*}|\widetilde{\mb{G}}|^2
dv dxd\tau.
\end{array}
\end{equation}

Since there is no dissipation for the density function, we want to get the estimation of $\na\phi$ next.
For this, by the system \eqref{ns-1} and \eqref{rare-s}, we obtain the following form for the system of
 the perturbation $(\phi,\psi,\z)$:
\begin{equation}\label{perturb-1}
\left\{
\begin{array}{ll}
\di \phi_t+u\cdot\na\phi+\r\div\psi+\psi\cdot\na\bar\r+\phi\div\bar u=0,\\
\di \psi_t+u\cdot\na\psi+\f{2\t}{3\r}\na\phi+\f23\na\z+\psi\cdot\na\bar u+\f23\Big(\f{\t}{\r}-\f{\bar\t}{\bar\r}\Big)\na\bar\r\\[3mm]
\di~\qquad\qquad
=-\f{1}{\r}\int v\otimes v\cdot\na_x\mb{G}dv,\\[3mm]
\di\zeta_t+u\cdot\na\z+\f23\t\div\psi+\psi\cdot\na\bar\t+\f23\z\div\bar u\\[3mm]
\di\qquad \qquad =\f{1}{\r}\Big[-\int\f12|v|^2v\cdot\na_x\mb{G}dv
+u\cdot\int v\otimes v\cdot\na_x\mb{G}dv
\Big].
\end{array}\right.
\end{equation}
Then from the structure of perturbation equation $\eqref{perturb-1}_2$, it can be seen that the gradient of pressure, i.e., $\na p$ can deduce the estimation of $\na\phi$,
which is quite different from compressible Navier-Stokes system case in \cite{LW}.   More precisely,
multiplying \eqref{perturb-1}$_2$ by $\na\phi$ and integrating over $[0,t]\times\mathbb{D}$ lead to
\begin{equation}\label{phi-1}
\begin{array}{ll}
\di \int\psi\cdot\na\phi dx\Big|_0^t+\int_0^t\int\f{2\t}{3\r}|\na\phi|^2dxd\tau
\\
\di =\int_0^t\int\div\psi(u\cdot\na\phi+\r\div\psi+\psi\cdot\na\bar\r+\phi\div\bar u)dxd\tau\\[2mm]
\di\quad -\int_0^t\int\Big[u\cdot\na\psi+\f23\na\z+\psi\cdot\na\bar u
+\f23\Big(\f{\t}{\r}-\f{\bar\t}{\bar\r}\Big)\na\bar\r\\
\di\qquad \quad\qquad
+\f{1}{\r}\int v\otimes v\cdot\na_x\mb{G}dv\Big]\cdot\na\phi\, dxd\tau,
\end{array}
\end{equation}
where in the above equality we have used \eqref{perturb-1}$_1$ and the following fact
\begin{equation*}
\begin{array}{ll}
\di\int\psi_t\cdot\na\phi \,dx=\f{d}{dt}\int \psi\cdot\na\phi \,dx-\int\psi\cdot\na\phi_t\,dx
=\f{d}{dt}\int \psi\cdot\na\phi \,dx+\int \div\psi\phi_t \,dx
\end{array}
\end{equation*}\begin{equation*}
\begin{array}{ll}\di=\f{d}{dt}\int \psi\cdot\na\phi \,dx-\int\div\psi(u\cdot\na\phi+\r\div\psi+\psi\cdot\na\bar\r+\phi\div\bar u)dx.
\end{array}
\end{equation*}
It follows from \eqref{phi-1} and Cauchy's inequality that
\begin{equation}\label{phi-x}
\begin{array}{ll}
\di \int_0^t\|\na\phi\|^2d\tau\leq C\|(\psi_0,\na\phi_0)\|^2+C\|(\psi,\na\phi)(t)\|^2+C\int_0^t\|(\na\psi,\na\z)\|^2d\tau\\
\di\quad +C\v\int_0^t\|\sqrt{\bar u_{1x_1}}(\phi,\psi_1,\z)\|^2d\tau+C\int_0^t\int\int\f{\nu(v)}{\mb{M}_*}|\na_x\mb{G}|^2dvdxd\tau.
\end{array}
\end{equation}

Then we derive the estimation of $(\phi_t,\psi_t,\z_t)$. Multiplying the equation \eqref{perturb-1}$_1$ by $\phi_t$,
\eqref{perturb-1}$_2$ by $\psi_t$,  \eqref{perturb-1}$_3$ by $\z_t$ respectively, then integrating over $[0,t]\times\mathbb{D}$, we have
\begin{equation}\label{phi-t}
\begin{array}{ll}
\di \int_0^t\|(\phi_t,\psi_t,\z_t)\|^2 d\tau\leq C\int_0^t\|\na(\phi,\psi,\z)\|^2d\tau
\\
\di\qquad +C\v\int_0^t\|\sqrt{\bar u_{1x_1}}(\phi,\psi_1,\z)\|^2d\tau+C\int_0^t\int\int\frac{\nu(| v|)}{\mathbf{M}_*}|\na_x\mathbf{G}|^2
d v dxd\tau.
\end{array}
\end{equation}
Combining \eqref{entropy-es}, \eqref{phi-x} and \eqref{phi-t} together gives that
\begin{equation}\label{le1}
\begin{array}{ll}
\di
\|(\phi,\psi,\z)(t)\|^2+\int_0^t\Big[\sum_{|\alpha'|=1}\|\partial^{\alpha'}(\phi,\psi,\z)\|^2
+\|\sqrt{\bar u_{1x_1}}(\phi,\psi_1,\z)\|^2\Big]d\tau\\
\di \leq C\|(\phi_0,\psi_0,\z_0,\na\phi_0)\|^2+C \v^{\f18}
+C\int_0^t\int\int\f{\nu(v)}{\mb{M}_*}|(\mb{G}_t,\na_x\mb{G})|^2dvdxd\tau\\[3mm]
\di
\qquad+C\|\na\phi(t)\|^2+C(\chi+\v)\int_0^t\int\int\f{\nu(v)}{\mb{M}_*}|\widetilde{\mb{G}}|^2 dv dxd\tau.
\end{array}
\end{equation}
Now we do the microscopic estimates. By \eqref{non-fluid} and \eqref{perturb}, $\widetilde{\mb{G}}$ satisfies the equation
\begin{equation}\label{tilde-G}
\di\widetilde{\mb{G}}_t-\mb{L}_{\mb{M}}\widetilde{\mb{G}}=
-\f{3}{2\t}\mb{P}_1\Big[v\cdot\Big(v\cdot\na\psi+\f{|v-u|^2}{2\t}\na\z\Big)\mb{M}\Big]
-\mb{P}_1(v\cdot\na_x\mb{G})+Q(\mb{G},\mb{G})-\bar{\mb{G}}_t.
\end{equation}
Multiplying the equation \eqref{tilde-G} by
$\f{\widetilde{\mb{G}}}{\mb{M_*}}$, and then integrating over $[0,t]\times\mathbb{D}\times\mathbb{R}^3$ yield that
\begin{equation}\label{Gle0}
\begin{array}{ll}
\di \int\int\f{|\wt{\mb{G}}|^2}{2\mb{M}_*}dvdx\Big|_0^t-\int_0^t\int\int\f{\wt{\mb{G}}}{\mb{M}_*}\mb{L}_{\mb{M}}\wt{\mb{G}}dvdxd\tau\\
\di =\int_0^t\int\int\Big\{-\f{3}{2\t}\mb{P}_1\Big[v\cdot\Big(v\cdot\na\psi+\f{|v-u|^2}{2\t}\na\z\Big)\mb{M}\Big]\\
\di \qquad\qquad\qquad -\mb{P}_1(v\cdot\na_x\mb{G})+Q(\mb{G},\mb{G})-\bar{\mb{G}}_t
\Big\}\f{\wt{\mb{G}}}{\mb{M}_*}dvdxd\tau:=\sum_{i=4}^7I_i.
\end{array}
\end{equation}
By Cauchy's inequality, we have
\be\label{I4I5}
\ba
|I_4|,|I_5|&\leq \f{\tilde\s}{8}\int_0^t\int\int\f{\nu(v)}{\mb{M}_*}|\wt{\mb{G}}|^2 dvdxd\tau
+C\int_0^t\|(\na\psi,\na\z)\|^2d\tau\\
&+C\int_0^t\int\int\f{\nu(v)}{\mb{M}_*}|\na_x\mb{G}|^2 dvdxd\tau.
\ea
\ee
Similar to \eqref{I33} and by Cauchy's inequality and Lemma \ref{Lemma 4.1}, one has
\be
\ba
|I_6|&\leq  \f{\tilde\s}{8}\int_0^t\int\int\f{\nu(v)}{\mb{M}_*}|\wt{\mb{G}}|^2 dvdxd\tau
+C\int_0^t\int\int\f{\nu(v)^{-1}}{\mb{M}_*}|Q(\mb{G},\mb{G})|^2 dvdxd\tau
\nonumber\ea
\ee\be\label{I6}
\ba&\leq \f{\tilde\s}{8}\int_0^t\int\int\f{\nu(v)}{\mb{M}_*}|\wt{\mb{G}}|^2 dvdxd\tau
+C\int_0^t\int\Big(\int\f{\nu(v)}{\mb{M}_*}|\mb{G}|^2dv\cdot\int\f{|\mb{G}|^2}{\mb{M}_*}dv\Big)dxd\tau\\
&\leq \Big(\f{\tilde\s}{8}+C(\chi+\v)\Big) \int_0^t\int\int\f{\nu(v)}{\mb{M}_*}|\wt{\mb{G}}|^2 dvdxd\tau
+C\v.
\ea
\ee
It follows from the Cauchy's inequality and Lemma \ref{appu} that
\be\label{I7}
\ba
|I_7|&\leq  \f{\tilde\s}{8}\int_0^t\int\int\f{\nu(v)}{\mb{M}_*}|\wt{\mb{G}}|^2 dvdxd\tau
+C\int_0^t\int\int\f{\nu(v)}{\mb{M}_*}|\bar{\mb{G}}_t|^2 dvdxd\tau\\
&\leq \f{\tilde\s}{8}\int_0^t\int\int\f{\nu(v)}{\mb{M}_*}|\wt{\mb{G}}|^2 dvdxd\tau\\
&
\qquad+C\int_0^t\int\Big(|(\bar u_{1x_1t},\bar\t_{x_1t})|^2+|(\r_t,u_t,\t_t)|^2|(\bar u_{1x_1},\bar\t_{x_1})|^2\Big)dxd\tau\\
&\leq \f{\tilde\s}{8} \int_0^t\int\int\f{\nu(v)}{\mb{M}_*}|\wt{\mb{G}}|^2 dvdxd\tau
+C\v\int_0^t\|(\phi_t,\psi_t,\z_t)\|^2d\tau+C\v.
\ea
\ee
Substituting \eqref{I4I5}-\eqref{I7} into \eqref{Gle0} gives that
\be\label{EG}
\ba
&\quad \int\int \f{|\widetilde{\mb{G}}|^2}{\mb{M}_*}(t,x,v)dxdv
 +\int_0^t\int\int\frac{\nu(v)}{\mb{M}_*}|\wt{\mb{G}}|^2dvdxd\tau\\
& \leq C\int\int \f{|\widetilde{\mb{G}}|^2}{\mb{M}_*}(0,x,v)dxdv
+C\v\int_0^t\|(\phi_t,\psi_t,\z_t)\|^2d\tau\\
&\quad+C\int_0^t\|(\na\psi,\na\z)\|^2d\tau+C\int_0^t\int\int\f{\nu(v)}{\mb{M}_*}|\na_x\mb{G}|^2dvdxd\tau+C\v,
\ea
\ee
which along with \eqref{le1} yields \eqref{lower}, and the proof of Proposition \ref{Prop3.1} is completed.
\qed

%
%
%

\subsection{Higher order estimates}\label{high-es}

In this subsection, we will consider the higher order energy estimates.

\begin{proposition} \label{Prop3.2}
Under the a priori assumption \eqref{assumption}, it holds that
\begin{equation}\label{higher}
\begin{array}{ll}
\di \|\na(\r, u, \t)(t)\|_{H^1}^2+\int\int\Big(\sum_{1\leq|\alpha|\leq2}\f{|\pa^{\alpha}\mb{G}|^2}{\mb{M}_*}
+\sum_{|\beta|=3}\f{|\pa^{\beta}f|^2}{\mb{M}_*}\Big)(t,x,v)dvdx\\
\di+\sum_{2\leq|\beta'|\leq3}\int_0^t\|\pa^{\beta'}(\r,u,\t)\|^2d\tau
+\sum_{1\leq|\gamma|\leq3}\int_0^t\int\int\f{\nu(v)}{\mb{M}_*}|\pa^{\gamma}\mb{G}|^2dvdxd\tau\\[5mm]
\di \leq C\|\na^3\r_0\|^2+C\int\int\Big(\sum_{1\leq|\alpha|\leq2}\f{|\pa^{\alpha}\mb{G}|^2}{\mb{M}_*}
+\sum_{|\beta|=3}\f{|\pa^{\beta}f|^2}{\mb{M}_*}\Big)(0,x,v)dvdx\\
\di\quad +C\|\na(\r_0,u_0,\t_0)\|_{H^1}^2+C(\chi+\v+\eta_0)\sum_{|\alpha'|=1}\int_0^t\|\pa^{\alpha'}(\phi,\psi,\z)\|^2d\tau\\
\di\quad
+C(\chi+\v)\int_0^t\int\int\f{\nu(v)}{\mb{M}_*}|\wt{\mb{G}}|^2dvdxd\tau+C\v^{\f12}.
\end{array}
\end{equation}
\end{proposition}
The proof of Proposition \ref{Prop3.2} is divided into the following five steps.

\vspace*{2pt}$\bullet$ \underline{{\bf Step 1}}. Estimate of $\di \sup_{t\in[0,T]}\|\na(\r,u,\t)(t)\|^2$:
\be 
\begin{array}{l}
\di \quad \|\na(\r,u,\t)(t)\|^2+\sum_{|\alpha|=2}\int_0^t\|\pa^{\alpha}(\r, u,\t)\|^2d\tau\nonumber
\end{array}
\ee
\be\label{fluid-1x}
\begin{array}{l}\di \leq C\|\na(\r_0,u_0,\t_0),\na^2\r_0\|^2+C(\chi+\v)\sum_{|\alpha'|=1}\int_0^t\|\pa^{\alpha'}(\phi,\psi,\z)\|^2d\tau\\
\di
\quad +C\v^{\f12}+C\|\na^2\r(t)\|^2+C\sum_{|\alpha|=2}\int_0^t\int\int\f{\nu(v)}{\mb{M_*}}|\pa^{\alpha}\mb{G}|^2dvdxd\tau\\
\di
\quad  +C(\chi+\v)\int_0^t\int\int\f{\nu(v)}{\mb{M_*}}|(\wt{\mb{G}},\mb{G}_t,\na_x\mb{G})|^2dvdxd\tau.
\end{array}
\ee

\vspace*{4pt}\noindent\emph{Proof of \eqref{fluid-1x}.} We simply use the system \eqref{ns-2} for $(\r,u,\t)$
with diffusion terms in order to obtain the estimates of $\int_0^t\|\pa^{\alpha}(\r,u,\t)\|d\tau$ for $|\alpha|=2$
rather than the system for the perturbation $(\phi,\psi,\z)$ due to the fact that the rarefaction wave itself has the decay-in-time $\|(\bar\rho,\bar u_1,\bar\theta)_{x_1x_1}\|^2\sim(1+t)^{-2}$,
which is integrable with respect to $t$. Therefore, we rewrite the system \eqref{ns-2} as
\be\label{re-ns-Pi}
\begin{cases}
\r_t+u\cdot\na\r+\r\div u=0,\\[1mm]
\di u_t+u\cdot\na u+\frac{2\t}{3\r}\na\r+\f23\na\t=\f{\mu(\t)}{\r}\big(\Delta u+\f13\na\div u\big)\\[1mm]
\di\qquad\qquad +\frac{\mu'(\t)}{\r}\na\t\cdot\big(\na u+(\na u)^t-\f23\mathbb{I}\div u\big)
-\f{1}{\r}\int v\otimes v\cdot\na_x\Pi dv,\\[1mm]
\di\t_t+u\cdot\na\t+\f23\t\div u=\f{\k(\t)}{\r}\Delta\t\\[1mm]
\di\qquad\qquad +\f{\k'(\t)}{\r}|\na\t|^2
+\f{\mu(\t)}{\r}\Big[\f{(\na u+(\na u)^t)^2}{2}-\f23(\div u)^2\Big]\\[1mm]
\di\qquad\qquad +\f{1}{\r}\Big(-\int\f12|v|^2v\cdot\na_x\Pi dv+u\cdot\int v\otimes v\cdot\na_x\Pi dv\Big).
\end{cases}
\ee
First we estimate  $\int_0^t\|\pa^{\alpha}(u,\t)\|d\tau$ for $|\alpha|=2$.
Applying the operator $\pa_i=\partial_{x_i}~(i=1,2,3)$ to the above system yields
\be\label{ns-1x}
\begin{cases}
\di \pa_i\r_t+u\cdot\na\pa_i\r+\r\pa_i\div u+\pa_i u\cdot\na\r+\pa_i\r\div u=0,\\
\di \pa_i u_t+u\cdot\na\pa_i u+\frac{2\t}{3\r}\na\pa_i\r+\f23\na\pa_i\t
+\pa_i u\cdot\na u+\pa_i\Big(\f{2\t}{3\r}\Big)\na\r\\
\di \quad=\f{\mu(\t)}{\r}\big(\Delta \pa_i u+\f13\na\pa_i\div u\big)
+\pa_i\Big(\f{\mu(\t)}{\r}\Big)\big(\Delta u+\f13\na\div u\big)
\\
\di\qquad +\pa_i\Big(\frac{\mu'(\t)}{\r}\Big)\na\t\cdot\big(\na u+(\na u)^t-\f23\mathbb{I}\div u\big)\\
\di\qquad
+\frac{\mu'(\t)}{\r}\na\pa_i\t\cdot\big(\na u+(\na u)^t-\f23\mathbb{I}\div u\big)\\
\di\qquad +\frac{\mu'(\t)}{\r}\na\t\cdot\pa_i\big(\na u+(\na u)^t-\f23\mathbb{I}\div u\big)
-\pa_i\Big(\f{1}{\r}\int v\otimes v\cdot\na_x\Pi dv\Big),\\
\di \pa_i\t_t+u\cdot\na\pa_i\t+\f23\t\pa_i\div u+\pa_iu\cdot\na\t+\f23\pa_i\t\div u\\
\di\quad=\f{\k(\t)}{\r}\Delta\pa_i\t+\pa_i\Big(\f{\k(\t)}{\r}\Big)\Delta\t
+\pa_i\Big(\f{\k'(\t)}{\r}\Big)|\na\t|^2+\f{2\k'(\t)}{\r}\na\t\cdot\pa_i\na\t\\
\di\qquad+\pa_i\Big(\f{\mu(\t)}{\r}\Big)\Big[\f{(\na u+(\na u)^t)^2}{2}-\f23(\div u)^2\Big]
\\
\di\qquad+\f{\mu(\t)}{\r}\pa_i\Big[\f{(\na u+(\na u)^t)^2}{2}-\f23(\div u)^2\Big]\\
\di\qquad +\pa_i\Big[\f{1}{\r}\Big(-\int\f12|v|^2v\cdot\na_x\Pi dv+u\cdot\int v\otimes v\cdot\na_x\Pi dv\Big)\Big].
\end{cases}
\ee
Multiplying the equation $\eqref{ns-1x}_1$ by $\f{2\t}{3\r^2}\pa_i\r$, $\eqref{ns-1x}_2$ by $\pa_i u$
and $\eqref{ns-1x}_3$ by $\f{\pa_i\t}{\t}$, respectively, and then adding them together and integrating over $[0,t]\times\mathbb{D}$ lead to
\be\label{long-1}
\begin{array}{ll}
\di\quad\|\na(\r,u,\t)(t)\|^2-\|\na(\r_0,u_0,\t_0)\|^2+\int_0^t\|\na^2(u,\t)\|^2d\tau\\[2mm]
\di \leq C\Big|\int_0^t\int\Big\{\pa_i u\cdot\na\Big(\f{2\t}{3\r}\Big)\pa_i\r
+\Big(\f{\t}{\r^2}\div u+\f{\t_t+u\cdot\na\t}{3\r^2}\Big)|\pa_i\r|^2+\f12\div u|\pa_i u|^2\\[3mm]
\di\quad+\Big[\Big(\f{1}{2\t}\Big)_t+\div\Big(\f{u}{2\t}\Big)\Big]|\pa_i\t|^2
-(\pa_i u\cdot\na\r+\pa_i\r\div u)\f{2\t}{3\r^2}\pa_i\r\\[3mm]
\di\quad-\Big[\pa_i u\cdot\na u+\pa_i\Big(\f{2\t}{3\r}\Big)\na\r\Big]\pa_i u
-[\pa_i u\cdot\na\t+\f23\pa_i\t\div u]\f{\pa_i\t}{\t}
\Big\} dxd\tau \Big|\\[3mm]
\di+C\Big|\int_0^t\int\Big\{
\pa_i\Big(\f{\mu(\t)}{\r}\Big)\big(\Delta u+\f13\na\div u\big)\pa_i u
+\pa_i\Big(\f{\k(\t)}{\r}\Big)\Delta\t\f{\pa_i\t}{\t}\\[3mm]
\di\quad
-\na\pa_i u\cdot\na\Big(\f{\mu(\t)}{\r}\Big)\pa_i u-\pa_i\div u\,\pa_i u\cdot\na\Big(\f{\mu(\t)}{\r}\Big)
-\na\pa_i\t\cdot\na\Big(\f{\k(\t)}{\r\t}\Big)\pa_i\t\\
\di\quad
+\f{\mu'(\t)}{\r}\pa_i\Big[\na\t\cdot\big(\na u+(\na u)^t-\f23\mathbb{I}\div u\big)\Big]\pa_i u
\\
\di\quad
+\f{2\k'(\t)}{\r}\na\t\cdot\pa_i\na\t\f{\pa_i\t}{\t}
+\f{\mu(\t)}{\r}\f{\pa_i\t}{\t}\pa_i\Big[\f{\big(\na u+(\na u)^t\big)^2}{2}-\f23(\div u)^2\Big]
\Big\}dxd\tau \Big|\\
\di +C\Big|\int_0^t\int\Big\{
\pa_i\Big(\f{\mu'(\t)}{\r}\Big)\na\t\cdot\big[\na u+(\na u)^t-\f23\mathbb{I}\div u\big]\pa_i u
+\pa_i\Big(\f{\k'(\t)}{\r}\Big)|\na\t|^2\f{\pa_i\t}{\t}\\
\di\quad +\pa_i\Big(\f{\mu(\t)}{\r}\Big)\Big[\f{\big(\na u+(\na u)^t\big)^2}{2}-\f23(\div u)^2\Big]\f{\pa_i\t}{\t}
\Big\}dxd\tau\Big|\\[3mm]
\di\quad
+C\Big|\int_0^t\int \Big\{\f{1}{\r}\Big(\int v\otimes v\cdot\na_x\Pi dv\Big)\pa_i^2 u\\[3mm]
\di\quad-\f{1}{\r}\Big(-\int\f12|v|^2v\cdot\na_x\Pi dv+ u\cdot\int v\otimes v\cdot\na_x\Pi dv\Big)\pa_i\Big(\f{\pa_i\t}{\t}\Big)
\Big\}dxd\tau\Big|\\[3mm]
\di :=\sum_{i=8}^{11}I_i.
\end{array}
\ee
%
We now estimate $I_8$ in \eqref{long-1}. By Lemma \ref{appu}, we can estimate the triple nonlinear terms as
\begin{equation*}
\ba
&\int_0^t\int|\na\r||\na u||\na\t|dxd\tau
\leq \int_0^t\int(|\na\phi|+|\bar\r_{x_1}|)(|\na\psi|+\bar u_{1x_1})(|\na\z|+|\bar\t_{x_1}|)dxd\tau\\
&\leq \int_0^t\int\big[|\na\phi||\na\psi||\na\z|+\bar u_{1x_1}(|\na\phi||\na\z|+|\na\psi||\na\z|+|\na\phi||\na\psi|)\\
&\quad\quad\quad\quad +\bar u_{1x_1}^2(|\na\phi|+|\na\psi|+|\na\z|)+\bar u_{1x_1}^3\big]dxd\tau \\
&\leq C(\chi+\v)\int_0^t\|\na(\phi,\psi,\z)\|^2d\tau+C\v^{\f12}.
\ea
\end{equation*}
Then we have
\be\label{I8}
I_8\leq C(\chi+\v)\int_0^t(\|\na(\phi,\psi,\z)\|^2+\|\z_t\|^2)d\tau+C\v.
\ee
We then estimate the typical term in $I_9$ as
\begin{equation}\label{I9}
\begin{array}{l}
\di\quad\int_0^t\int|\na^2 u||\na\t||\na u|dxd\tau\leq C\int_0^t\int|\na^2 u|(|\na\z|+|\bar\t_{x_1}|)(|\na\psi|+\bar u_{x_1})dxd\tau\\
\di\leq C\int_0^t\int|\na^2 u|[|\na\z||\na\psi|+\bar u_{1x_1}(|\na\psi|+|\na\z|+\bar u^2_{1x_1})]dxd\tau\\
\di\leq C(\chi+\v)\int_0^t(\|\na^2 u\|^2+\|\na(\psi,\z)\|^2)d\tau+C\v^{\f12}.
\end{array}
\end{equation}
The term $I_{10}$ is the higher nonlinear term compared with $I_8$ and $I_9$ and is easier to estimate. Similar to \eqref{I3+}, it holds that
\be\label{Pix}
\begin{array}{ll}
\di\quad\Big|\int v\otimes v\cdot\na_x\Pi dv\Big| \\
\di \leq C\Big(\int\f{\nu(v)}{\mb{M}_*}|(\na_x \mb{G}_t, \na_x^2\mb{G} )|^2dv\Big)^{\f12}+C(\chi+\v) \Big(\int\f{\nu(v)}{\mb{M}_*}|\wt{\mb{G}}|^2dv\Big)^{\f12}
\\ \di\quad +C(\chi+\v) \Big(\int\f{\nu(v)}{\mb{M}_*}|(\mb{G}_t, \na_x\mb{G})|^2dv\Big)^{\f12}+C|(\bar u_{1x_1},\bar \t_{x_1})|^2,
\end{array}
\ee
which together with Cauchy's inequality gives
\be\label{I11}
\ba
I_{11}&\leq C\int_0^t(\|\na^2(u,\t)\|+\|\na\t\|^2_{L^4})\|\int v\otimes v\cdot\na_x\Pi dv\|d\tau\\
&\leq \f14\int_0^t\|\na^2(u,\t)\|^2d\tau+C\int_0^t\int\int\f{\nu(v)}{\mb{M}_*}|(\na_x\mb{G}_t,\na_x^2\mb{G})|^2dvdxd\tau+C\v\\
&\quad +C(\chi+\v)\int_0^t\|\na\z\|^2d\tau+C(\chi+\v)\int_0^t\int\int\f{\nu(v)}{\mb{M}_*}|(\mb{G}_t,\na_x\mb{G},\wt{\mb{G}})|^2dvdxd\tau.
\ea
\ee
Substituting \eqref{I8}, \eqref{I9} and \eqref{I11} into \eqref{long-1} yields
\be\label{de-1}
\ba
&\quad \|\na(\r,u,\t)(t)\|^2+\int_0^t\|\na^2(u,\t)\|^2d\tau\\
& \leq C\|\na(\r_0,u_0,\t_0)\|^2+C\v^{\f12}+C(\chi+\v)\int_0^t\int\int\f{\nu(v)}{\mb{M_*}}|(\wt{\mb{G}},\mb{G}_t,\na_x\mb{G})|^2dvdxd\tau
\\
& +C(\chi+\v)\sum_{|\alpha'|=1}\int_0^t\|\pa^{\alpha'}(\phi,\psi,\z)\|^2d\tau
+C\int_0^t\int\int\f{\nu(v)}{\mb{M_*}}|(\na_x\mb{G}_t,\na_x^2\mb{G})|^2dvdxd\tau.
\ea
\ee
Next, we estimate $\di\int_0^t\|\nabla^2\r\|^2d\tau$.
For this, we use the system \eqref{ns-1} and rewrite it as
\be\label{re-ns-G}
\begin{cases}
\r_t+u\cdot\na\r+\r\div u=0,\\
\di u_t+u\cdot\na u+\frac{2\t}{3\r}\na\r+\f23\na\t=-\f{1}{\r}\int v\otimes v\cdot\na_x\mb{G}dv,\\[3mm]
\di\t_t+u\cdot\na\t+\f23\t\div u= \f{1}{\r}\big(-\int\f12|v|^2v\cdot\na_x\mb{G}dv+u\cdot\int v\otimes v\cdot\na_x\mb{G} dv\big).
\end{cases}
\ee
Then applying $\pa_i=\partial_{x_i}~(i=1,2,3)$ to the above system leads to
\be\label{ns-1x-G}
\begin{cases}
\di \pa_i\r_t+u\cdot\na\pa_i\r+\r\pa_i\div u+\pa_i u\cdot\na\r+\pa_i\r\div u=0,\\
\di \pa_i u_t+u\cdot\na\pa_i u+\frac{2\t}{3\r}\na\pa_i\r+\f23\na\pa_i\t
+\pa_i u\cdot\na u+\pa_i\Big(\f{2\t}{3\r}\Big)\na\r\\
\di \qquad=\f{\pa_i\r}{\r^2}\int v\otimes v\cdot\na_x\mb{G} dv-\f{1}{\r}\int v\otimes v\cdot\na_x\pa_i\mb{G} dv,\\
\di \pa_i\t_t+u\cdot\na\pa_i\t+\f23\t\pa_i\div u+\pa_iu\cdot\na\t+\f23\pa_i\t\div u\\[3mm]
\di\qquad =-\f{\pa_i\r}{\r^2}\Big(-\int\f12|v|^2v\cdot\na_x\mb{G} dv+u\cdot\int v\otimes v\cdot\na_x\mb{G} dv\Big)\\[4mm]
\di\qquad\quad +\f{1}{\r}\Big(-\int\f12|v|^2v\cdot\na_x\pa_i\mb{G} dv+u\cdot\int v\otimes v\cdot\na_x\pa_i\mb{G} dv\\
\di\qquad\qquad\quad\
+\pa_i u\cdot\int v\otimes v\cdot\na_x\mb{G}dv\Big).
\end{cases}
\ee
Multiplying the equation \eqref{ns-1x-G}$_2$ by $\na\pa_i\r$ and integrating over $[0,t]\times\mathbb{D}$ yield
\begin{equation}\label{phi-2}
\begin{array}{ll}
\di\quad \int\pa_i u\cdot\na\pa_i\r \,dx\Big|_0^t+\int_0^t\int\f{2\t}{3\r}|\na\pa_i\r|^2dxd\tau\\
\di=\int_0^t\int\pa_i\div u(u\cdot\na\pa_i\r+\r\pa_i\div u+\pa_i u\cdot\na\r+\pa_i\r\div u)dxd\tau\\[2mm]
\di-\int_0^t\int\Big[u\cdot\na\pa_i u+\f23\na\pa_i\t+\pa_i u\cdot\na u+\pa_i\Big(\f{2\t}{3\r}\Big)\na\r\Big]\cdot\na\pa_i\r\,dxd\tau\\
\di+\int_0^t\int\Big(\f{\pa_i\r}{\r^2}\int v\otimes v\cdot\na_x\mb{G} dv-\f{1}{\r}\int v\otimes v\cdot\na_x\pa_i\mb{G} dv\Big)\cdot\na\pa_i\r\,dxd\tau,
\end{array}
\end{equation}
where in the above equality we have used \eqref{ns-1x-G}$_1$ and the following identity
\be
\begin{array}{ll}
\di\int\pa_i u_t\cdot\na\pa_i\r \,dx=\f{d}{dt}\int \pa_i u\cdot\na\pa_i\r \,dx-\int\pa_i u\cdot\na\pa_i\r_t\,dx\\
\di\quad
=\f{d}{dt}\int \pa_i u\cdot\na\pa_i\r \,dx+\int \pa_i\div u\pa_i\r_t \,dx\\
\di\quad=\f{d}{dt}\int \pa_i u\cdot\na\pa_i\r \,dx
-\int\pa_i\div u(u\cdot\na\pa_i\r+\r\pa_i\div u+\pa_i u\cdot\na\r+\pa_i\r\div u)dx.
\end{array}
\ee
By \eqref{phi-2} and Cauchy's inequality, it holds that
\be\label{rho-2x}
\begin{array}{l}
\di \int_0^t\|\na^2\r\|^2d\tau \leq C\|(\na u_0, \na^2\r_0)\|^2+C\v+C\|(\na u,\na^2\r)(t)\|^2
\\
\di\quad+C\int_0^t\|\na^2(u,\t)\|^2d\tau
+C\int_0^t\int\int\f{\nu(v)}{\mb{M}_*}|\na_x^2\mb{G}|^2dvdxd\tau\\
\di\quad
+C(\chi+\v)\int_0^t\|\na(\phi,\psi,\z)\|^2d\tau+C(\chi+\v)\int_0^t\int\int\f{\nu(v)}{\mb{M_*}}|\na_x\mb{G}|^2dvdxd\tau,
\end{array}
\ee
which together with \eqref{de-1} implies
\be\label{1x}
\begin{array}{l}
\di \|\na(\r,u,\t)(t)\|^2+\int_0^t\|\na^2(\r, u,\t)\|^2d\tau\\
\di \leq C\|\big(\na\r_0,\na u_0,\na\t_0,\na^2\r_0\big)\|^2+C\v^{\f12}+C\|\na^2\r(t)\|^2\\
\di+C(\chi+\v)\int_0^t\|\na(\phi,\psi,\z)\|^2d\tau+C\sum_{|\alpha|=2}\int_0^t\int\int\f{\nu(v)}{\mb{M_*}}|\pa^{\alpha}\mb{G}|^2dvdxd\tau
\end{array}
\ee\be
\begin{array}{l}\di+C(\chi+\v)\int_0^t\int\int\f{\nu(v)}{\mb{M_*}}|(\wt{\mb{G}},\mb{G}_t,\na_x\mb{G})|^2dvdxd\tau.\nonumber
\end{array}
\ee
Then we will derive the estimations of $\|\na(\r_t, u_t,\t_t)\|$ and $\|(\r_{tt}, u_{tt}, \t_{tt})\|$.
First, multiplying the equation \eqref{ns-1x-G}$_1$ by $\pa_i\r_t$, \eqref{ns-1x-G}$_2$ by $\pa_i u_t$,
\eqref{ns-1x-G}$_3$ by $\pa_i\t_t~(i=1,2,3)$ respectively, then integrating over $[0,t]\times\mathbb{D}$, we have
\begin{equation}\label{phi-xt}
\begin{array}{l}
\di \int_0^t\|\na(\r_t, u_t,\t_t)\|^2 d\tau \leq C\int_0^t\int\int\frac{\nu(v)}{\mathbf{M}_*}|\na^2_x\mathbf{G}|^2
d v dxd\tau+C\int_0^t\|\na^2(\r,u,\t)\|^2d\tau
 \\
\di\quad+C(\chi+\v)\int_0^t\|\na(\phi,\psi,\z)\|^2d\tau +C(\chi+\v)\int_0^t\int\int\f{\nu(v)}{\mb{M}_*}|\na_x\mb{G}|^2dvdxd\tau+C\v.
\end{array}
\end{equation}

Next, applying $\pa_t$ to the system \eqref{re-ns-G} gives
\be\label{ns-2t-G}
\begin{cases}
\di\r_{tt}+u\cdot\na\r_t+\r\div u_t+ u_t\cdot\na\r+\r_t\div u=0,\\
\di u_{tt}+u\cdot\na u_t+\frac{2\t}{3\r}\na\r_t+\f23\na\t_t
+ u_t\cdot\na u+\Big(\f{2\t}{3\r}\Big)_t\na\r\\
\di \quad=\f{\r_t}{\r^2}\int v\otimes v\cdot\na_x\mb{G} dv-\f{1}{\r}\int v\otimes v\cdot\na_x\mb{G}_t dv,\\
\di\t_{tt}+u\cdot\na\t_t+\f23\t\div u_t+ u_t\cdot\na\t+\f23\t_t\div u\\[3mm]
\di\quad=-\f{\r_t}{\r^2}\Big(-\int\f12|v|^2v\cdot\na_x\mb{G} dv+u\cdot\int v\otimes v\cdot\na_x\mb{G} dv\Big)\\[3mm]
\di\qquad +\f{1}{\r}\Big(-\int\f12|v|^2v\cdot\na_x\mb{G}_t dv+u\cdot\int v\otimes v\cdot\na_x\mb{G}_t dv
+ u_t\cdot\int v\otimes v\cdot\na_x\mb{G}dv\Big).
\end{cases}
\ee
Then, multiplying the equation \eqref{ns-2t-G}$_1$ by $\r_{tt}$, \eqref{ns-2t-G}$_2$ by $ u_{tt}$,
\eqref{ns-2t-G}$_3$ by $\t_{tt}$ respectively, then adding them together and integrating over $[0,t]\times\mathbb{D}$ lead to
\begin{equation}\label{phi-tt}
\begin{array}{l}
\di\int_0^t\|(\r_{{tt}}, u_{tt},\t_{tt})\|^2 d\tau \leq C\int_0^t\|\na(\r_t,u_t,\t_t)\|^2d\tau\\
\di\quad
 +C(\chi+\v)\sum_{|\alpha'|=1}\int_0^t\|\pa^{\alpha'}(\phi,\psi,\z)\|^2d\tau+C\int_0^t\int\int\frac{\nu(| v|)}{\mathbf{M}_*}|\na_x\mathbf{G}_t|^2
d v dxd\tau\\
\di\quad
+C(\chi+\v)\int_0^t\int\int\f{\nu(v)}{\mb{M}_*}|\na_x\mb{G}|^2dvdxd\tau+C\v.
\end{array}
\end{equation}
where we have used that
\be
\ba
 \int_0^t\int|u_t\cdot\na\r\,\r_{tt}|dx d\tau&\leq C\int_0^t\int\big[|\psi_t||\na\phi|+(|\psi_t|
 +|\na\phi|)\bar u_{1x_1}+\bar u^2_{1x_1}\big]|\r_{tt}|dxd\tau\\
&\leq C(\chi+\v)\int_0^t(\|\r_{tt}\|^2+\|(\psi_t,\na\phi)\|^2)d\tau+C\v.
\ea
\ee
Now \eqref{1x}, \eqref{phi-xt} and \eqref{phi-tt} can deduce \eqref{fluid-1x}directly. \hfill $\Box$

$\bullet$ \underline{{\bf Step 2}}.  Estimate of $\di \sup_{t\in[0,T]}\int\int \f{|\pa^{\alpha'}\mb{G}|^2}{\mb{M}_*}(t,x,v)dv dx$~($|\alpha'|=1$):
\be\label{eG-1x}
\begin{array}{l}
\di \int\int \f{|\pa^{\alpha'}\mb{G}|^2}{\mb{M}_*}(t,x,v)dv dx
+\int_0^t\int\int\f{\nu(v)}{\mb{M}_*}|\pa^{\alpha'}\mb{G}|^2dvdxd\tau\\
\di \leq C \int\int \f{|\pa^{\alpha'}\mb{G}|^2}{\mb{M}_*}(0,x,v)dv dx+C\int_0^t\|\pa^{\alpha'}\na(u,\t)\|^2d\tau
\end{array}
\ee\be
\begin{array}{l}\di
+C(\chi+\v)\int_0^t\|\pa^{\alpha'}(\phi,\psi,\z)\|^2d\tau+C\sum_{|\alpha|=2}\int_0^t\int\int\f{\nu(v)}{\mb{M}_*}|\pa^{\alpha}\mb{G}|^2dvdxd\tau\nonumber\\
\di
+C(\chi+\v)\int_0^t\int\int\f{\nu(v)}{\mb{M}_*}|\wt{\mb{G}}|^2dvdxd\tau+C\v.\nonumber
\end{array}
\ee

\vspace*{4pt}\noindent\emph{Proof of \eqref{eG-1x}.} Applying $\pa^{\alpha'}$ to the equation \eqref{non-fluid} for $|\alpha'|=1$, one has,
\begin{equation}\label{G-1x}
\begin{array}{l}
\di \pa^{\alpha'}\mb{G}_t -\mb{L}_\mb{M}\pa^{\alpha'}\mb{G}\\
\di =
-\pa^{\alpha'}\mb{P}_1(v\cdot\nabla_x\mb{M})-\pa^{\alpha'}\mb{P}_1(v\cdot\nabla_x\mb{G})
+2Q(\pa^{\alpha'}\mb{M},\mb{G})+ 2Q(\pa^{\alpha'}\mb{G}, \mb{G}).
\end{array}
\end{equation}
Multiplying the above equation by $\f{\pa^{\alpha'}{\mb{G}}}{\mb{M_*}}$, and integrating over $[0,t]\times\mathbb{D}\times\mathbb{R}^3$ gives
\begin{equation}\label{es-G-1x}
\begin{array}{ll}
\di \int\int\f{|\pa^{\alpha'}\mb{G}|^2}{2\mb{M}_*}dvdx\Big|_0^t
-\int_0^t\int\int\f{\pa^{\alpha'}\mb{G}}{\mb{M}_*}\mb{L}_{\mb{M}}\pa^{\alpha'}\mb{G}dvdxd\tau\\
\di=\int_0^t\int\int\Big[
-\pa^{\alpha'}\mb{P}_1(v\cdot\nabla_x\mb{M})-\pa^{\alpha'}\mb{P}_1(v\cdot\nabla_x\mb{G})
+2Q(\pa^{\alpha'}\mb{M},\mb{G})\\
\di\qquad \qquad\qquad + 2Q(\pa^{\alpha'}\mb{G}, \mb{G}) \Big]\f{\pa^{\alpha'}\mb{G}}{\mb{M}_*}dvdxd\tau:=\sum_{i=12}^{15}I_i.
\end{array}
\end{equation}
Note that $\mb{P}_1(v\cdot\nabla_x\mb{M})$ does not contain the density gradient $\na\r$,
thus by Cauchy's inequality, we have
\be\label{I12}
\ba
|I_{12}|& \leq C\int_0^t\int \big[|\pa^{\alpha'}\na(u,\t)|+|\na(\r,u,\t)||\pa^{\alpha'}(\r,u,\t)|\big]
\Big(\int\f{\nu(v)}{\mb{M}_*}|\pa^{\alpha'}\mb{G}|^2dv\Big)^{\f12}dxd\tau\\
& \leq \f{\tilde\s}{8}\int_0^t\int\int\f{\nu(v)}{\mb{M}_*}|\pa^{\alpha'}\mb{G}|^2 dvdxd\tau
+C\int_0^t\|\pa^{\alpha'}\na(u,\t)\|^2d\tau\\
&\quad+C(\chi+\v)\int_0^t\|\pa^{\alpha'}(\phi,\psi,\z)\|^2d\tau+C\v.
\ea
\ee
Similarly, one has,
\be\label{I13}
\begin{array}{l}
\di |I_{13}| \leq \f{\tilde\s}{8}\int_0^t\int\int\f{\nu(v)}{\mb{M}_*}|\pa^{\alpha'}\mb{G}|^2 dvdxd\tau
+C\int_0^t\int\int\f{\nu(v)}{\mb{M}_*}|\pa^{\alpha'}\na_x\mb{G}|^2 dvdxd\tau\\
\di~\quad +C\int_0^t\int |\pa^{\alpha'}(\r,u,\t)|
\Big(\int\f{\nu(v)}{\mb{M}_*}|\na_x\mb{G}|^2dv\Big)^{\f12}
\Big(\int\f{\nu(v)}{\mb{M}_*}|\pa^{\alpha'}\mb{G}|^2dv\Big)^{\f12}dxd\tau\\[3mm]
\di \quad\leq \Big(\f{\tilde\s}{8}+C(\chi+\v)\Big)\int_0^t\int\int\f{\nu(v)}{\mb{M}_*}|\pa^{\alpha'}\mb{G}|^2 dvdxd\tau\\
\di~\quad
+C\int_0^t\int\int\f{\nu(v)}{\mb{M}_*}|\pa^{\alpha'}\na_x\mb{G}|^2 dvdxd\tau.
\end{array}
\ee
The Cauchy's inequality, Lemma \ref{Lemma 4.1} and assumption \eqref{assum-1} yield
\be\label{I14}
\begin{array}{l}
\di |I_{14}|\leq  \f{\tilde\s}{8}\int_0^t\int\int\f{\nu(v)}{\mb{M}_*}|\pa^{\alpha'}\mb{G}|^2 dvdxd\tau
+C\int_0^t\int\int\f{\nu(v)^{-1}}{\mb{M}_*}|Q(\pa^{\alpha'}\mb{M},\mb{G})|^2 dvdxd\tau\\
\di \leq \f{\tilde\s}{8}\int_0^t\int\int\f{\nu(v)}{\mb{M}_*}|\pa^{\alpha'}\mb{G}|^2 dvdxd\tau
\\\di ~+C\int_0^t\int\Big(\int\f{\nu(v)}{\mb{M}_*}|\pa^{\alpha'}\mb{M}|^2dv\cdot\int\f{\nu(v)}{\mb{M}_*}|\mb{G}|^2dv\Big)dxd\tau
\end{array}
\ee\be
\begin{array}{l}\di \leq \f{\tilde\s}{8}\int_0^t\int\int\f{\nu(v)}{\mb{M}_*}|\pa^{\alpha'}\mb{G}|^2 dvdxd\tau
+C\int_0^t\int|\pa^{\alpha'}(\r,u,\t)|^2\Big(\int\f{\nu(v)}{\mb{M}_*}|\mb{G}|^2dv\Big)dxd\tau\nonumber\\
\di \leq \f{\tilde\s}{8} \int_0^t\int\int\f{\nu(v)}{\mb{M}_*}|\pa^{\alpha'}\mb{G}|^2 dvdxd\tau
+C(\chi+\v)\int_0^t\int\int\f{\nu(v)}{\mb{M}_*}|\wt{\mb{G}}|^2dvdxd\tau\nonumber\\
\di \quad+C(\chi+\v)\int_0^t\|\pa^{\alpha'}(\phi,\psi,\z)\|^2d\tau+C\v\nonumber
\end{array}
\ee
and
\be\label{I15}
\ba
|I_{15}|&\leq  \f{\tilde\s}{8}\int_0^t\int\int\f{\nu(v)}{\mb{M}_*}|\pa^{\alpha'}\mb{G}|^2 dvdxd\tau
+C\int_0^t\int\int\f{\nu(v)^{-1}}{\mb{M}_*}|Q(\pa^{\alpha'}\mb{G},\mb{G})|^2 dvdxd\tau\\
&\leq \f{\tilde\s}{8}\int_0^t\int\int\f{\nu(v)}{\mb{M}_*}|\pa^{\alpha'}\mb{G}|^2 dvdxd\tau
+C\int_0^t\int\Big(\int\f{\nu(v)}{\mb{M}_*}|\pa^{\alpha'}\mb{G}|^2dv\cdot\int\f{|\mb{G}|^2}{\mb{M}_*}dv\\
&\quad+\int\f{|\pa^{\alpha'}\mb{G}|^2}{\mb{M}_*}dv\cdot\int\f{\nu(v)}{\mb{M}_*}|\mb{G}|^2dv\Big)dxd\tau\\
&\leq \f{\tilde\s}{8} \int_0^t\int\int\f{\nu(v)}{\mb{M}_*}|\pa^{\alpha'}\mb{G}|^2 dvdxd\tau\\
&\di \quad
+C(\chi+\v)\int_0^t\int\int\f{\nu(v)}{\mb{M}_*}|(\pa^{\alpha'}\mb{G},\wt{\mb{G}})|^2 dvdxd\tau.
\ea
\ee
Substituting \eqref{I12}-\eqref{I15} into \eqref{es-G-1x} gives \eqref{eG-1x}.
\qed

Combining \eqref{fluid-1x} and \eqref{eG-1x} together, it holds that
\be\label{order-1}
\begin{array}{l}
\di \|\na(\r,u,\t)(t)\|^2
+\sum_{|\alpha'|=1}\int\int \f{|\pa^{\alpha'}\mb{G}|^2}{\mb{M}_*}(t,x,v)dv dx
+\sum_{|\alpha|=2}\int_0^t\|\pa^{\alpha}(\r, u,\t)\|^2d\tau\\
\di\quad+\sum_{|\alpha'|=1}\int_0^t\int\int\f{\nu(v)}{\mb{M}_*}|\pa^{\alpha'}\mb{G}|^2dvdxd\tau\\
\di \leq C\|\na\r_0,\na u_0,\na\t_0,\na^2\r_0\|^2+C \sum_{|\alpha'|=1}\int\int \f{|\pa^{\alpha'}\mb{G}|^2}{\mb{M}_*}(0,x,v)dv dx
+C\|\na^2\r(t)\|^2\\
\di\quad+C(\chi+\v)\sum_{|\alpha'|=1}\int_0^t\|\pa^{\alpha'}(\phi,\psi,\z)\|^2d\tau
+C\sum_{|\alpha|=2}\int_0^t\int\int\f{\nu(v)}{\mb{M_*}}|\pa^{\alpha}\mb{G}|^2dvdxd\tau\\
\di\quad+C(\chi+\v)\int_0^t\int\int\f{\nu(v)}{\mb{M_*}}|\wt{\mb{G}}|^2dvdxd\tau+C\v^{\f12}.
\end{array}
\ee
$\bullet$ \underline{{\bf Step 3}}.  Estimate of $\di \sup_{t\in[0,T]} \|\na^2(\r,u,\t)(t)\|^2$:
\be\label{fluid-2x}
\begin{array}{l}
\di \|\na^2(\r,u,\t)(t)\|^2+\sum_{|\beta|=3}\int_0^t\|\na^{\beta}(\r,u,\t)\|^2d\tau
\leq C\|(\na^2\r_0,\na^2 u_0,\na^2 \t_0,\na^3\r_0)\|^2\\
\di\quad+C(\chi+\v)\sum_{1\leq|\alpha|\leq2}\int_0^t\|\pa^{\alpha}(\phi,\psi,\z)\|^2d\tau
+C\sum_{|\beta|=3}\int_0^t\int\int\f{\nu(v)}{\mb{M}_*}|\pa^{\beta}\mb{G}|^2dvdxd\tau\\
\di\quad+C(\chi+\v)\int_0^t\int\int\f{\nu(v)}{\mb{M}_*}\Big(\sum_{1\leq|\alpha|\leq2}|\pa^{\alpha}\mb{G}|^2
+|\wt{\mb{G}}|^2\Big)dvdxd\tau+C\|\na^3\r(t)\|^2+C\v.
\end{array}
\ee
\vspace*{4pt}\noindent\emph{Proof of \eqref{fluid-2x}.}
First, applying $\pa_j~(j=1,2,3)$ to system \eqref{ns-1x} leads to
\be\label{ns-2x-Pi}
\begin{cases}
\di \pa_j\pa_i\r_t+u\cdot\na\pa_j\pa_i\r+\r\pa_j\pa_i\div u
+(\pa_j u\cdot\na\pa_i\r+\pa_i u\cdot\na\pa_j\r+\div u\pa_j\pa_i\r)\\
~+(\pa_j\r\pa_i\div u+\pa_j\pa_i u\cdot\na\r+\pa_i\r\pa_j\div u)=0,\\
\di \pa_j\pa_i u_t+u\cdot\na\pa_j\pa_i u+\frac{2\t}{3\r}\na\pa_j\pa_i\r+\f23\na\pa_j\pa_i\t
+\pa_j u\cdot\na\pa_i u+\pa_j\pa_i u\cdot\na u\\
\di~+\pa_i u\cdot\na\pa_j u+\pa_j\Big(\f{2\t}{3\r}\Big)\na\pa_i\r+\pa_j\pa_i\Big(\f{2\t}{3\r}\Big)\na\r+\pa_i\Big(\f{2\t}{3\r}\Big)\na\pa_j\r\\
\di
=\f{\mu(\t)}{\r}\big(\Delta\pa_j\pa_i u+\f13\na\pa_j\pa_i\div u\big)+\pa_j\Big(\f{\mu(\t)}{\r}\Big)\big(\Delta\pa_i u+\f13\na\pa_i\div u\big)\\
\di~
+\pa_j\Big\{\pa_i\Big(\f{\mu(\t)}{\r}\Big)\big(\Delta u+\f13\na\div u\big)\\
\di~+\pa_i\Big[\f{\mu'(\t)}{\r}\na\t\cdot\big(\na u+(\na u)^t-\f23\mathbb{I}\div u\big)\Big]
\Big\} -\pa_j\pa_i\Big(\f{1}{\r}\int v\otimes v\cdot\na_x\Pi dv\Big)\\
\di \pa_j\pa_i\t_t+u\cdot\na\pa_j\pa_i\t+\f23\t\pa_j\pa_i\div u+\pa_ju\cdot\na\pa_i\t+\pa_j\pa_i u\cdot\na\t
+\pa_i u\cdot\na\pa_j\t\\
\di~ +\f23\pa_j\t\pa_i\div u+\f23\pa_j\pa_i\t\div u+\f23\pa_i\t\pa_j\div u
=\f{\k(\t)}{\r}\Delta\pa_j\pa_i\t+\pa_j\Big(\f{\k(\t)}{\r}\Big)\Delta\pa_i\t \\
\di~ +\pa_j\Big\{\pa_i\Big(\f{\k(\t)}{\r}\Big)\Delta\t
+\pa_i\Big(\f{\k'(\t)}{\r}|\na\t|^2\Big)
+\pa_i\Big[\f{\mu(\t)}{\r}\Big(\f{(\na u+(\na u)^t)^2}{2}-\f23(\div u)^2\Big)\Big]
\Big\}\\
\di~+\pa_j\pa_i\Big[\f{1}{\r}\Big(-\int\f12|v|^2v\cdot\na_x\Pi dv+u\cdot\int v\otimes v\cdot\na_x\Pi dv
\Big)\Big].
\end{cases}
\ee
Multiplying the equation $\eqref{ns-2x-Pi}_1$ by $\f{2\t}{3\r^2}\pa_j\pa_i\r$, $\eqref{ns-2x-Pi}_2$ by $\pa_j\pa_i u$
and $\eqref{ns-2x-Pi}_3$ by $\f{\pa_j\pa_i\t}{\t}$,  and then adding them together and integrating over $[0,t]\times\mathbb{D}$, we have
\be\label{w3}
\begin{array}{l}
\di \|\na^2(\r,u,\t)(t)\|^2+\int_0^t\|\na^3(u,\t)\|^2d\tau\leq C\|\na^2(\r_0,u_0,\t_0)\|^2+C\v\\
\di\quad+C(\chi+\v)\int_0^t\|\na^2(\r,u,\t)\|^2d\tau+C\sum_{|\beta|=3}\int_0^t\int\int\f{\nu(v)}{\mb{M}_*}|\pa^{\beta}\mb{G}|^2dvdxd\tau\\
\di\quad +C(\chi+\v)\int_0^t\int\int\f{\nu(v)}{\mb{M}_*}\Big(\sum_{1\leq|\alpha|\leq 2}|\pa^{\alpha}\mb{G}|^2
+|\wt{\mb{G}}|^2\Big)dvdxd\tau,
\end{array}
\ee
where we have used the fact that
\begin{equation*}
\begin{array}{l}
\di\big|\int_0^t\int-\pa_j\pa_i\Big(\f{1}{\r}\int v\otimes v\cdot\na_x\Pi dv\Big)\pa_j\pa_i u\,dxd\tau\big|\\
\di= \big|\int_0^t\int\pa_i\Big(\f{1}{\r}\int v\otimes v\cdot\na_x\Pi dv\Big)\pa_j\pa_j\pa_i u\,dxd\tau\big|\\
\di \leq C\int_0^t\Big(\|\na\r\|_{L^{\infty}}\|\int v\otimes v\cdot\na_x\Pi dv\|
+ \|\int v\otimes v\cdot\na_x\pa_i \Pi dv\| \Big)\|\na^3 u\|d\tau\\
\di \leq \big(\f18+C(\chi+\v)\big)\int_0^t\|\na^3 u\|^2d\tau+C(\chi+\v)\int_0^t\|\int v\otimes v\cdot\na_x\Pi dv\|^2d\tau\\
\di\quad+C\int_0^t\|\int v\otimes v\cdot\na_x\pa_i \Pi dv\|d\tau\\
\di \leq \big(\f18+C(\chi+\v)\big)\int_0^t\|\na^3 u\|^2d\tau
+C\sum_{|\beta|=3}\int_0^t\int\int\f{\nu(v)}{\mb{M}_*}|\pa^{\beta}\mb{G}|^2dvdxd\tau+C\v\\
\di\quad +C(\chi+\v)\int_0^t\int\int\f{\nu(v)}{\mb{M}_*}\Big(\sum_{1\leq|\alpha|\leq 2}|\pa^{\alpha}\mb{G}|^2
+|\wt{\mb{G}}|^2\Big)dvdxd\tau
\end{array}
\end{equation*}\begin{equation*}
\begin{array}{l}\di\quad
+C(\chi+\v)\int_0^t\|\na^2(\r,u,\t)\|^2d\tau.
\end{array}
\end{equation*}
Next, we derive the estimate of $\di \int_0^t\|\na^3\r\|d\tau$. Applying $\pa_j~(j=1,2,3)$ to the system \eqref{ns-1x-G} gives
\be\label{ns-2x-G}
\begin{cases}
\di \pa_j\pa_i\r_t+u\cdot\na\pa_j\pa_i\r+\r\pa_j\pa_i\div u
+(\pa_j u\cdot\na\pa_i\r+\pa_i u\cdot\na\pa_j\r+\div u\pa_j\pa_i\r)\\
\quad+(\pa_j\r\pa_i\div u+\pa_j\pa_i u\cdot\na\r+\pa_i\r\pa_j\div u)=0,\\
\di \pa_j\pa_i u_t+u\cdot\na\pa_j\pa_i u+\frac{2\t}{3\r}\na\pa_j\pa_i\r+\f23\na\pa_j\pa_i\t
+\pa_j u\cdot\na\pa_i u+\pa_j\pa_i u\cdot\na u\\
\di\quad +\pa_i u\cdot\na\pa_j u+\pa_j\Big(\f{2\t}{3\r}\Big)\na\pa_i\r
+\pa_j\pa_i\Big(\f{2\t}{3\r}\Big)\na\r+\pa_i\Big(\f{2\t}{3\r}\Big)\na\pa_j\r\\
\di\quad
=-\pa_j\pa_i\Big(\f{1}{\r}\int v\otimes v\cdot\na_x\mb{G} dv\Big),\\
\di \pa_j\pa_i\t_t+u\cdot\na\pa_j\pa_i\t+\f23\t\pa_j\pa_i\div u+\pa_ju\cdot\na\pa_i\t+\pa_j\pa_i u\cdot\na\t
+\pa_i u\cdot\na\pa_j\t\\
\di\quad +\f23\pa_j\t\pa_i\div u+\f23\pa_j\pa_i\t\div u+\f23\pa_i\t\pa_j\div u\\[2mm]
\di\quad
=\pa_j\pa_i\Big[\f{1}{\r}\Big(-\int\f12|v|^2v\cdot\na_x\mb{G}dv+u\cdot\int v\otimes v\cdot\na_x\mb{G}dv
\Big)\Big].
\end{cases}
\ee
Multiplying the equation \eqref{ns-2x-G}$_2$ by $\na\pa_j\pa_i\r$ and then integrating over $[0,t]\times\mathbb{D}$ yield
\begin{equation}\label{phi-3}
\begin{array}{ll}
\di\quad \int\pa_j\pa_i u\cdot\na\pa_j\pa_i\r \,dx\Big|_0^t+\int_0^t\int\f{2\t}{3\r}|\na\pa_j\pa_i\r|^2dxd\tau\\
\di=\int_0^t\int\pa_j\pa_i\div u\big[u\cdot\na\pa_j\pa_i\r+\r\pa_j\pa_i\div u+(\pa_j u\cdot\na\pa_i\r+\pa_i u\cdot\na\pa_j\r+\div u\pa_j\pa_i\r)\\
\di\qquad\qquad\qquad\qquad +(\pa_j\r\pa_i\div u+\pa_j\pa_i u\cdot\na\r+\pa_i\r\pa_j\div u)\big]dxd\tau\\
\di
-\int_0^t\int(u\cdot\na\pa_j\pa_i u+\f23\na\pa_j\pa_i\t)\na\pa_j\pa_i\r dxd\tau\\[2mm]
\di-\int_0^t\int\Big[\pa_j u\cdot\na\pa_i u+\pa_j\pa_i u\cdot\na u+\pa_i u\cdot\na\pa_j u
+\pa_j\Big(\f{2\t}{3\r}\Big)\na\pa_i\r+\pa_j\pa_i\Big(\f{2\t}{3\r}\Big)\na\r\\
\di+\pa_i\Big(\f{2\t}{3\r}\Big)\na\pa_j\r\Big]\cdot\na\pa_j\pa_i\r\,dxd\tau
-\int_0^t\int\pa_j\pa_i\Big(\f{1}{\r}\int v\otimes v\cdot\na_x\mb{G} dv\Big)\cdot\na\pa_j\pa_i\r\,dxd\tau,
\end{array}
\end{equation}
where in the above equations we have used \eqref{ns-2x-G}$_1$ and the following fact
\be
\begin{array}{ll}
\di\int \pa_j\pa_i u_t\cdot\na\pa_j\pa_i\r \,dx=\f{d}{dt}\int \pa_j\pa_i u\cdot\na\pa_j\pa_i\r \,dx
-\int\pa_j\pa_i u\cdot\na\pa_j\pa_i\r_t\,dx\\
\di\quad
=\f{d}{dt}\int \pa_j\pa_i u\cdot\na\pa_j\pa_i\r \,dx+\int \pa_j\pa_i\div u\pa_j\pa_i\r_t \,dx\\
\di\quad=\f{d}{dt}\int \pa_j\pa_i u\cdot\na\pa_j\pa_i\r \,dx
-\int\pa_j\pa_i\div u\big[u\cdot\na\pa_j\pa_i\r+\r\pa_j\pa_i\div u\\
\di\qquad\qquad +(\pa_j u\cdot\na\pa_i\r+\pa_i u\cdot\na\pa_j\r+\div u\pa_j\pa_i\r)\\
\di
\qquad\qquad +(\pa_j\r\pa_i\div u+\pa_j\pa_i u\cdot\na\r+\pa_i\r\pa_j\div u)\big]dx.
\end{array}
\ee
By \eqref{phi-3} and  Cauchy's inequality, it holds that
\be\label{rho-3x}
\begin{array}{l}
\di \int_0^t\|\na^3\r\|^2d\tau \leq C\|(\na^2 u_0, \na^3\r_0)\|^2+C\|(\na^2 u,\na^3\r)(t)\|^2
\\
\di\qquad+C(\chi+\v)\int_0^t\|\na^2(\r,u,\t)\|^2d\tau
+C\int_0^t\int\int\f{\nu(v)}{\mb{M}_*}|\na_x^3\mb{G}|^2dvdxd\tau
\end{array}
\ee\be
\begin{array}{l}\di\qquad +C(\chi+\v)\int_0^t\int\int\f{\nu(v)}{\mb{M_*}}|\na^2_x\mb{G}|^2dvdxd\tau+C\int_0^t\|\na^3(u,\t)\|^2d\tau,\nonumber
\end{array}
\ee
which together with \eqref{w3} gives
\be\label{pa-3x}
\begin{array}{l}
\di \|\na^2(\r,u,\t)(t)\|^2+\int_0^t\|\na^3(\r,u,\t)\|^2d\tau
\leq C\|(\na^2\r_0,\na^2u_0,\na^2\t_0,\na^3\r_0)\|^2+C\v\\
\di\quad+C(\chi+\v)\int_0^t\|\na^2(\r,u,\t)\|^2d\tau
+C\sum_{|\beta|=3}\int_0^t\int\int\f{\nu(v)}{\mb{M}_*}|\pa^{\beta}\mb{G}|^2dvdxd\tau\\
\di\quad +C(\chi+\v)\int_0^t\int\int\f{\nu(v)}{\mb{M}_*}\Big(\sum_{1\leq|\alpha|\leq 2}|\pa^{\alpha}\mb{G}|^2
+|\wt{\mb{G}}|^2\Big)dvdxd\tau+C\|\na^3\r(t)\|^2.
\end{array}
\ee
Then we want to estimate $\na^2(\r_t,u_t,\t_t)$, $\na(\r_{tt}, u_{tt}, \t_{tt})$ and
$\pa^3_t(\r,u,\t)$. First, we estimate $\na^2(\r_t,u_t,\t_t)$. Multiplying the equation \eqref{ns-2x-G}$_1$ by $\pa_j\pa_i\r_t$,
\eqref{ns-2x-G}$_2$ by $\pa_j\pa_i u_t$ and \eqref{ns-2x-G}$_3$ by $\pa_j\pa_i\t_t$, respectively,
and then integrating over $[0,t]\times\mathbb{D}$, one has
\begin{equation}\label{xxt}
\begin{array}{l}
\di\int_0^t\|\na^2(\r_t, u_t,\t_t)\|^2 d\tau \leq C\int_0^t\|\na^3_x(\r,u,\t)\|^2d\tau
 +C(\chi+\v)\int_0^t\|\na^2(\r,u,\t)\|^2d\tau\\
\di\quad+C\int_0^t\int\int\frac{\nu(| v|)}{\mathbf{M}_*}|\na^3_x\mathbf{G}|^2 dv dxd\tau
+C(\chi+\v)\int_0^t\int\int\f{\nu(v)}{\mb{M}_*}|\na^2_x\mb{G}|^2dvdxd\tau.
\end{array}
\end{equation}

Next, applying $\pa_t$ to the system \eqref{ns-1x-G}, and multiplying
the three equations by $\pa_i\r_{tt}$,  $\pa_i u_{tt}$ and $\pa_i\t_{tt}$, respectively,
then integrating over $[0,t]\times\mathbb{D}$, we obtain
\begin{equation}\label{xtt}
\begin{array}{l}
\di\quad\int_0^t\|\na(\r_{tt}, u_{tt},\t_{tt})\|^2 d\tau\\
\di \leq C\int_0^t\|\na^2_x(\r_t,u_t,\t_t)\|^2d\tau
 +C(\chi+\v)\int_0^t\Big(\|\na^2(\r,u,\t)\|^2+\|\na(\r_t, u_t, \t_t)\|^2\Big)d\tau\\
\di\quad+C(\chi+\v)\sum_{|\alpha'|=1}\int_0^t\|\pa^{\alpha'}(\phi,\psi,\z)\|^2d\tau
+C\sum_{|\beta|=3}\int_0^t\int\int\f{\nu(v)}{\mb{M}_*}|\pa^{\beta}\mb{G}|^2dvdxd\tau\\
\di\quad +C(\chi+\v)\sum_{|\alpha|=2}\int_0^t\int\int\f{\nu(v)}{\mb{M}_*}|\pa^{\alpha}\mb{G}|^2 dvdxd\tau+C\v.
\end{array}
\end{equation}
Finally, applying $\pa_t$ to system \eqref{ns-2t-G}, and multiplying
the three equations by $\pa^3_t\r$,  $\pa^3_t u$ and $\pa^3_t\t$, respectively,
then integrating over $[0,t]\times\mathbb{D}$, it holds that
\begin{equation}\label{3t}
\begin{array}{l}
\di\quad\int_0^t\|\pa^3_t(\r, u,\t)\|^2 d\tau\leq C\int_0^t\|\na_x(\r_{tt},u_{tt},\t_{tt})\|^2d\tau\\
\di\quad  +C(\chi+\v)\int_0^t\Big(\|\na(\r_t,u_t,\t_t)\|^2+\|(\r_{tt}, u_{tt}, \t_{tt})\|^2\Big)d\tau\\
\di \quad
+C(\chi+\v)\sum_{|\alpha'|=1}\int_0^t\|\pa^{\alpha'}(\phi,\psi,\z)\|^2d\tau
+C\int_0^t\int\int\f{\nu(v)}{\mb{M}_*}|\na_x\mb{G}_{tt}|^2dvdxd\tau\\
\di\quad +C(\chi+\v)\int_0^t\int\int\f{\nu(v)}{\mb{M}_*}|\na_x\mb{G}_t|^2 dvdxd\tau+C\v.
\end{array}
\end{equation}

Combining \eqref{pa-3x}, \eqref{xxt}, \eqref{xtt} and \eqref{3t} together leads to \eqref{fluid-2x} directly.
\qed

$\bullet$ \underline{{\bf Step 4}}.  Estimate of $\di \sup_{t\in[0,T]} \int\int \f{|\pa^{\alpha}\mb{G}|^2}{\mb{M}_*}(t,x,v)dv dx$~($|\alpha|=2$):
\be\label{eG-2x}
\begin{array}{l}
\di \int\int \f{|\pa^{\alpha}\mb{G}|^2}{\mb{M}_*}(t,x,v)dv dx
+\int_0^t\int\int\f{\nu(v)}{\mb{M}_*}|\pa^{\alpha}\mb{G}|^2dvdxd\tau\\
\di
\leq C \int\int \f{|\pa^{\alpha}\mb{G}|^2}{\mb{M}_*}(0,x,v)dv dx+C\int_0^t\|\pa^{\alpha}\na(u,\t)\|^2d\tau+C\v\\
\di
+C(\chi+\v)\int_0^t\Big(\|\pa^{\alpha}(\r,u,\t)\|^2
+\sum_{|\alpha'|=1}\|\pa^{\alpha'}(\phi,\psi,\z)\|^2\Big)d\tau\\
\di +C\sum_{|\beta|=3}\int_0^t\int\int\f{\nu(v)}{\mb{M}_*}|\pa^{\beta}\mb{G}|^2dvdxd\tau\\
\di
+C(\chi+\v)\int_0^t\int\int\f{\nu(v)}{\mb{M}_*}|(\mb{G}_t,\na_x\mb{G},\wt{\mb{G}})|^2dvdxd\tau.
\end{array}
\ee

\vspace*{4pt}\noindent\emph{Proof of \eqref{eG-2x}.} Applying $\pa^{\alpha}$ to \eqref{non-fluid} for $|\alpha|=2$, one has
\begin{equation}\label{G-2x}
\ba
\pa^{\alpha}\mb{G}_t -\mb{L}_\mb{M}\pa^{\alpha}\mb{G}&=
-\pa^{\alpha}\mb{P}_1(v\cdot\nabla_x\mb{M})-\pa^{\alpha}\mb{P}_1(v\cdot\nabla_x\mb{G})
+2\sum_{|\alpha'|=1}Q(\pa^{\alpha'}\mb{M},\pa^{\alpha'}\mb{G})\\
&\quad +2\sum_{|\alpha'|=1}Q(\pa^{\alpha'}\mb{G},\pa^{\alpha'}\mb{G})
+2Q(\pa^{\alpha}\mb{M},\mb{G})+ 2Q(\mb{G}, \pa^{\alpha}\mb{G}).
\ea
\end{equation}
Multiplying the above equation by $\f{\pa^{\alpha}{\mb{G}}}{\mb{M_*}}$, and integrating over $[0,t]\times\mathbb{D}\times\mathbb{R}^3$ gives
\begin{equation}\label{es-G-2x}
\begin{array}{ll}
\di \int\int\f{|\pa^{\alpha}\mb{G}|^2}{2\mb{M}_*}dvdx\Big|_0^t
-\int_0^t\int\int\f{\pa^{\alpha}\mb{G}}{\mb{M}_*}\mb{L}_{\mb{M}}\pa^{\alpha}\mb{G}dvdxd\tau\\
\di=\int_0^t\int\int\Big[
-\pa^{\alpha}\mb{P}_1(v\cdot\nabla_x\mb{M})-\pa^{\alpha}\mb{P}_1(v\cdot\nabla_x\mb{G})
+2\sum_{|\alpha'|=1}Q(\pa^{\alpha'}\mb{M},\pa^{\alpha'}\mb{G})\\
\di\quad  +2\sum_{|\alpha'|=1}Q(\pa^{\alpha'}\mb{G},\pa^{\alpha'}\mb{G})+2Q(\pa^{\alpha}\mb{M},\mb{G})
+4Q(\mb{G}, \pa^{\alpha}\mb{G})\Big]\f{\pa^{\alpha}\mb{G}}{\mb{M}_*}dvdxd\tau\\
\di :=\sum_{i=16}^{21}I_i.
\end{array}
\end{equation}
Note $\mb{P}_1(v\cdot\nabla_x\mb{M})$ does not contain the density gradient $\na\r$,
thus by Cauchy's inequality, one has
\be\label{I16}
\ba
|I_{16}|& \leq C\int_0^t\int \Big(\sum_{|\alpha'|=1}\big(
|\pa^{\alpha'}\na(\r,u,\t)||\pa^{\alpha'}(\r,u,\t)|+|\na(\r,u,\t)||\pa^{\alpha}(\r,u,\t)|\\
& \quad +|\pa^{\alpha}\na(u,\t)|+|\na(\r,u,\t)||\pa^{\alpha'}(\r,u,\t)|^2\big)|\Big)
\Big(\int\f{\nu(v)}{\mb{M}_*}|\pa^{\alpha}\mb{G}|^2dv\Big)^{\f12}dxd\tau\\
& \leq \f{\tilde\s}{8}\int_0^t\int\int\f{\nu(v)}{\mb{M}_*}|\pa^{\alpha}\mb{G}|^2 dvdxd\tau
+C\int_0^t\|\pa^{\alpha}\na(u,\t)\|^2d\tau\\
&\quad+C(\chi+\v)\int_0^t\big(\|\pa^{\alpha}(\r,u,\t)\|^2
+\sum_{|\alpha'|=1}\|\pa^{\alpha'}(\phi,\psi,\z)\|^2\big)d\tau+C\v.
\ea
\ee
It follows from Cauchy's inequality and assumption \eqref{assum-1} that
\be\label{I17}
\ba
|I_{17}|& \leq \f{\tilde\s}{8}\int_0^t\int\int\f{\nu(v)}{\mb{M}_*}|\pa^{\alpha}\mb{G}|^2dvdxd\tau
+C\int_0^t\int\int\f{\nu(v)}{\mb{M}_*}|\pa^{\alpha}\na_x\mb{G}|^2dvdxd\tau\\
&\quad+C\int_0^t\int\big(|\pa^{\alpha}(\r,u,\t)|+|\pa^{\alpha'}(\r,u,\t)|^2\big)
\Big(\int\f{|\na_x\mb{G}|^2}{\mb{M}_*}dv\Big)^{\f12}\\
&\qquad\qquad\quad\cdot
\Big(\int\f{\nu(v)}{\mb{M}_*}|\pa^{\alpha'}\na_x\mb{G}|^2dv\Big)^{\f12}dxd\tau\\
& \leq \Big(\f{\tilde\s}{8}+C(\chi+\v)\Big)\int_0^t\int\int\f{\nu(v)}{\mb{M}_*}|\pa^{\alpha}\mb{G}|^2 dvdxd\tau\\
&\quad\di
+C\int_0^t\int\int\f{\nu(v)}{\mb{M}_*}|\pa^{\alpha}\na_x\mb{G}|^2dvdxd\tau\\
&\quad+C(\chi+\v)\int_0^t\Big(\|\pa^{\alpha}(\r,u,\t)\|^2+\int\int\f{\nu(v)}{\mb{M}_*}|\na_x\mb{G}|^2dvdx\Big)d\tau.
\ea
\ee
By Cauchy's inequality, Lemma \ref{Lemma 4.1} and \eqref{assum-1}, one has
{\small\be\label{I18}
\ba
|I_{18}|& \leq C\int_0^t\int\Big(\int\f{\nu(v)}{\mb{M}_*}|\pa^{\alpha}\mb{G}|^2dv\Big)^{\f12}
\sum_{|\alpha'|=1}\Big(\int\f{\nu(v)^{-1}}{\mb{M}_*}|Q(\pa^{\alpha'}\mb{M},\pa^{\alpha'}\mb{G})|^2dv\Big)^{\f12}dxd\tau\\
&\leq C\int_0^t\int\Big(\int\f{\nu(v)}{\mb{M}_*}|\pa^{\alpha}\mb{G}|^2dv\Big)^{\f12}\\
&
\qquad\qquad\sum_{|\alpha'|=1}\Big(\int\f{\nu(v)}{\mb{M}_*}|\pa^{\alpha'}\mb{M}|^2dv
\cdot \int\f{\nu(v)}{\mb{M}_*}|\pa^{\alpha'}\mb{G}|^2dv \Big)^{\f12}dxd\tau\\
&\leq C(\chi+\v)\int_0^t\int\int \f{\nu(v)}{\mb{M}_*}\Big(|\pa^{\alpha}\mb{G}|^2
+\sum_{|\alpha'|=1}|\pa^{\alpha'}\mb{G}|^2\Big)dvdxd\tau,
\ea
\ee}and
\be\label{I19}
\ba
|I_{19}|& \leq C\int_0^t\int\Big(\int\f{\nu(v)}{\mb{M}_*}|\pa^{\alpha}\mb{G}|^2dv\Big)^{\f12}
\sum_{|\alpha'|=1}\Big(\int\f{\nu(v)^{-1}}{\mb{M}_*}|Q(\pa^{\alpha'}\mb{G},\pa^{\alpha'}\mb{G})|^2dv\Big)^{\f12}dxd\tau\\
&\leq C\int_0^t\int\Big(\int\f{\nu(v)}{\mb{M}_*}|\pa^{\alpha}\mb{G}|^2dv\Big)^{\f12}
\\
&\qquad\qquad \sum_{|\alpha'|=1}\Big(\int\f{|\pa^{\alpha'}\mb{G}|^2}{\mb{M}_*}dv
\cdot \int\f{\nu(v)}{\mb{M}_*}|\pa^{\alpha'}\mb{G}|^2dv \Big)^{\f12}dxd\tau\\
&\leq C(\chi+\v)\int_0^t\int\int \f{\nu(v)}{\mb{M}_*}\Big(|\pa^{\alpha}\mb{G}|^2
+\sum_{|\alpha'|=1}|\pa^{\alpha'}\mb{G}|^2\Big)dvdxd\tau.
\ea
\ee
By Lemma \ref{Lemma 4.1} and the fact
\be\label{hard-1}
\begin{array}{ll}
\di\int|\pa^{\alpha}(\r,u,\t)|^2\Big(\int\f{\nu(v)}{\mb{M}_*}|\mb{G}|^2dv\Big)dx
\di \leq \|\pa^{\alpha}(\r,u,\t)\|^2\Big\|\int\f{\nu(v)}{\mb{M}_*}|\mb{G}|^2dv\Big\|_{L^{\infty}_x}\\
\di \leq C\|\pa^{\alpha}(\r,u,\t)\|^2\int\int\f{\nu(v)}{\mb{M}_*}(|\mb{G}|^2+|\na_x\mb{G}|^2+|\na_x^2\mb{G}|^2)dv dx\\
\di\leq C(\chi+\v)\|\pa^{\alpha}(\r,u,\t)\|^2+C(\chi+\v)\int\int\f{\nu(v)}{\mb{M}_*}|(\wt{\mb{G}},\na_x\mb{G},\na^2_x\mb{G})|^2dv dx,
\end{array}
\ee
it holds that
\be\label{I20}
\begin{array}{ll}
\di |I_{20}| \leq \f{\tilde\s}{8}\int_0^t\int\int\f{\nu(v)}{\mb{M}_*}|\pa^{\alpha}\mb{G}|^2dvdxd\tau
+C\int_0^t\int\int \f{\nu(v)^{-1}}{\mb{M}_*}|Q(\pa^{\alpha}\mb{M},\mb{G})|^2dvdxd\tau\\
\di\quad \leq \f{\tilde\s}{8}\int_0^t\int\int\f{\nu(v)}{\mb{M}_*}|\pa^{\alpha}\mb{G}|^2dvdxd\tau\\
\di\qquad
+C\int_0^t\int \Big(\int\f{\nu(v)}{\mb{M}_*}|\pa^{\alpha}\mb{M}|^2dv
\cdot\int\f{\nu(v)}{\mb{M}_*}|\pa^{\alpha}\mb{G}|^2dv\Big)dxd\tau\\
\di\quad \leq \f{\tilde\s}{8}\int_0^t\int\int\f{\nu(v)}{\mb{M}_*}|\pa^{\alpha}\mb{G}|^2dvdxd\tau
\\
\di\qquad+C\int_0^t\int\big(|\pa^{\alpha}(\r,u,\t)|^2+\sum_{|\alpha'|=1}|\pa^{\alpha'}(\r,u,\t)|^4\big)
\Big(\int\f{\nu(v)}{\mb{M}_*}|\mb{G}|^2dv\Big)dxd\tau\\
\quad\di\leq \f{\tilde\s}{8}\int_0^t\int\int\f{\nu(v)}{\mb{M}_*}|\pa^{\alpha}\mb{G}|^2dvdxd\tau
+C(\chi+\v)\sum_{1\leq|\alpha|\leq2}\int_0^t\|\pa^{\alpha}(\phi,\psi,\z)\|^2d\tau\\
\di\qquad
+C(\chi+\v)\int_0^t\int\int \f{\nu(v)}{\mb{M}_*}|(\na_x^2\mb{G},\na_x\mb{G},\wt{\mb{G}})|^2dvdxd\tau+C\v.
\end{array}
\ee
Similarly,  by \eqref{assum-1}, Cauchy's inequality and Sobolev's inequality, we have
\be\label{hard}
\begin{array}{l}
\di\int\Big(\int\f{\nu(v)}{\mb{M}_*}|\mb{G}|^2dv\cdot\int\f{|\pa^{\alpha}\mb{G}|^2}{\mb{M}_*}dv\Big)dx
\leq \Big\|\int\f{\nu(v)}{\mb{M}_*}|\mb{G}|^2dv\Big\|_{L^{\infty}_x}\int\int\f{|\pa^{\alpha}\mb{G}|^2}{\mb{M}_*} dv dx\\
\di\leq C\int\int\f{\nu(v)}{\mb{M}_*}(|\mb{G}|^2+|\na_x\mb{G}|^2+|\na_x^2\mb{G}|^2) dvdx\cdot\int\int\f{|\pa^{\alpha}\mb{G}|^2}{\mb{M}_*} dv dx\\
\di\leq C(\chi+\v)\int\int\f{\nu(v)}{\mb{M}_*}|(\wt{\mb{G}},\na_x\mb{G},\pa^{\alpha}\mb{G})|^2dv dx.
\end{array}
\ee
Then it follows from Lemma \ref{Lemma 4.1} and \eqref{hard} that
\be\label{I21}
\ba
|I_{21}|& \leq \f{\tilde\s}{8}\int_0^t\int\int\int\f{\nu(v)}{\mb{M}_*}|\pa^{\alpha}\mb{G}|^2dvdxd\tau
+C\int_0^t\int\int\f{\nu(v)^{-1}}{\mb{M}_*}|Q(\mb{G},\pa^{\alpha}\mb{G})|^2dvdxd\tau\\
&\leq \f{\tilde\s}{8}\int_0^t\int\int\int\f{\nu(v)}{\mb{M}_*}|\pa^{\alpha}\mb{G}|^2dvdxd\tau
+C\int_0^t\int\Big(\int\f{|\mb{G}|^2}{\mb{M}_*}dv\cdot\int\f{\nu(v)}{\mb{M}_*}|\pa^{\alpha}\mb{G}|^2dv\\
&\quad +\int\f{\nu(v)}{\mb{M}_*}|\mb{G}|^2dv\cdot\int\f{|\pa^{\alpha}\mb{G}|^2}{\mb{M}_*}dv\Big)dxd\tau\\
&\leq \Big(\f{\tilde\s}{8}+C(\chi+\v)\Big)\int_0^t\int\int \f{\nu(v)}{\mb{M}_*}|\pa^{\alpha}\mb{G}|^2dvdxd\tau\\
&\quad\di
+C(\chi+\v)\int_0^t\int\int \f{\nu(v)}{\mb{M}_*}|(\wt{\mb{G}},\na_x\mb{G})|^2dvdxd\tau.
\ea
\ee

Substituting \eqref{I16}-\eqref{I19}, \eqref{I20} and \eqref{I21} into \eqref{es-G-2x} gives \eqref{eG-2x}.
\hfill $\Box$

Combining \eqref{fluid-2x} and \eqref{eG-2x} together, we can derive that
\be\label{order-2}
\begin{array}{l}
\di \|\na^2(\r,u,\t)(t)\|^2
+\sum_{|\alpha|=2}\int\int \f{|\pa^{\alpha}\mb{G}|^2}{\mb{M}_*}(t,x,v)dv dx
+\sum_{|\beta|=3}\int_0^t\|\pa^{\beta}(\r, u,\t)\|^2d\tau\\
\di\quad+\sum_{|\alpha|=2}\int_0^t\int\int\f{\nu(v)}{\mb{M}_*}|\pa^{\alpha}\mb{G}|^2dvdxd\tau
 \leq C\|\na^2(\r_0,u_0,\t_0),\na^3\r_0\|^2
\end{array}
\ee\be
\begin{array}{l} \di\quad+C \sum_{|\alpha|=2}\int\int \f{|\pa^{\alpha}\mb{G}|^2}{\mb{M}_*}(0,x,v)dv dx
+C\sum_{|\beta|=3}\int_0^t\int\int\f{\nu(v)}{\mb{M_*}}|\pa^{\beta}\mb{G}|^2dvdxd\tau\nonumber\\
\di\quad+C(\chi+\v)\int_0^t\Big(\sum_{|\alpha|=2}\|\pa^{\alpha}(\r,u,\t)\|^2
+\sum_{|\alpha'|=1}\|\pa^{\alpha'}(\phi,\psi,\z)\|^2\Big)d\tau
\nonumber\\
\di\quad+C(\chi+\v)\int_0^t\int\int\f{\nu(v)}{\mb{M_*}}\Big(\sum_{|\alpha'|=1}|\pa^{\alpha'}\mb{G}|^2+|\wt{\mb{G}}|^2\Big)dvdxd\tau
+C\v+C\|\na^3\r(t)\|^2.\nonumber
\end{array}
\ee

Finally, in order to obtain the third order derivatives (with respect to $x$ and/or $t$) estimates on $\mb{G}$, and
the third order spatial derivative estimates on the density $\na^3\r$, we need to work on the original Boltzmann equation \eqref{B}.

$\bullet$ \underline{\bf Step 5}.  Estimate of $\di \sup_{t\in[0,T]} \int\int\f{|\pa^{\beta}f|^2}{\mb{M}_*}(t,x,v)dv dx$~($|\beta|=3$):
\be\label{es-f}
\begin{array}{l}
\di\quad \int\int\f{|\pa^{\beta}f|^2}{\mb{M}_*}(t,x,v)dv dx+\int_0^t\int\int\f{\nu(v)}{\mb{M}_*}|\pa^{\beta}\mb{G}|^2dvdxd\tau\\
\di \leq C\int\int\f{|\pa^{\beta}f|^2}{\mb{M}_*}(0,x,v)dv dx
+C(\chi+\v+\eta_0)\sum_{1\leq|\gamma|\leq3}\int_0^t\|\pa^{\gamma}(\phi,\psi,\z)\|^2d\tau\\
\di\quad+C(\chi+\v)\sum_{1\leq|\alpha|\leq 2}\int_0^t\int\int\f{\nu(v)}{\mb{M}_*}|\pa^{\alpha}\mb{G}|^2dvdxd\tau+C\v.
\end{array}
\ee
\vspace*{4pt}\noindent\emph{Proof of \eqref{es-f}.} Applying $\pa^{\beta}$ to equation \eqref{B} gives
\be
\begin{array}{l}
\di \pa^{\beta}f_t+v\cdot\na_x\pa^{\beta}f=\mb{L}_{\mb{M}}\pa^{\beta}\mb{G}\\
\di \quad
+2\sum_{0\leq|\alpha|\leq|\beta|-1}C_{\beta}^{\alpha}Q(\pa^{\beta-\alpha}\mb{M},\pa^{\alpha}\mb{G})
+\sum_{0\leq|\alpha|\leq|\beta|}C_{\beta}^{\alpha}Q(\pa^{\beta-\alpha}\mb{G},\pa^{\alpha}\mb{G}).
\end{array}
\ee
Multiplying the above equation by $\f{\pa^{\beta}f}{\mb{M_*}}$ and integrating over $[0,t]\times\mathbb{D}\times\mathbb{R}^3$,
we have
\be\label{3f}
\begin{array}{l}
\di\int\int\f{|\pa^{\beta}f|^2}{2\mb{M}_*}dv dx\Big|_0^t
-\int_0^t\int\int\f{\pa^{\beta}\mb{G}}{\mb{M}_*}\mb{L}_{\mb{M}}\pa^{\beta}\mb{G}dvdxd\tau\\
\di
=\int_0^t\int\int\f{\pa^{\beta}\mb{M}}{\mb{M}_*}\mb{L}_{\mb{M}}\pa^{\beta}\mb{G}dvdxd\tau+\int_0^t\int\int\Big[2\sum_{0\leq|\alpha|\leq|\beta|-1}C_{\beta}^{\alpha}Q(\pa^{\beta-\alpha}\mb{M},\pa^{\alpha}\mb{G})\\
\di
\qquad+\sum_{0\leq|\alpha|\leq|\beta|}C_{\beta}^{\alpha}Q(\pa^{\beta-\alpha}\mb{G},\pa^{\alpha}\mb{G})
\Big]\f{\pa^{\beta}f}{\mb{M}_*}dvdxd\tau.
\end{array}
\ee
Since $\mb{M}_t$, $\na_x\mb{M}\in \mathfrak{N}$, $\mb{P}_1(\pa^{\beta}\mb{M})$ does not contain $\pa^{\beta}(\r,u,\t)$. Thus we have
\be\label{f-1}
\begin{array}{l}
\di \Big|\int_0^t\int\int \f{\mb{L}_{\mb{M}}\pa^{\beta}\mb{G}}{\mb{M}}\cdot\pa^{\beta}\mb{M} \,dvdxd\tau\Big|
= \Big|\int_0^t\int\int \f{\mb{L}_{\mb{M}}\pa^{\beta}\mb{G}}{\mb{M}}\cdot\mb{P}_1(\pa^{\beta}\mb{M}) \,dvdxd\tau\Big|\\
\di\leq \f{\tilde\s}{8}\int_0^t\int\int\f{\nu(v)}{\mb{M}_*}|\pa^{\beta}\mb{G}|^2dvdxd\tau
+C(\chi+\v)\sum_{1\leq|\alpha|\leq2}\int_0^t\|\pa^{\alpha}(\phi,\psi,\z)\|^2d\tau+C\v
\end{array}
\ee
and
\be
\begin{array}{l}
\di \Big|\int_0^t\int\int\mb{L}_{\mb{M}}\pa^{\beta}\mb{G}\cdot\pa^{\beta}\mb{M}
\Big(\f{1}{\mb{M}}-\f{1}{\mb{M}_*}\Big) dvdxd\tau\Big|\nonumber
\end{array}
\ee\be
\begin{array}{l}\di\leq \f{\tilde\s}{8}\int_0^t\int\int\f{\nu(v)}{\mb{M}_*}|\pa^{\beta}\mb{G}|^2dvdxd\tau
+C\eta_0\int_0^t\|\pa^{\beta}(\r,u,\t)\|^2d\tau\\
\di\quad+C(\chi+\v)\sum_{1\leq|\alpha|\leq2}\int_0^t\|\pa^{\alpha}(\phi,\psi,\z)\|^2d\tau+C\v,
\end{array}
\ee
which together with \eqref{f-1} implies that
\be\label{f-2}
\begin{array}{l}
\di \Big|\int_0^t\int\int\f{\pa^{\beta}\mb{M}}{\mb{M}_*}\mb{L}_{\mb{M}}\pa^{\beta}\mb{G}dvdxd\tau\Big|\\
\di= \Big|\int_0^t\int\int \f{\mb{L}_{\mb{M}}\pa^{\beta}\mb{G}}{\mb{M}}\cdot\mb{P}_1(\pa^{\beta}\mb{M}) \,dvdxd\tau\Big|\\
\di\quad
+\Big|\int_0^t\int\int\mb{L}_{\mb{M}}\pa^{\beta}\mb{G}\cdot\pa^{\beta}\mb{M}
\Big(\f{1}{\mb{M}}-\f{1}{\mb{M}_*}\Big) dvdxd\tau\Big|\\
\di \leq \f{\tilde\s}{4}\int_0^t\int\int\f{\nu(v)}{\mb{M}_*}|\pa^{\beta}\mb{G}|^2dvdxd\tau
+C\eta_0\int_0^t\|\pa^{\beta}(\r,u,\t)\|^2d\tau\\
\di\quad+C(\chi+\v)\sum_{1\leq|\alpha|\leq2}\int_0^t\|\pa^{\alpha}(\phi,\psi,\z)\|^2d\tau+C\v.
\end{array}
\ee
Using the similar technique as in obtaining \eqref{hard-1} and \eqref{I20}, we have
\be\label{f-3}
\begin{array}{l}
\di\Big|\int_0^t\int\int\sum_{0\leq|\alpha|\leq |\beta|-1}C_{\beta}^{\alpha}
Q(\pa^{\beta-\alpha}\mb{M},\pa^{\alpha}\mb{G})\f{\pa^{\beta}\mb{M}}{\mb{M}_*} dvdxd\tau\Big|\\
\di \leq C(\chi+\v)\int_0^t\Big(\|\pa^{\beta}(\r,u,\t)\|^2+\sum_{1\leq|\alpha|\leq 2}\|\pa^{\alpha}(\phi,\psi,\z)\|^2\Big)d\tau\\
\di\quad+C(\chi+\v)\int_0^t\int\int\f{\nu(v)}{\mb{M}_*}\big(\sum_{1\leq|\alpha|\leq 3}|\pa^{\alpha}\mb{G}|^2+|\wt{\mb{G}}|\big)dvdxd\tau+C\v.
\end{array}
\ee
Similar to \eqref{hard}, by Sobolev inequality one has,
\be\label{hard-2}
\begin{array}{l}
\di\int\Big(\int\f{|\pa^{\beta}\mb{M}|^2}{\mb{M}_*}dv\cdot\int\f{\nu(v)}{\mb{M}_*}|\mb{G}|^2dv\Big)dx
\leq \Big\|\int\f{\nu(v)}{\mb{M}_*}|\mb{G}|^2dv\Big\|_{L^{\infty}_x}\int\int\f{|\pa^{\beta}\mb{M}|^2}{\mb{M}_*} dv dx\\
\di\leq C\int\int\f{\nu(v)}{\mb{M}_*}|(\mb{G},\na_x\mb{G},\na_x^2\mb{G})|^2dvdx\cdot\int\int\f{|\pa^{\beta}\mb{M}|^2}{\mb{M}_*} dv dx\\
\di\leq C(\chi+\v)\int\int\f{\nu(v)}{\mb{M}_*}|(\wt{\mb{G}},\na_x\mb{G},\na^2_x\mb{G})|^2 dv dx
+C(\chi+\v)\int_0^t\|\pa^{\beta}(\r,u,\t)\|^2d\tau\\
\di\quad+C(\chi+\v)\sum_{1\leq|\alpha|\leq2}\|\pa^{\alpha}(\phi,\psi,\z)\|^2d\tau+C\v,
\end{array}
\ee
which together with Lemma \ref{Lemma 4.1} and \eqref{assum-1} leads to
\be
\begin{array}{l}
\di\Big|\int_0^t\int\int\sum_{0\leq|\alpha|\leq |\beta|-1}C_{\beta}^{\alpha}
Q(\pa^{\beta-\alpha}\mb{M},\pa^{\alpha}\mb{G})\f{\pa^{\beta}\mb{G}}{\mb{M}_*} dvdxd\tau\Big|\nonumber\\
\di \leq \f{\tilde\s}{8}\int_0^t\int\int\f{\nu(v)}{\mb{M}_*}|\pa^{\beta}\mb{G}|^2dvdxd\tau
+C\int_0^t\int\int \f{\nu(v)^{-1}}{\mb{M}_*}
|Q(\pa^{\beta-\alpha}\mb{M},\pa^{\alpha}\mb{G})|^2dvdxd\tau\nonumber\\
\di\leq \f{\tilde\s}{8}\int_0^t\int\int\f{\nu(v)}{\mb{M}_*}|\pa^{\beta}\mb{G}|^2dvdxd\tau
+C\int_0^t\int \Big(\int\f{\nu(v)}{\mb{M}_*}|\pa^{\beta-\alpha}\mb{M}|^2dv\cdot\int\f{|\pa^{\alpha}\mb{G}|^2}{\mb{M}_*}dv\nonumber\\
\di\quad +\int\f{|\pa^{\beta-\alpha}\mb{M}|^2}{\mb{M}_*}dv\cdot\int\f{\nu(v)}{\mb{M}_*}|\pa^{\alpha}\mb{G}|^2dv\Big)dx d\tau\nonumber
\end{array}
\ee\be\label{f-4}
\begin{array}{l}\di\leq \f{\tilde\s}{8}\int_0^t\int\int\f{\nu(v)}{\mb{M}_*}|\pa^{\beta}\mb{G}|^2dvdxd\tau+C(\chi+\v)\int_0^t\sum_{1\leq|\alpha|\leq 2}\|\pa^{\alpha}(\phi,\psi,\z)\|^2d\tau+C\v
\\
\di\quad +C(\chi+\v)\int_0^t\Big[\int\int\f{\nu(v)}{\mb{M}_*}(\sum_{1\leq|\gamma|\leq 3}|\pa^{\gamma}\mb{G}|^2+|\wt{\mb{G}}|^2)dvdx+\|\pa^{\beta}(\r,u,\t)\|^2\Big]d\tau.
\end{array}
\ee
Similar to \eqref{hard-2} and \eqref{f-4}, we have
\be\label{f-5}
\begin{array}{l}
\di\Big|\int_0^t\int\int\sum_{0\leq|\alpha|\leq |\beta|}C_{\beta}^{\alpha}
Q(\pa^{\beta-\alpha}\mb{G},\pa^{\alpha}\mb{G})\f{\pa^{\beta}\mb{G}}{\mb{M}_*} dvdxd\tau\Big|\\
\di\leq \f{\tilde\s}{8}\int_0^t\int\int\f{\nu(v)}{\mb{M}_*}|\pa^{\beta}\mb{G}|^2dvdxd\tau+C(\chi+\v)\int_0^t\sum_{1\leq|\alpha|\leq 2}\|\pa^{\alpha}(\phi,\psi,\z)\|^2d\tau+C\v
\\
\di\quad +C(\chi+\v)\int_0^t\Big[\int\int\f{\nu(v)}{\mb{M}_*}(\sum_{1\leq|\gamma|\leq 3}|\pa^{\gamma}\mb{G}|^2+|\wt{\mb{G}}|^2)dvdx+\|\pa^{\beta}(\r,u,\t)\|^2\Big]d\tau.\end{array}
\ee

Substituting \eqref{f-2}, \eqref{f-3}, \eqref{f-4} and \eqref{f-5} into \eqref{3f} implies \eqref{es-f}. \hfill $\Box$

We can deduce from \eqref{order-2} and \eqref{es-f} that
\be\label{2nd+3rd}
\begin{array}{l}
\di \|\na^2(\r,u,\t)(t)\|^2+\int\int\Big(\sum_{|\beta|=3}\f{|\pa^{\beta}f|^2}{\mb{M}_*}
+\sum_{|\alpha|=2}\f{|\pa^{\alpha}\mb{G}|^2}{\mb{M}_*}\Big)(t,x,v)dv dx\\
\di
+\sum_{|\beta|=3}\int_0^t\|\pa^{\beta}(\r, u,\t)\|^2d\tau+\sum_{2\leq|\beta'|\leq3}\int_0^t\int\int\f{\nu(v)}{\mb{M}_*}|\pa^{\beta'}\mb{G}|^2dvdxd\tau\\
\di
 \leq C\|(\na^2\r_0,\na^2 u_0,\na^2\t_0,\na^3\r_0)\|^2+C\v\\
 \di +C\int\int\Big(\sum_{|\beta|=3}\f{|\pa^{\beta}f|^2}{\mb{M}_*}
+\sum_{|\alpha|=2}\f{|\pa^{\alpha}\mb{G}|^2}{\mb{M}_*}\Big)(0,x,v)dv dx\\
\di
+C(\chi+\v+\eta_0)\sum_{1\leq|\alpha|\leq2}\int_0^t\|\pa^{\alpha}(\phi,\psi,\z)\|^2d\tau\\
\di+C(\chi+\v)\int_0^t\int\int\f{\nu(v)}{\mb{M_*}}\Big(\sum_{|\alpha'|=1}|\pa^{\alpha'}\mb{G}|^2+|\wt{\mb{G}}|^2\Big)dvdxd\tau.
\end{array}
\ee
\vspace*{4pt}\noindent\emph{Proof of Proposition \ref{Prop3.2}.} Multiplying \eqref{2nd+3rd} by a large constant $C$, then adding \eqref{order-1} together implies \eqref{higher}, thus we can finish the proof of Proposition \ref{Prop3.2}.
\qed

\section{Appendix: Local-in-time existence}\label{appendix}
\setcounter{equation}{0}

In this appendix, we prove Lemma \ref{local} for the local-in-time existence of the solution to the 3D Boltzmann equation \eqref{B}-\eqref{far-field}.

\vspace*{4pt}\noindent\emph{Proof of Lemma \ref{local}.} For this, we first define a functional space:
\be
\di \mathcal{S}([0,T])=\Big\{ F(t,x,v)\in H^3_x\Big(L^2_v\Big(\frac{1}{\sqrt{\mb{M}_*}}\Big)\Big):
\|F\|_{\mathcal{S}}\leq \Xi \Big\},
\ee
where $\Xi>0$, $T>0$ are time- and the $\|\cdot\|_{\mathcal{S}}$ is defined as:
\be
\begin{array}{l}
\di\|F\|^2_{\mathcal{S}}:=\sum_{0\leq|\beta|\leq3}\Big(\sup_{0\leq t\leq T}\int\int\f{|\pa^{\beta}F(t,x,v)|^2}{\mb{M}_*}dv dx\\
\di \qquad\qquad\qquad\qquad+\int_0^T\int\int\f{1+|v|}{\mb{M}_*}|\pa^{\beta}F(t,x,v)|^2dvdxd\tau\Big).
\end{array}
\ee
where $\mb{M}_*$ is a global Maxwellian $\mb{M_*}=\mb{M}_{[\r_*,u_*,\t_*]}$ with $\r_*>0,~\t_*>0$ as in Theorem \ref{thm}.

Set $g(t,x,v)=f(t,x,v)-\bar{\mb{M}}(t,x ,v)$, where $\bar{\mb{M}}(t,x_1,v)=\mb{M}_{[\bar\r,\bar u,\bar\t]}(t,x_1,v)$
with $(\bar\r,\bar u,\bar\t)$ being the $3$-rarefaction wave defined in \eqref{au}, then $g(t,x,v)$ solves
\be
\begin{cases}
g_t+v\cdot\na_x g=\mb{L}_{\bar{\mb{M}}}g+Q(g,g)-(\bar{\mb{M}}_t+v\cdot\na_x\bar{\mb{M}}),\\
g(t,x,v)|_{t=0}=g_0(x,v),
\end{cases}
\ee
where the linearized collision operator $\mb{L}_{\bar{\mb{M}}}$ is defined in \eqref{lm} with $\mb{M}$ replaced by $\bar{\mb{M}}$ as follows,
\be\label{bar-lm}
(\mb{L}_{\bar{\mb{M}}}g)(v)=-\nu_{\bar{\mb{M}}}(v)g(v)+\sqrt{\bar{\mb{M}}(v)}K_{\bar{\mb{M}}}
\left(\left(\frac{g}{\sqrt{\bar{\mb{M}}}}\right)(v)\right).
\ee
Here the collision frequency $\nu_{\bar{\mb{M}}}(v)\sim(1+|v|)$ as $|v|\rightarrow+\infty$ and $K_{\bar{\mb{M}}}(\cdot)=K_{2\bar{\mb{M}}}(\cdot)-K_{1\bar{\mb{M}}}(\cdot)$ with the kernel $k_{i\bar{\mb{M}}}(v,v_*)$ $(i=1,2)$
are all defined in \eqref{nkk}. Then we have the following lemma.
\begin{lemma}\label{k}
If $\bar\t/2<\t_ *<\bar\t$, then there exist positive constants $\bar\sigma=\bar\sigma(\bar\r,\bar u,\bar \t; \r_ *,\break u_ *,$ $\t_ *)$
and $\bar\eta_0=\bar\eta_0(\bar\r,\bar u,\bar \t;\r_ *,u_ *,\t_ *)$ such that if $|\bar\r-\r_ *|+|\bar u-u_ *|+|\bar\t-\t_ *|<\bar\eta_0$,
then for $g(v)\in  \mathfrak{N}^\bot$,
$$
-\int\f{g\mb{L}_{\bar{\mb{M}}}g}{\mb{M}_ *}dv\geq
\bar\sigma\int\f{\nu_{\bar{\mb{M}}}(v)g^2}{\mb{M}_ *}dv.
$$
\end{lemma}
The proof of this lemma can be followed step by step as Lemma 4.2 in \cite{Liu-Yang-Yu-Zhao}, we omit it for brevity.

Before we prove the local existence of 3D Boltzmann equation \eqref{B}-\eqref{far-field}, we should note that the solution $f(t,x,v)$
is non-negative for any $t>0$ if the initial data $f_0(x,v)$ is non-negative, whose proof can be found in \cite{LY}.

Under the above preparations, we now turn to construct local-in-time solution to the 3D Boltzmann equation \eqref{B}-\eqref{far-field}.
For this purpose, we consider the following iteration sequence $g^n(t,x,v)$ $(n\geq0)$
\be\label{local-g}
\begin{cases}
g^{n+1}_t+v\cdot\na_x g^{n+1}=\mb{L}_{\bar{\mb{M}}}g^{n+1}+Q(g^n,g^n)-(\bar{\mb{M}}_t+v\cdot\na_x\bar{\mb{M}}),\\
g^{n+1}(t,x,v)|_{t=0}=g_0(x,v).
\end{cases}
\ee

We want to show by induction that if $\|g_0\|_{\mathcal{S}}\leq \f{\Xi}{2\sqrt{C_0}}$ with $C_0:=\f{1}{\min\{1,\bar\sigma\}}\geq1$,
then $\|g^n\|_{\mathcal{S}}\leq \Xi$ for all $n$, provided that $\Xi$ and $T$ are chosen suitably small.  In fact, multiplying \eqref{local-g} by $\f{g^{n+1}}{\mb{M}_*}$ and then integrating over $[0,t]\times\mathbb{D}\times\mathbb{R}^3$ give
\be\label{gn}
\begin{array}{l}
\di\quad \int\int \f{|g^{n+1}|^2}{2\mb{M}_*}dv dx\Big|^{\tau=t}_{\tau=0}-\int_0^t\int\int \f{g^{n+1}}{\mb{M}_*}\mb{L}_{\bar{\mb{M}}}g^{n+1}dvdxd\tau\\
\di =\int_0^t\int\int \f{g^{n+1}}{\mb{M}_*}Q(g^n,g^n)dvdxd\tau-\int_0^t\int\int\f{g^{n+1}}{\mb{M}_*}(\bar{\mb{M}}_t+v\cdot\na_x\bar{\mb{M}})dvdxd\tau\\
\di
:=\sum_{i=1}^2J_i.
\end{array}
\ee
By Sobolev inequality, we have
\be\label{a-1}
\begin{array}{l}
\di\Big\|\int\f{|g^n|^2}{\mb{M}_*}dv\Big\|_{L^{\infty}_x}
\leq C\int\int\f{|g^n|^2+|\na_xg^n|^2+|\na_x^2g^n|^2}{\mb{M}_*}dvdx
\leq C\Xi^2.
\end{array}
\ee
It follows from Cauchy's inequality, Lemma \ref{Lemma 4.1} and \eqref{a-1} that
\be\label{J1}
\ba
J_1 & \leq \f{\bar\sigma}{4}\int_0^t\int\int\f{\nu_{\bar{\mb{M}}}(v)}{\mb{M}_*}|g^{n+1}|^2dv dxd\tau\\
&\quad
+C\int_0^t\int\int\f{\nu_{\bar{\mb{M}}}(v)^{-1}}{\mb{M_*}}|Q(g^n,g^n)|^2dvdxd\tau\\
& \leq \f{\bar\sigma}{4}\int_0^t\int\int\f{\nu_{\bar{\mb{M}}}(v)}{\mb{M}_*}|g^{n+1}|^2dv dxd\tau\\
&
\quad+C\int_0^t\int \Big(\int\f{|g^n|^2}{\mb{M}_*}dv\cdot\int\f{\nu_{\bar{\mb{M}}}(v)}{\mb{M}_*}|g^n|^2dv\Big)dxd\tau\\
& \leq \f{\bar\sigma}{4}\int_0^t\int\int\f{\nu_{\bar{\mb{M}}}(v)}{\mb{M}_*}|g^{n+1}|^2dv dxd\tau\\
&
\quad+C\sup_{0\leq t\leq T}\Big\|\int\f{|g^n|^2}{\mb{M}_*}dv\Big\|_{L^{\infty}_x}\int_0^t\int \int\f{\nu_{\bar{\mb{M}}}(v)}{\mb{M}_*}|g^n|^2dvdxd\tau\\
& \leq \f{\bar\sigma}{4}\int_0^t\int\int\f{\nu_{\bar{\mb{M}}}(v)}{\mb{M}_*}|g^{n+1}|^2dv dxd\tau+C\Xi^4.
\ea
\ee
By Cauchy's inequality and Lemma \ref{appu}, it holds
\be\label{J2}
\ba
J_2 & \leq \f{\bar\sigma}{4}\int_0^t\int\int\f{\nu_{\bar{\mb{M}}}(v)}{\mb{M}_*}|g^{n+1}|^2dv dxd\tau
+C\int_0^t\|(\bar\r_{x_1},\bar u_{1x_1},\bar \t_{x_1})\|^2d\tau\\[3mm]
& \leq \f{\bar\sigma}{4}\int_0^t\int\int\f{\nu_{\bar{\mb{M}}}(v)}{\mb{M}_*}|g^{n+1}|^2dv dxd\tau+C\ln (1+t).
\ea
\ee
Substituting \eqref{J1}-\eqref{J2} into \eqref{gn}, which along with Lemma \ref{k} leads to
\be\label{g0}
\begin{array}{l}
 \di\quad \int\int \f{|g^{n+1}(t,x,v)|^2}{\mb{M_*}}dv dx +\int_0^t\int\int \f{\nu_{\bar{\mb{M}}}(v)}{\mb{M}_*}|g^{n+1}|^2dvdxd\tau\\[3mm]
\di \leq C_0\int\int\f{|g_0(x,v)|^2}{\mb{M}_*}dv dx+C(\Xi^4+\ln(1+t)).
\end{array}
\ee

For the corresponding estimate for $\pa^{\beta}g^{n+1}(t,x,v)$ with $|\beta|=j~(j=1,2,3)$, since
\be
\begin{array}{l}
\di \pa^{\beta}(\mb{L}_{\bar{\mb{M}}}g^{n+1})=\mb{L}_{\bar{\mb{M}}}(\pa^{\beta}g^{n+1})
+2\sum_{0<|\alpha|\leq j}C_{\beta}^{\alpha}Q(\pa^{\alpha}\bar{\mb{M}},\pa^{\beta-\alpha}g^{n+1}).
\end{array}
\ee

Applying $\pa^{\beta}$ to \eqref{local-g}, which together with the above fact leads to
\be\label{g-x}
\begin{cases}
\di(\pa^{\beta}g^{n+1})_t+v\cdot\na_x (\pa^{\beta}g^{n+1})
=\mb{L}_{\bar{\mb{M}}}(\pa^{\beta}g^{n+1})\\[3mm]
\di \qquad +2\sum_{0<|\alpha|\leq j}C_{\beta}^{\alpha}Q(\pa^{\alpha}\bar{\mb{M}},\pa^{\beta-\alpha}g^{n+1}) +2Q(\pa^{\beta}g^n,g^n)\\[3mm]
\di\qquad
+2\sum_{0<|\alpha|<j}C_{\beta}^{\alpha}Q(\pa^{\alpha}g^n,\pa^{\beta-\alpha}g^n)-(\pa^{\beta}\bar{\mb{M}}_t+v\cdot\na_x\pa^{\beta}\bar{\mb{M}}),\\[3mm]
\di \pa^{\beta}g^{n+1}(t,x,v)|_{t=0}=\pa^{\beta}g_0(x,v).
\end{cases}
\ee
Multiplying the above equation by $\f{\pa^{\beta}g^{n+1}}{\mb{M}_*}$, then integrating over $[0,t]\times\mathbb{D}\times\mathbb{R}^3$, one leads to
{\small\be\label{gn+}
\begin{array}{l}
\di\quad \int\int \f{|\pa^{\beta}g^{n+1}|^2}{2\mb{M}_*}dv dx\Big|^{\tau=t}_{\tau=0}
-\int_0^t\int\int \f{\pa^{\beta}g^{n+1}}{\mb{M}_*}\mb{L}_{\bar{\mb{M}}}(\pa^{\beta}g^{n+1})dvdxd\tau\\
\di =\int_0^t\int\int\f{\pa^{\beta}g^{n+1}}{\mb{M}_*}\Big(2\sum_{0<|\alpha|\leq j}C_{\beta}^{\alpha}Q(\pa^{\alpha}\bar{\mb{M}},\pa^{\beta-\alpha}g^{n+1})+2Q(\pa^{\beta}g^n,g^n)\Big)dvdxd\tau\\
\di+\int_0^t\int\int\f{\pa^{\beta}g^{n+1}}{\mb{M}_*}\Big(2\sum_{0<|\alpha|< j}C_{\beta}^{\alpha}Q(\pa^{\alpha}g^n,\pa^{\beta-\alpha}g^{n})
-(\pa^{\beta}\bar{\mb{M}}_t+v\cdot\na_x\pa^{\beta}\bar{\mb{M}})\Big)dvdxd\tau\\
\di :=\sum_{i=3}^4J_i.
\end{array}
\ee}By Cauchy's inequality and Lemma \ref{Lemma 4.1}, one has
\be\label{J3-1}
\ba
J^1_3 & \leq \f{\bar\sigma}{8}\int_0^t\int\int\f{\nu_{\bar{\mb{M}}}(v)}{\mb{M}_*}|\pa^{\beta}g^{n+1}|^2dv dxd\tau\\
&\quad
+C\sum_{0<|\alpha|\leq j}\int_0^t\int\int\f{\nu_{\bar{\mb{M}}}(v)^{-1}}{\mb{M_*}}|Q(\pa^{\alpha}\bar{\mb{M}},\pa^{\beta-\alpha}g^{n+1})|^2dvdxd\tau\\
& \leq \f{\bar\sigma}{8}\int_0^t\int\int\f{\nu_{\bar{\mb{M}}}(v)}{\mb{M}_*}|\pa^{\beta}g^{n+1}|^2dv dxd\tau\\
&\quad
+C\v\sum_{0<|\alpha|\leq j}\int_0^t\int\int\f{\nu_{\bar{\mb{M}}}(v)}{\mb{M}_*}|\pa^{\beta-\alpha}g^{n+1}|^2dv dxd\tau.
\ea
\ee
Similar to \eqref{hard}, we have
\be\label{a-2}
\begin{array}{l}
\di\int\Big(\int\f{\nu_{\bar{\mb{M}}}(v)}{\mb{M}_*}|g^n|^2dv\cdot\int\f{|\pa^{\beta}g^n|^2}{\mb{M}_*}dv\Big)dx\\
\di
\leq \Big\|\int\f{\nu_{\bar{\mb{M}}}(v)}{\mb{M}_*}|g^n|^2dv\Big\|_{L^{\infty}_x}\int\int\f{|\pa^{\beta}g^n|^2}{\mb{M}_*} dv dx\\
\di \leq C\int\int\f{\nu_{\bar{\mb{M}}}(v)}{\mb{M}_*}(|g^n|^2+|\na_xg^n|^2+|\na_x^2g^n|^2)dvdx\cdot\int\int\f{|\pa^{\beta}g^n|^2}{\mb{M}_*} dv dx\\
\di\leq C\Xi^2\int\int\f{\nu_{\bar{\mb{M}}}(v)}{\mb{M}_*}|(g^n,\na_x g^n,\na^2_x g^n)|^2dv dx\leq C\Xi^4.
\end{array}
\ee
It follows from Lemma \ref{Lemma 4.1}, \eqref{a-1} and \eqref{a-2} that
\be\label{J3-2}
\ba
J_3^2 & \leq \f{\bar\sigma}{8}\int_0^t\int\int\f{\nu_{\bar{\mb{M}}}(v)}{\mb{M}_*}|\pa^{\beta}g^{n+1}|^2dv dxd\tau\\
&\quad
+C\int_0^t\int\int\f{\nu_{\bar{\mb{M}}}(v)^{-1}}{\mb{M}_*}|Q(\pa^{\beta}g^n,g^n)|^2dvdxd\tau\\
& \leq \f{\bar\sigma}{8}\int_0^t\int\int\f{\nu_{\bar{\mb{M}}}(v)}{\mb{M}_*}|\pa^{\beta}g^{n+1}|^2dv dxd\tau\\
&\quad
+C\int_0^t\int \Big(\int\f{\nu_{\bar{\mb{M}}}(v)}{\mb{M}_*}|\pa^{\beta}g^n|^2dv\Big)
\cdot\Big(\int\f{|g^n|^2}{\mb{M}_*}dv\Big)dxd\tau\\
& \quad +C\int_0^t\int \Big(\int\f{|\pa^{\beta}g^n|^2}{\mb{M}_*}dv\Big)
\cdot\Big(\int\f{\nu_{\bar{\mb{M}}}(v)}{\mb{M}_*}|g^n|^2dv\Big)dxd\tau\\
& \leq \f{\bar\sigma}{8}\int_0^t\int\int\f{\nu_{\bar{\mb{M}}}(v)}{\mb{M}_*}|\pa^{\beta}g^{n+1}|^2dv dxd\tau+C\Xi^4.
\ea
\ee
Similarly, we have
\be\label{J4-1}
J_4^1\leq \f{\bar\sigma}{8}\int_0^t\int\int\f{\nu_{\bar{\mb{M}}}(v)}{\mb{M}_*}|\pa^{\beta}g^{n+1}|^2dv dxd\tau+C\Xi^4.
\ee
By Cauchy's inequality and Lemma \ref{appu}, we have
\be\label{J4-2}
J_4^2\leq \f{\bar\sigma}{8}\int_0^t\int\int\f{\nu_{\bar{\mb{M}}}(v)}{\mb{M}_*}|\pa^{\beta}g^{n+1}|^2dv dxd\tau+C\v.
\ee

Substituting \eqref{J3-1}, \eqref{J3-2}-\eqref{J4-2} into \eqref{gn+}, we obtain
\be\label{gx-L2}
\begin{array}{l}
\di \int\int \f{|\pa^{\beta}g^{n+1}(t,x,v)|^2}{\mb{M}_*}dv dx
+\int_0^t\int\int\f{\nu_{\bar{\mb{M}}}(v)}{\mb{M}_*}|\pa^{\beta}g^{n+1}|^2dvdxd\tau\\[2mm]
\di \leq C_0\int\int\f{|\pa^{\beta}g_0(x,v)|^2}{\mb{M}_*}dv dx\\[2mm]
\di\quad
+C\v\sum_{0<|\alpha|\leq j}\int_0^t\int\int\f{\nu_{\bar{\mb{M}}}(v)}{\mb{M}_*}|\pa^{\beta-\alpha}g^{n+1}|^2dv dxd\tau
+C(\Xi^4+\v),
\end{array}
\ee
which along with \eqref{g0} and choosing $\v$ suitable small such that $C\v<\f12$ yield
\be\label{g-s}
\begin{array}{l}
\di \|g^{n+1}\|^2_{\mathcal{S}}\leq 2C_0\|g_0\|^2_{\mathcal{S}}+C(\Xi^4+\ln(1+T)+\v)\\[2mm]
\di \qquad\qquad \leq \f12\Xi^2+C(\Xi^4+\ln(1+T)+\v)\leq \Xi^2,
\end{array}
\ee
provided that we choose $\Xi>0$, $T>0$ and $\v>0$ are suitably small, such that
$$
C(\Xi^4+\ln(1+T)+\v)\leq \f12 \Xi^2.
$$
In the following we will show that $\{g^n(t,x,v)\}$ is a Cauchy sequence in $H^3_x\big(L_v^2\break\big(\f{1}{\sqrt{\mb{M}_*}}\big)\big)$. Set
$$
h^n(t,x,v)=g^{n+1}(t,x,v)-g^n(t,x,v),\quad n\geq0,
$$
then $h^n(t,x,v)$ $(n\geq1)$ satisfies
\be\label{h}
\begin{cases}
h^{n}_t+v\cdot\na_x h^{n}=\mb{L}_{\bar{\mb{M}}}h^{n}+Q(g^n,h^{n-1})+Q(g^{n-1},h^{n-1}),\\[3mm]
h^{n}(t,x,v)|_{t=0}=0.
\end{cases}
\ee

First, we multiply $\f{h^n}{\mb{M}_*}$ by the above equation, then integrate by parts over $[0,t]\times\mathbb{D}\times\mathbb{R}^3$ gives
\be\label{h-1}
\begin{array}{l}
\di \quad \int\int\f{|h^n|^2}{2\mb{M}_*}(t,x,v)dvdx-\int_0^t\int\int\f{h^n}{\mb{M}_*}\mb{L}_{\bar{\mb{M}}}h^ndvdxd\tau\\[3mm]
\di =\int_0^t\int\int \f{h^n}{\mb{M}_*}\big(Q(g^n,h^{n-1})+Q(g^{n-1},h^{n-1})\big)dvdxd\tau.
\end{array}
\ee
By Cauchy's inequality and Lemma \ref{Lemma 4.1}, one has
\be
\begin{array}{l}
\di \quad \int_0^t\int\int\f{h^n}{\mb{M}_*}Q(g^n,h^{n-1})dvdxd\tau\\[1mm]
\di \leq \int_0^t\int \Big(\int\f{\nu_{\bar{\mb{M}}}(v)}{\mb{M}_*}|h^n|^2dv\Big)^{\f12}
\Big(\int\f{\nu_{\bar{\mb{M}}}(v)^{-1}}{\mb{M}_*}|Q(g^n,h^{n-1})|^2dv\Big)^{\f12}dxd\tau \\[1mm]
\di \leq  \int_0^t\int \Big(\int\f{\nu_{\bar{\mb{M}}}(v)}{\mb{M}_*}|h^n|^2dv\Big)^{\f12}
\Big(\int\f{|g^n|^2}{\mb{M}_*}dv\cdot\int\f{\nu_{\bar{\mb{M}}}(v)}{\mb{M}_*}|h^{n-1}|^2dv\Big)^{\f12}dxd\tau\\[1mm]
\di \quad +\int_0^t\int \Big(\int\f{\nu_{\bar{\mb{M}}}(v)}{\mb{M}_*}|h^n|^2dv\Big)^{\f12}
\Big(\int\f{\nu_{\bar{\mb{M}}}(v)}{\mb{M}_*}|g^{n}|^2dv\cdot\int\f{|h^{n-1}|^2}{\mb{M}_*}dv\Big)^{\f12}dxd\tau\\[1mm]
\di :=\sum_{i=1}^2K_i.
\end{array}
\ee
It follows from Cauchy's inequality and \eqref{a-1} that
\be\label{K1}
\ba
K_1 & \leq \sup_{0\leq t\leq T}\Big\|\int\f{|g^n|^2}{\mb{M}_*}dv\Big\|^{\f12}_{L^{\infty}_x}
\\&\qquad\cdot
\int_0^t\int\Big(\int \f{\nu_{\bar{\mb{M}}}(v)}{\mb{M_*}}|h^n|^2dv\Big)^{\f12}\cdot
\Big(\int \f{\nu_{\bar{\mb{M}}}(v)}{\mb{M_*}}|h^{n-1}|^2dv\Big)^{\f12} dxd\tau\\[1mm]
& \leq C\Xi\int_0^t\int\int \f{\nu_{\bar{\mb{M}}}(v)}{\mb{M}_*}|(h^n,h^{n-1})|^2dvdxd\tau.
\ea
\ee
By Sobolev inequality, one has
\be\label{K2}
\ba
K_2 & \leq\int_0^t \Big\|\int\f{\nu_{\bar{\mb{M}}}(v)}{\mb{M}_*}|g^n|^2dv\Big\|^{\f12}_{L^{\infty}_x}
\int\Big(\int \f{\nu_{\bar{\mb{M}}}(v)}{\mb{M_*}}|h^n|^2dv\Big)^{\f12}\cdot
\Big(\int \f{|h^{n-1}|^2}{\mb{M_*}}dv\Big)^{\f12} dxd\tau\\
& \leq C\int_0^t\Big(\int\int \f{\nu_{\bar{\mb{M}}}(v)}{\mb{M}_*} |(g^n,\na_x g^n,\na^2_x g^n)|^2dvdx\Big)^{\f12}\\
& \quad \times \Big(\int\int \f{\nu_{\bar{\mb{M}}}(v)}{\mb{M}_*} | h^n|^2dvdx\Big)^{\f12}
\Big(\int\int \f{|h^{n-1}|^2}{\mb{M}_*} dvdx\Big)^{\f12} d\tau\\
& \leq C\sup_{0\leq t\leq T}\Big(\int\int \f{|h^{n-1}|^2}{\mb{M}_*}(t,x,v) dvdx\Big)^{\f12}
\Big(\int_0^t\int\int \f{\nu_{\bar{\mb{M}}}(v)}{\mb{M}_*}|h^n|^2dvdxd\tau\Big)^{\f12}\\
&\quad  \times\Big(\int_0^t\int\int \f{\nu_{\bar{\mb{M}}}(v)}{\mb{M}_*}|(g^n,\na_xg^n,\na_x^2g^n)|^2dvdxd\tau\Big)^{\f12}\\
& \leq C\Xi\sup_{0\leq t\leq T}\int\int\f{|h^{n-1}|^2}{\mb{M}_*}(t,x,v) dvdx
+C\Xi\int_0^t\int\int \f{\nu_{\bar{\mb{M}}}(v)}{\mb{M}_*}|h^n|^2dvdxd\tau.
\ea
\ee
The other term on the right hand side of \eqref{h-1} is the same as estimate \eqref{K1} and \eqref{K2},
so we deduce from \eqref{h-1}, \eqref{K1} and \eqref{K2} that
\be\label{es-h1}
\begin{array}{l}
\di\quad \int\int \f{|h^n|^2}{\mb{M}_*}(t,x,v)dvdx+\int_0^t\int\int \f{\nu_{\bar{\mb{M}}}(v)}{\mb{M}_*}|h^n|^2dvdxd\tau\\
\di \leq C\Xi \Big(\sup_{0\leq t\leq T}\int\int \f{|h^{n-1}|^2}{\mb{M}_*}(t,x,v)dvdx+\int_0^t\int\int \f{\nu_{\bar{\mb{M}}}(v)}{\mb{M}_*}|h^{n-1}|^2dvdxd\tau\Big).
\end{array}
\ee

Next, applying $\pa^{\beta}$ to \eqref{h} with $|\beta|=j~(j=1,2,3)$, then multiplying $\f{\pa^{\beta}h^n}{\mb{M}_*}$ to the resulted equation
and integrating by parts over $[0,t]\times\mathbb{D}\times\mathbb{R}^3$ implies
\be\label{h-2}
\begin{array}{l}
\di\quad
 \int\int\f{|\pa^{\beta}h^n|^2}{2\mb{M}_*}(t,x,v)dvdx-\int_0^t\int\int\f{\pa^{\beta}h^n}{\mb{M}_*}\mb{L}_{\bar{\mb{M}}}(\pa^{\beta}h^n)dvdxd\tau\\
\di =\int_0^t\int\int \f{\pa^{\beta}h^n}{\mb{M}_*}\Big(2\sum_{0<|\alpha|\leq j}C_{\beta}^{\alpha}Q(\pa^{\alpha}\bar{\mb{M}},\pa^{\beta-\alpha}h^n)
+Q(g^n,\pa^{\beta}h^{n-1})\\
\di+Q(g^{n-1},\pa^{\beta}h^{n-1})+2\sum_{0<|\alpha|\leq j}C_{\beta}^{\alpha}Q(\pa^{\alpha}(g^n+g^{n-1}),\pa^{\beta-\alpha}h^{n-1})\Big)dvdxd\tau.
\end{array}
\ee

We just estimate the most difficult term as follows and the other terms can be done similarly.
\be
\begin{array}{l}
\di \quad \int_0^t\int\int\f{\pa^{\beta}h^n}{\mb{M}_*}Q(g^n,\pa^{\beta}h^{n-1})dvdxd\tau\\
\di \leq \int_0^t\int \Big(\int\f{\nu_{\bar{\mb{M}}}(v)}{\mb{M}_*}|\pa^{\beta}h^n|^2dv\Big)^{\f12}
\Big(\int\f{\nu_{\bar{\mb{M}}}(v)^{-1}}{\mb{M}_*}|Q(g^n,\pa^{\beta}h^{n-1})|^2dv\Big)^{\f12}dxd\tau \\
\di \leq  \int_0^t\int \Big(\int\f{\nu_{\bar{\mb{M}}}(v)}{\mb{M}_*}|\pa^{\beta}h^n|^2dv\Big)^{\f12}
\Big(\int\f{|g^n|^2}{\mb{M}_*}dv\cdot\int\f{\nu_{\bar{\mb{M}}}(v)}{\mb{M}_*}|\pa^{\beta}h^{n-1}|^2dv\Big)^{\f12}dxd\tau\\
\di \quad +\int_0^t\int \Big(\int\f{\nu_{\bar{\mb{M}}}(v)}{\mb{M}_*}|\pa^{\beta}h^n|^2dv\Big)^{\f12}
\Big(\int\f{\nu_{\bar{\mb{M}}}(v)}{\mb{M}_*}|g^{n}|^2dv\cdot\int\f{|\pa^{\beta}h^{n-1}|^2}{\mb{M}_*}dv\Big)^{\f12}dxd\tau\\
\di:=\sum_{i=3}^4K_i.
\end{array}
\ee
It follows from Cauchy's inequality and \eqref{a-1} that
\be\label{K3}
\ba
K_3 & \leq \sup_{0\leq t\leq T}\Big\|\int\f{|g^n|^2}{\mb{M}_*}dv\Big\|^{\f12}_{L^{\infty}_x}\\&\qquad\cdot
\int_0^t\int\Big(\int \f{\nu_{\bar{\mb{M}}}(v)}{\mb{M_*}}|\pa^{\beta}h^n|^2dv\Big)^{\f12}\cdot
\Big(\int \f{\nu_{\bar{\mb{M}}}(v)}{\mb{M_*}}|\pa^{\beta}h^{n-1}|^2dv\Big)^{\f12} dxd\tau\\
& \leq C\Xi\int_0^t\int\int \f{\nu_{\bar{\mb{M}}}(v)}{\mb{M}_*}|(\pa^{\beta}h^n,\pa^{\beta}h^{n-1})|^2dvdxd\tau.
\ea
\ee
By Sobolev inequality, one has
\be\label{K4}
\ba
K_4 & \leq\int_0^t \Big\|\int\f{\nu_{\bar{\mb{M}}}(v)}{\mb{M}_*}|g^n|^2dv\Big\|^{\f12}_{L^{\infty}_x}\\
&\qquad\cdot
\int\Big(\int \f{\nu_{\bar{\mb{M}}}(v)}{\mb{M_*}}|\pa^{\beta}h^n|^2dv\Big)^{\f12}\cdot
\Big(\int \f{|\pa^{\beta}h^{n-1}|^2}{\mb{M_*}}dv\Big)^{\f12} dxd\tau\\
& \leq C\int_0^t\Big(\int\int \f{\nu_{\bar{\mb{M}}}(v)}{\mb{M}_*} |(g^n,\na_x g^n,\na^2_x g^n)|^2dvdx\Big)^{\f12}\\
& \quad \times \Big(\int\int \f{\nu_{\bar{\mb{M}}}(v)}{\mb{M}_*} |\pa^{\beta} h^n|^2dvdx\Big)^{\f12}
\Big(\int\int \f{|\pa^{\beta}h^{n-1}|^2}{\mb{M}_*} dvdx\Big)^{\f12} d\tau\\
& \leq C\sup_{0\leq t\leq T}\Big(\int\int \f{|\pa^{\beta}h^{n-1}|^2}{\mb{M}_*}(t,x,v) dvdx\Big)^{\f12}
\Big(\int_0^t\int\int \f{\nu_{\bar{\mb{M}}}(v)}{\mb{M}_*}|\pa^{\beta}h^n|^2dvdxd\tau\Big)^{\f12}
\ea
\ee\be
\ba&\quad  \times\Big(\int_0^t\int\int \f{\nu_{\bar{\mb{M}}}(v)}{\mb{M}_*}|(g^n,\na_xg^n,\na_x^2g^n)|^2dvdxd\tau\Big)^{\f12}\nonumber\\
& \leq C\Xi\sup_{0\leq t\leq T}\int\int\f{|\pa^{\beta}h^{n-1}|^2}{\mb{M}_*}(t,x,v) dvdx
+C\Xi\int_0^t\int\int \f{\nu_{\bar{\mb{M}}}(v)}{\mb{M}_*}|\pa^{\beta}h^n|^2dvdxd\tau.\nonumber
\ea
\ee
It follows from \eqref{h-2}, \eqref{K3} and \eqref{K4} that
\be\label{es-h2}
\begin{array}{l}
\di\quad \int\int \f{|\pa^{\beta}h^n|^2}{\mb{M}_*}(t,x,v)dvdx+\int_0^t\int\int \f{\nu_{\bar{\mb{M}}}(v)}{\mb{M}_*}|\pa^{\beta}h^n|^2dvdxd\tau\\
\di \leq C\v\sum_{0<|\alpha|\leq j}\int_0^t\int\int\f{\nu_{\bar{M}}(|v|)}{\mb{M}_*}|\pa^{\beta-\alpha}h^n|^2dvdxd\tau
+C\Xi \|h^{n-1}\|^2_{\mathcal{S}}.
\end{array}
\ee

Therefore, we deduce from \eqref{es-h1} and \eqref{es-h2} that
\be
\|h^n\|^2_{\mathcal{S}}\leq \f14\|h^{n-1}\|^2_{\mathcal{S}},\quad n\geq 1,
\ee
provided that $\Xi>0$ and $\v>0$ suitable small, thus we have proved that $\{g^n(t,x,v)\}$ is Cauchy sequence in $H^3_x\Big(L_v^2\Big(\f{1}{\sqrt{\mb{M}_*}}\Big)\Big)$.   \hfill $\Box$

%
%
%

%
%
%

\section*{Acknowledgments}
The authors would like to thank the anonymous referees for the valuable comments and suggestions,
which greatly improved the presentation of the manuscript.


\end{document}